\documentclass[aop,preprint]{imsart}
\setattribute{journal}{name}{}

\usepackage[T1]{fontenc}	
\usepackage{graphicx}
\usepackage{amsmath, amsthm, amssymb}
\usepackage{upgreek}
\usepackage{mathrsfs}
\usepackage{pdfsync}
\usepackage[a4paper,top=3cm,bottom=3cm,left=2.5cm,right=2.5cm, bindingoffset=0mm]{geometry}

\usepackage[usenames,dvipsnames]{xcolor}
\usepackage[inline,shortlabels]{enumitem}
\usepackage[hyphens]{url}									
\usepackage{verbatim}								
\usepackage{dsfont}									

\usepackage{mathtools}
\usepackage{xparse}

\newcommand{\texorpdfstring}[2]{#1}

\usepackage{mathtools}

\allowdisplaybreaks

\arxiv{}

\startlocaldefs

\newcommand{\Alpha}{\mathrm{A}}

\usepackage{xparse}

\ExplSyntaxOn
\NewDocumentCommand{\makeabbrev}{mmm}
 {
  \yoruk_makeabbrev:nnn { #1 } { #2 } { #3 }
 }

\cs_new_protected:Npn \yoruk_makeabbrev:nnn #1 #2 #3
 {
  \clist_map_inline:nn { #3 }
   {
    \cs_new_protected:cpn { #2 } { #1 { ##1 } }
   }
 }
 \ExplSyntaxOff

\makeabbrev{\textbf}{tbf#1}{a,b,c,d,e,f,g,h,i,j,k,l,m,n,o,p,q,r,s,t,u,v,w,x,y,z,A,B,C,D,E,F,G,H,I,J,K,L,M,N,O,P,Q,R,S,T,U,V,W,X,Y,Z}

\makeabbrev{\textbf}{bf#1}{a,b,c,d,e,f,g,h,i,j,k,l,m,n,o,p,q,r,s,t,u,v,w,x,y,z,A,B,C,D,E,F,G,H,I,J,K,L,M,N,O,P,Q,R,S,T,U,V,W,X,Y,Z}

\makeabbrev{\textsf}{tsf#1}{a,b,c,d,e,f,g,h,i,j,k,l,m,n,o,p,q,r,s,t,u,v,w,x,y,z,A,B,C,D,E,F,G,H,I,J,K,L,M,N,O,P,Q,R,S,T,U,V,W,X,Y,Z}

\makeabbrev{\mathsf}{mss#1}{a,b,c,d,e,f,g,h,i,j,k,l,m,n,o,p,q,r,s,t,u,v,w,x,y,z,A,B,C,D,E,F,G,H,I,J,K,L,M,N,O,P,Q,R,S,T,U,V,W,X,Y,Z}

\makeabbrev{\mathfrak}{mf#1}{a,b,c,d,e,f,g,h,i,j,k,l,m,n,o,p,q,r,s,t,u,v,w,x,y,z,A,B,C,D,E,F,G,H,I,J,K,L,M,N,O,P,Q,R,S,T,U,V,W,X,Y,Z}

\makeabbrev{\mathrm}{mrm#1}{a,b,c,d,e,f,g,h,i,j,k,l,m,n,o,p,q,r,s,t,u,v,w,x,y,z,A,B,C,D,E,F,G,H,I,J,K,L,M,N,O,P,Q,R,S,T,U,V,W,X,Y,Z}

\makeabbrev{\mathbf}{mbf#1}{a,b,c,d,e,f,g,h,i,j,k,l,m,n,o,p,q,r,s,t,u,v,w,x,y,z,A,B,C,D,E,F,G,H,I,J,K,L,M,N,O,P,Q,R,S,T,U,V,W,X,Y,Z}

\makeabbrev{\mathcal}{mc#1}{A,B,C,D,E,F,G,H,I,J,K,L,M,N,O,P,Q,R,S,T,U,V,W,X,Y,Z}

\makeabbrev{\mathbb}{mbb#1}{A,B,C,D,E,F,G,H,I,J,K,L,M,N,O,P,Q,R,S,T,U,V,W,X,Y,Z}

\makeabbrev{\mathscr}{ms#1}{A,B,C,D,E,F,G,H,I,J,K,L,M,N,O,P,Q,R,S,T,U,V,W,X,Y,Z}

\makeabbrev{\mathrm}{#1}{
Id,id,ran,rk,diag,stab,ann,conv,pr,ev,tr,End,Hom,sgn,im,op,can,fin,ext,red,tot,Leb,
rot,usc,lsc,Lip,lip,bSymLip,osc,AC,loc,coz,z,
supp,Opt,Adm,Cpl,Geo,GeoOpt,GeoAdm,GeoCpl,reg,res,
bd,co,Ric,Exp,dExp,dist,seg,Seg,cut,fcut,Cut,SDiff,Iso,Isom,diam,cl,Homeo,Diff,Der,vol,dvol,inj,relint, Graph, sub,
var,law,Var,Poi,Gam,pa,so,iso,fs,inv,pqi,mix,erg,
TestF,
}

\makeabbrev{\mathsf}{#1}{CD,BE,MCP,Ent,wMTW,MTW,Ch,RCD,EVI,Rad,dRad,SL,cSL,dSL,ScL,Irr,SC,wFe,VA}

\makeabbrev{\mathsc}{#1}{mmaf,cg}

\newcommand{\T}{\tau} 
\newcommand{\Bo}[1]{\msB(#1)} 

\newcommand{\spec}[1]{\mathrm{spec}(#1)}

\newcommand{\defeq}{\eqqcolon}
\renewcommand{\complement}{\mathrm{c}}
\newcommand{\acts}{\,\raisebox{\depth}{\scalebox{1}[-1]{$\circlearrowleft$}}\,}
\newcommand{\mathsc}[1]{\text{\textsc{#1}}}
\newcommand{\emparg}{{\,\cdot\,}}
\newcommand{\dint}[2][]{\;\sideset{^{\scriptstyle{#1}}\!\!\!\!}{_{#2}^{\scriptscriptstyle\oplus}}\int}
\newcommand{\forallae}[1]{{\textrm{\,for ${#1}$-a.e.~}}}
\newcommand{\as}[1]{\quad #1\text{-a.e.}}
\renewcommand{\Cap}{\mathrm{cap}}
\newcommand{\dom}[1]{\msD(#1)}
\newcommand{\domext}[1]{\msD(#1)_e}

\DeclareMathOperator{\eqdef}{\coloneqq}

\let\epsilon\varepsilon

\newcommand{\longrar}{\longrightarrow}
\newcommand{\rar}{\rightarrow}

\newcommand{\nlim}{\lim_{n}}								

\newcommand{\diff}{\mathop{}\!\mathrm{d}}						
		
\newcommand{\abs}[1]{\left\lvert#1\right\rvert}						
\newcommand{\tnorm}[1]{\big\lVert#1\big\rVert}					
\newcommand{\norm}[1]{\left\lVert#1\right\rVert}					
\newcommand{\ttnorm}[1]{\lVert#1\rVert}					
\newcommand{\set}[1]{\left\{#1\right\}}							
\newcommand{\ttset}[1]{\{#1\}}									
\newcommand{\tonde}[1]{\left(#1\right)}							
\newcommand{\ttonde}[1]{\big({#1}\big)}
\newcommand{\class}[2][]{\left[#2\right]_{#1}}						
\newcommand{\tclass}[2][]{\big [#2\big]_{#1}}						
\newcommand{\ttclass}[2][]{[#2]_{#1}}						

\newcommand{\rep}[1]{\hat{#1}}								
\newcommand{\reptwo}[1]{\tilde{#1}}							
\newcommand{\tquadre}[1]{\big[#1\big]}							
\newcommand{\tscalar}[2]{\big\langle #1 \, \big |\, #2\big\rangle}			
\newcommand{\scalar}[2]{\left\langle #1 \,\middle |\, #2\right\rangle}		

\newcommand{\otym}[1]{{\scriptscriptstyle{\otimes #1}}}
\newcommand{\sym}[1]{{\scriptscriptstyle{(#1)}}}
	
\newcommand{\seq}[1]{\tonde{#1}}								
\newcommand{\tseq}[1]{{\big(#1\big)}}
\newcommand{\ttseq}[1]{(#1)}
\newcommand{\Cc}{\mcC_c}									
\newcommand{\Cz}{\mcC_0}									

\newcommand{\pfwd}{\sharp}
\DeclareMathOperator*{\esssup}{esssup}

\DeclareMathOperator{\car}{\mathds 1}

\DeclareMathOperator{\emp}{\varnothing}
\newcommand{\N}{{\mathbb N}}
\newcommand{\R}{{\mathbb R}}
\DeclareMathOperator{\Q}{{\mathbb Q}}
\newcommand{\C}{{\mathbb C}} 	

\newcommand{\restr}{\big\lvert}

\newcommand{\iref}[1]{\ref{#1}}

\newcommand{\comma}{\,\,\mathrm{,}\;\,}
\newcommand{\semicolon}{\,\,\mathrm{;}\;\,}
\newcommand{\fstop}{\,\,\mathrm{.}}

\newcommand{\cdc}{\Gamma}
\DeclareMathOperator{\zero}{{\mathbf 0}}
\DeclareMathOperator{\uno}{{\mathbf 1}}

\usepackage{scrextend}						

\newcommand{\n}[1]{{\overline{#1}}}

\let\temp\phi
\let\phi\varphi
\let\varphi\temp

\newcommand{\PP}{{\pi}}

\numberwithin{equation}{section}
\theoremstyle{plain}
\newtheorem{theorem}{Theorem}[section]
\newtheorem*{theorem*}{Theorem}
\newtheorem{proposition}[theorem]{Proposition}
\newtheorem{lemma}[theorem]{Lemma}
\newtheorem{corollary}[theorem]{Corollary}

\theoremstyle{definition}
\newtheorem{definition}[theorem]{Definition}

\theoremstyle{remark}
\newtheorem{remark}[theorem]{Remark}
\newtheorem{example}[theorem]{Example}
\newtheorem{assumption}[theorem]{Assumption}

\endlocaldefs

\begin{document}

\begin{frontmatter}
\title{Ergodic Decomposition of Dirichlet Forms\\via Direct Integrals and Applications\thanksref{T1}
}
\runtitle{Ergodic Decomposition of Dirichlet Forms}

\begin{aug}
\author{\fnms{Lorenzo} \snm{Dello Schiavo}\thanksref{t2}\ead[label=e1]{lorenzo.delloschiavo@ist.ac.at}}

\thankstext{t2}{The author is grateful to Professors Sergio Albeverio and Andreas Eberle, and to Dr.~Kohei Suzuki, for fruitful conversations on the subject of the present work, and for respectively pointing out the references~\cite{AlbKonRoe97},~\cite{Ebe96}, and~\cite{Kuw11, AlbRoe90}.
Finally, he is especially grateful to an anonymous Reviewer for their very careful reading and their suggestions which improved the readability of the paper.
}

\thankstext{T1}{A large part of this work was written while the author was a member of the Institut f\"ur angewandte Mathematik of the University of Bonn. Research supported by the Collaborative Research Center 1060.
The author gratefully acknowledges funding of his current position by the Austrian Science Fund (FWF) grant F65, and by the European Research Council (ERC, grant No.~716117, awarded to Prof.\ Dr.~Jan Maas).}

\runauthor{L.~Dello~Schiavo}

\affiliation{IST Austria}

\address{Lorenzo~~Dello Schiavo
\\
IST Austria,\\
Am Campus 1,\\
3400 Klosterneuburg,\\
Austria\\
\printead{e1}
}

\end{aug}

\begin{abstract}
We study superpositions and direct integrals of quadratic and Dirichlet forms.
We show that each quasi-regular Dirichlet space over a probability space admits a unique representation as a direct integral of irreducible Dirichlet spaces, quasi-regular for the same underlying topology.
The same holds for each quasi-regular strongly local Dirichlet space over a metrizable Luzin, Radon measure space, and admitting carr\'e du champ operator. In this case, the representation is only projectively unique.
\end{abstract}

\vspace{.5cm}
\today
\vspace{.5cm}

\begin{keyword}[class=MSC]
\kwd{37A30} \kwd{secondary: 31C25, 60J25, 60J35}
\end{keyword}

\begin{keyword} \kwd{direct integral} \kwd{ergodic decomposition} \kwd{ergodic theorem} \kwd{Dirichlet forms.}
\end{keyword}

\end{frontmatter}


\section{Introduction}\label{s:Intro}
Let~$(X,\T,\mcX,\mu)$ be a locally compact Polish Radon measure space, and~$(E,\dom{E})$ be a regular Dirichlet form on~$(X,\T,\mcX,\mu)$. As it is well-known,~$(E,\dom{E})$ is properly associated with a right process
\begin{align*}
\mbfM\eqdef \ttonde{\Omega, \mcF, \seq{M_t}_{t\geq 0}, \seq{P_x}_{x\in X_\partial},\xi}
\end{align*}
with state space~$X$, life-time~$\xi$, and cemetery~$\partial$.
For a $\mu$-measurable subset~$A\subset X$, we say that
\begin{itemize}
\item $A$ is \emph{$\mbfM$-invariant} (e.g.~\cite[Dfn.~IV.6.1]{MaRoe92}) if there exists~$\Omega_{A^\complement}\in\mcF$ with~$P_x\Omega_{A^\complement}=0$ for every~$x\in A$~and
\begin{align*}
\Omega_{A^\complement} \supset \set{\omega\in \Omega: A^\complement \cap \overline{\set{M_s(\omega): s\in [0,t]}} \neq \emp \text{~for some~} 0\leq t<\xi} \semicolon
\end{align*}

\item $A$ is \emph{$E$-invariant} (also cf.~Dfn.~\ref{d:Invariant} below) if~$\car_A f\in \dom{E}$ for any~$f\in \dom{E}$ and
\begin{align*}
E(f,g)=E(\car_A f,\car_A g)+E(\car_{A^\complement}f,\car_{A^\complement} g)\comma \qquad f,g\in \dom{E}
\end{align*}
\end{itemize}

If the form~$(E,\dom{E})$ is additionally strongly local, then the process~$\mbfM$ is a Markov diffusion, and the following are $\mu$-essentially equivalent, see Rmk.~\ref{r:QuasiClopen}~below:
\begin{itemize}
\item $A$ is {$\mbfM$-invariant};
\item $A$ is {$E$-invariant};
\item $A$ is $E$-\emph{quasi-clopen}, i.e., simultaneously $E$-quasi-open and $E$-quasi-closed (see e.g.,~\cite[p.~70]{FukOshTak11}).
\end{itemize}

We say that a set~$A\subset X$ is \emph{$\mu$-trivial} if it is $\mu$-measurable and either~$\mu A=0$ or~$\mu A^\complement=0$. The process~$\mbfM$ is \emph{irreducible} if every $\mbfM$-invariant set is $\mu$-trivial.
When~$\mbfM$ is not irreducible, it is natural ---~in the study of the pathwise properties of~$\mbfM$~--- to restrict our attention to ``minimal'' $\mbfM$-invariant subsets of~$X$. In the local case, thanks to the quasi-topological characterization of $\mbfM$-invariance, such sets may be thought of as the ``connected components'' of the space~$X$ as seen by~$\mbfM$.

This description is in fact purely heuristic, since it may happen that all such ``minimal'' $\mbfM$-invariant sets are $\mu$-negligible. The question arises, whether these ideas for the study of $\mbfM$-invariance can be made rigorous by resorting to the Dirichlet form~$(E,\dom{E})$ associated with~$\mbfM$.
More precisely, we look for a decomposition~$\ttonde{(E_\zeta,\dom{E_\zeta})}_{\zeta\in Z}$ of~$(E,\dom{E})$ over some index set~$Z$, and we require that
\begin{itemize}
\item $(E_\zeta,\dom{E_\zeta})$ is a Dirichlet form on~$(X,\T)$ additionally irreducible (Dfn.~\ref{d:Invariant}) for every~$\zeta\in Z$;
\item we may reconstruct~$(E,\dom{E})$ from~$\ttonde{(E_\zeta,\dom{E_\zeta})}_{\zeta\in Z}$ in a unique way.
\end{itemize}
Because of the first property, such a decomposition ---~if any~--- would deserve the name of \emph{ergodic decomposition} of~$(E,\dom{E})$.

For instance, let us consider the standard Dirichlet form~$E^g$ on a (second-countable) Riemannian manifold~$(M,g)$, i.e.\ the one generated by the (negative halved) Laplace--Beltrami operator~$\Delta_g$ and properly associated with the Brownian motion on~$M$. In this case, we expect that~$Z$ is a discrete space, indexing the connected components of~$M$, and that
\begin{align*}
E^g=\bigoplus_{\zeta\in Z} E^g_\zeta\comma
\end{align*}
where~$(E^g_\zeta,\dom{E^g_\zeta})$ is but the standard form of the connected component of index~$\zeta$.
This simple example suggests that, in the general case of our interest, we should expect that~$(E,\dom{E})$ is recovered from the decomposition~$\ttonde{(E_\zeta,\dom{E_\zeta})}_{\zeta\in Z}$ as a ``direct integral'',
\begin{align*}
E=\dint[]{Z} E_\zeta \fstop
\end{align*}

Our purpose is morefold:
\begin{itemize}
\item to introduce a notion of \emph{direct integral of Dirichlet forms}, and to compare it with the existing notions of \emph{superposition of Dirichlet forms}~\cite[\S{V.3.1}]{BouHir91} (also cf.~\cite[\S3.1($2^\circ$), p.~113]{FukOshTak11} and~\cite{Tom80}), and of \emph{direct integral of quadratic forms}~\cite{AlbRoe90};

\item to discuss an Ergodic Decomposition Theorem for quasi-regular Dirichlet forms, a counterpart for Dirichlet forms to the Ergodic Decomposition Theorems for group actions, e.g.~\cite{GreSch00,Buf14, deJRoz17};

\item to provide rigorous justification to the assumption ---~quite standard in the literature about \mbox{(quasi-)} regular Dirichlet forms~---, that one may consider irreducible forms with no loss of generality;

\item to establish tools for the generalization to arbitrary (quasi-regular) Dirichlet spaces of results currently available only in the irreducible case, e.g.\ the study~\cite{LenSchWir18}
 of invariance under order-isomorphism, cf.~\cite{LzDSWir21}.
\end{itemize}

For strongly local Dirichlet forms, our ergodic decomposition result takes the following form.
\begin{theorem*}
Let~$(X,\T,\mcX,\mu)$ be a metrizable Luzin Radon measure space, and consider a quasi-regular strongly local Dirichlet form $(E,\dom{E})$ on $L^2(\mu)$ admitting carr\'e du champ operator. Then, there exist
\begin{enumerate}[$(i)$]
\item a $\sigma$-finite measure space~$(Z,\mcZ,\nu)$;
\item a family of measures~$\seq{\mu_\zeta}_{\zeta\in Z}$ so that~$(X,\T,\mcX,\mu_\zeta)$ is a (metrizable Luzin) Radon measure space for~$\nu$-a.e.~$\zeta\in Z$, the map~$\zeta\mapsto \mu_\zeta A$ is $\nu$-measurable for every~$A\in\mcX$ and
\begin{align*}
\mu A=\int_Z \mu_\zeta A \diff\nu(\zeta)\comma \qquad A\in\mcX\comma
\end{align*}
and for~$\nu^{\otimes 2}$-a.e.~$(\zeta,\zeta')$, with~$\zeta\neq\zeta'$, the measures~$\mu_\zeta$ and~$\mu_{\zeta'}$ are mutually singular;
\item a family of quasi-regular strongly local irreducible Dirichlet forms~$(E_\zeta,\dom{E_\zeta})$ on~$L^2(\mu_\zeta)$;
\end{enumerate}
so that
\begin{align*}
L^2(\mu)=\dint[]{Z} L^2(\mu_\zeta)\diff\nu(\zeta) \qquad \text{and} \qquad E=\dint[]{Z} E_\zeta \diff\nu(\zeta) \fstop
\end{align*}

Additionally, the disintegration is ($\nu$-essentially) projectively unique, and unique if~$\mu$ is totally finite.
\end{theorem*}

\paragraph{Plan of the work} Firstly, we shall discuss the notion of direct integral of (non-negative definite) quadratic forms on abstract Hilbert spaces,~\S\ref{ss:DirIntQ}, and specialize it to direct integrals of Dirichlet forms,~\S\ref{ss:DirIntE}, via disintegration of measures. In~\S\ref{ss:ErgProbab} we show existence and uniqueness of the ergodic decomposition for regular and quasi-regular (not necessarily local) Dirichlet forms on probability spaces. The results are subsequently extended to strongly local quasi-regular Dirichlet forms on $\sigma$-finite spaces and admitting carr\'e du champ operator,~\S\ref{ss:ErgSigmaFinite}.
Examples are discussed in~\S\ref{ss:Examples}; an application is discussed in~\ref{ss:Appl}.

\paragraph{Bibliographical note} Our reference of choice for direct integrals of Hilbert spaces is the monograph~\cite{Dix81} by J.~Dixmier. For some results we shall however need the more general concept of direct integrals of Banach spaces, after~\cite{HayMirYve91, deJRoz17}. For the sake of simplicity, all such results are confined to the Appendix.

\section{Direct Integrals}\label{s:DirInt}
Every Hilbert space is assumed to be \emph{separable} and a real vector space.

\subsection{Quadratic forms} Let~$(H,\norm{\emparg})$ be a Hilbert space with scalar product~$\scalar{\emparg}{\emparg}$.
%
By a \emph{quadratic form}~$(Q,D)$ on~$H$ we shall always mean a symmetric positive semi-definite ---~if not otherwise stated, densely defined~--- bilinear form. To~$(Q,D)$ we associate the non-relabeled quadratic functional~$Q\colon H\rar \R\cup\set{+\infty}$ defined by
\begin{align*}
Q(u)\eqdef \begin{cases}
Q(u,u) & \text{if } u\in D
\\
+\infty & \text{otherwise}
\end{cases} \comma \qquad u\in H\fstop
\end{align*}

Additionally, for every~$\alpha>0$, we set
\begin{align*}
Q_\alpha(u,v)\eqdef&\ Q(u,v)+\alpha\scalar{u}{v}\comma \qquad u,v\in D\comma
\\
Q_\alpha(u)\eqdef&\ Q(u)+\alpha\norm{u}^2\comma \qquad\qquad u\in H \fstop
\end{align*}
For~$\alpha>0$, we let~$\dom{Q}_\alpha$ be the completion of~$D$, endowed with the Hilbert norm~$Q_\alpha^{1/2}$.
The following result is well-known.

\begin{lemma}\label{l:WLSC} Let~$(Q, D)$ be a quadratic form on~$H$. The following are equivalent:
\begin{enumerate}[$(i)$]
\item\label{i:l:WLSC:1} $(Q,D)$ is closable, say, with closure~$(Q,\dom{Q})$;
\item\label{i:l:WLSC:2} the identical inclusion~$\iota\colon D\rar H$ extends to a continuous injection~$\iota_\alpha\colon \dom{Q}_\alpha\rar H$ satisfying~$\norm{\iota_\alpha}_\op\leq \alpha^{-1}$;
\item\label{i:l:WLSC:3} $Q$ is lower semi-continuous w.r.t. the strong topology of~$H$;
\item\label{i:l:WLSC:4} $Q$ is lower semi-continuous w.r.t. the weak topology of~$H$.
\end{enumerate}
\end{lemma}
\begin{proof}
\iref{i:l:WLSC:1}$\iff$\iref{i:l:WLSC:3} is~\cite[Lem.~VIII.3.14a, p.~461]{Kat95}. \iref{i:l:WLSC:1}$\iff$\iref{i:l:WLSC:2} is noted in~\cite[Rmk.~I.3.2.(ii)]{MaRoe92}. For~\iref{i:l:WLSC:3}$\iff$\iref{i:l:WLSC:4} note that every convex subset of a Hilbert space is weakly closed if and only if it is strongly closed and apply this fact to the sublevel sets of~$Q\colon H\rar \R\cup\set{+\infty}$.
\end{proof}

To every closed quadratic form~$(Q,\dom{Q})$ we associate a non-negative self-adjoint operator~$-L$, with domain defined by the equality~$\dom{\sqrt{-L}}=\dom{Q}$, such that~$Q(u,v) = \scalar{-Lu}{v}$ for all~$u,v\in\dom{L}$. We denote the associated strongly continuous contraction semigroup by~$T_t\eqdef e^{tL}$,~$t>0$, and the associated strongly continuous contraction resolvent by~$G_\alpha\eqdef (\alpha-L)^{-1}$,~$\alpha> 0$. By Hille--Yosida Theorem, e.g.~\cite[p.~27]{MaRoe92},
\begin{subequations}\label{eq:Hille--Yosida}
\begin{align}
\label{eq:Hille--YosidaG}
Q_\alpha(G_\alpha u,v)=&\scalar{u}{v}_H\comma\qquad u\in H\comma v\in\dom{Q}\comma 
\\
\label{eq:Hille--YosidaT}
T_t=&H\text{-}\!\!\lim_{\alpha\rar \infty} e^{t\alpha (\alpha G_\alpha-1)} \fstop
\\
\label{eq:Hille--YosidaE}
Q(u,v)=\lim_{\beta\rar\infty}Q^{(\beta)}(u,v)\eqdef& \lim_{\beta\rar\infty} \scalar{\beta u-\beta G_\beta u}{v}_H\comma \qquad u,v\in H
\fstop
\end{align}
\end{subequations}

\subsection{Direct integrals}
Let~$\seq{H_\zeta}_{\zeta\in Z}$ be a family of Hilbert spaces indexed by some index set~$Z$.
If~$Z$ is at most countable, the \emph{direct sum} of the Hilbert spaces~$H_\zeta$ is defined as
\begin{align}\label{eq:DirectSum}
\bigoplus_{\zeta\in Z} H_\zeta\eqdef \set{\seq{h_\zeta}_{\zeta\in Z} : h_\zeta\in H_\zeta \text{ for all~$\zeta\in Z$ and } \sum_{\zeta\in Z} \norm{h_\zeta}_{H_\zeta}^2<\infty} \fstop
\end{align}

The \emph{direct integral} of a family of Hilbert spaces is a natural generalization of the concept of direct sum of Hilbert spaces to the case when the indexing set~$Z$ is more than countable.
In this case, the requirement of $\ell_2$-summability in the definition of direct sum is replaced by a requirement of $L^2$-integrability (see below for the precise definitions), which implies that~$Z$ should be taken to be a measure space~$(Z,\mcZ,\nu)$.
When~$Z$ is an at most countable discrete space, and~$\nu$ is the counting measure, then the direct integral of the Hilbert spaces~$\seq{H_\zeta}_\zeta$ is isomorphic to their direct sum~\eqref{eq:DirectSum}.
Since their introduction by J.~von Neumann in~\cite[\S3]{Neu49}, direct integrals have become a main tool in operator theory, and in particular in the classification of von Neumann algebras.

In order to make the definition of direct integral precise, let us first introduce some measure-theoretical notions.

\begin{definition}[Measure spaces]\label{d:Standard} A measurable space~$(X,\mcX)$ is
\begin{itemize}
\item \emph{separable} if~$\mcX$ contains all singletons in~$X$, i.e.~$\set{x}\in\mcX$ for each~$x\in X$;

\item \emph{separated} if~$\mcX$ separates points in~$X$, i.e.\ for every~$x,y\in X$ there exists~$A\in\mcX$ with~$x\in A$ and~$y\in A^\complement$;

\item \emph{countably separated} if there exists a countable family of sets in~$\mcX$ separating points in~$X$;

\item \emph{countably generated} if there exists a countable family of sets in~$\mcX$ generating~$\mcX$ as a $\sigma$-algebra;

\item a \emph{standard Borel space} if there exists a Polish topology~$\T$ on~$X$ so that~$\mcX$ coincides with the Borel $\sigma$-algebra induced by~$\T$.
\end{itemize}

For any subset~$X_0$ of a measurable space~$(X,\mcX)$, the trace $\sigma$-algebra on~$X_0$ is $\mcX\cap X_0\eqdef\set{A\cap X_0: A\in\mcX}$.
A $\sigma$-finite measure space~$(X,\mcX,\mu)$ is \emph{standard} if there exists $X_0\in\mcX$, $\mu$-conegligible and so that~$(X_0,\mcX\cap X_0)$ is a standard Borel space. 
We denote by~$(X,\mcX^\mu,\hat\mu)$ the (Carath\'eodory) completion of~$(X,\mcX,\mu)$. A $[-\infty,\infty]$-valued function is called \emph{$\mu$-measurable} if it is measurable w.r.t.~$\mcX^\mu$. For measures~$\mu_1$,~$\mu_2$ we write~$\mu_1\sim \mu_2$ to indicate that~$\mu_1$ and~$\mu_2$ are equivalent, i.e.\ mutually absolutely continuous.
A \emph{$\sigma$-ideal}~$\mcN$ of a measure space~$(X,\mcX,\mu)$ is any sub-$\sigma$-algebra of~$\mcX$ closed w.r.t.~$\cap$, i.e.\ satisfying~$A\cap N\in \mcN$ whenever~$A\in\mcX$ and~$N\in\mcN$.
In particular, the family~$\mcN_\mu$ of $\mu$-negligible subsets of~$(X,\mcX)$ is always a $\sigma$-ideal of~$(X,\mcX^\mu,\hat\mu)$.

For functions~$f,g\colon X\to\R\cup\set{\pm\infty}$ we denote by~$f^+\eqdef 0\vee f$, resp.~$f^-\eqdef -(0\wedge f)$, the positive, resp.\ negative part of~$f$, and by~$f\wedge g$, resp.~$f\vee g$, the pointwise minimum, resp.\ maximum, of~$f$ and~$g$.
\end{definition}

We now recall the main definitions concerning direct integrals of separable Hilbert spaces, referring to~\cite[\S\S{II.1}, II.2]{Dix81} for a systematic treatment.

\begin{definition}[Measurable fields,~{\cite[\S{II.1.3}, Dfn.~1,~p.~164]{Dix81}}]\label{d:DirInt}
\
Let~$(Z,\mcZ,\nu)$ be a $\sigma$-finite measure space,~$\seq{H_\zeta}_{\zeta\in Z}$ be a family of separable Hilbert spaces, and~$F$ be the linear space~$F\eqdef \prod_{\zeta\in Z} H_\zeta$. We say that~$\zeta\mapsto H_\zeta$ is a \emph{$\nu$-measurable field of Hilbert spaces} (\emph{with underlying space~$S$}) if there exists a linear subspace~$S$ of~$F$ with
\begin{enumerate}[$(a)$]
\item\label{i:d:DirInt1} for every~$u\in S$, the function~$\zeta\mapsto \norm{u_\zeta}_{\zeta}$ is $\nu$-measurable;
\item\label{i:d:DirInt2} if~$v\in F$ is such that $\zeta\mapsto \scalar{u_\zeta}{v_\zeta}_{\zeta}$ is $\nu$-measurable for every~$u\in S$, then~$v\in S$;
\item\label{i:d:DirInt3} there exists a sequence~$\seq{u_n}_n\subset S$ such that~$\seq{u_{n,\zeta}}_n$ is a total sequence\footnote{A sequence~$\seq{u_n}_n$ in a Banach space~$B$ is called \emph{total} if the strong closure of its linear span coincides with~$B$.} in~$H_\zeta$ for every~$\zeta\in Z$.
\end{enumerate}

Any such $S$ is called a \emph{space of $\nu$-measurable vector fields}. Any sequence in~$S$ possessing property~\iref{i:d:DirInt3} is called a \emph{fundamental} sequence.
\end{definition}

\begin{proposition}[{\cite[\S{II.1.4, Prop.~4, p.~167}]{Dix81}}]\label{p:Dix} Let~$\mcS$ be a subfamily of~$F$ satisfying both Definition~\ref{d:DirInt}\iref{i:d:DirInt1} and~\iref{i:d:DirInt3} with~$\mcS$ in place of~$S$. Then, there exists exactly one space of $\nu$-measurable vector fields~$S$ so that~$\mcS\subset S$.
\end{proposition}

\begin{definition}[Direct integrals, {\cite[\S{II.1.5}, Prop.~5,~p.~169]{Dix81}}]\label{d:DirIntH}
Let~$\zeta\mapsto H_\zeta$ be a $\nu$-measurable field of Hilbert spaces with underlying space~$S$.
A $\nu$-measurable vector field~$u\in S$ is called ($\nu$-)\emph{square-integrable} if
\begin{align}\label{eq:NormH}
\norm{u}\eqdef \tonde{\int_Z \norm{u_\zeta}_{\zeta}^2 \, \diff\nu(\zeta)}^{1/2}<\infty\fstop
\end{align}

Two square-integrable vector fields~$u$,~$v$ are called ($\nu$-)\emph{equivalent} if~$\norm{u-v}=0$. The space~$H$ of equivalence classes of square-integrable vector fields, endowed with the non-relabeled quotient norm~$\norm{\emparg}$, is a Hilbert space~\cite[\S{II.1.5}, Prop.~5(i), p.~169]{Dix81}, called the \emph{direct integral of~$\zeta\mapsto H_\zeta$} (\emph{with underlying space~$S$}) and denoted by
\begin{align}\label{eq:DirInt}
H=\dint[S]{Z} H_\zeta \diff\nu(\zeta) \fstop
\end{align}
The superscript~`$S$' is omitted whenever~$S$ is apparent from context.
\end{definition}

In the following, it will occasionally be necessary to distinguish an element~$u$ of~$H$ from one of its representatives modulo $\nu$-equivalence, say~$\rep u$ in~$S$.
In this case, we shall write~$u=\class[H]{\rep u}$.
When the specification of the variable~$\zeta$ is necessary, given~$u\in H$, resp.~$\rep u\in S$, we shall write~$\zeta\mapsto u_\zeta$ in place of~$u$, resp.~$\zeta\mapsto \rep u_\zeta$ in place of~$\rep u$.
In most cases, the distinction between~$u$ and~$\rep u$ is however immaterial, similarly to the distinction between the class of a function in $L^2(\nu)$ and any of its $\nu$-representatives.
Therefore in most cases we shall simply write~$u$ in place of both~$u$ and~$\rep u$.

\begin{lemma} Let~$(Z,\mcZ,\nu)$ be a $\sigma$-finite countably generated measure space. Then, the space~$H$ in~\eqref{eq:DirInt} is separable.
\end{lemma}
\begin{proof}
It suffices to note that~$L^2(\nu)$ is separable, e.g.~\cite[365X(p)]{Fre00}. Then, the proof of~\cite[\S{II.1.6}, Cor., p.~172]{Dix81} applies \emph{verbatim}.
\end{proof}

\begin{remark}\label{r:S0} In general, the space~$H$ in~\eqref{eq:DirInt} depends on~$S$, cf.~\cite[p.~169, after Dfn.~3]{Dix81}. It is nowadays standard to define \emph{the} direct integral of~$\zeta\mapsto H_\zeta$ as the one with underlying space~$S$ generated (in the sense of Proposition~\ref{p:Dix}) by an algebra~$\mcS$ of `simple functions', see e.g.~\cite[\S6.1, p.~61]{HayMirYve91} or the Appendix. Here, we prefer the original definition in~\cite{Dix81}, since we shall make a (possibly) different choice of~$S$, more natural when addressing direct integrals of Dirichlet forms.
\end{remark}

Let~$H$ be a direct integral Hilbert space defined as in~\eqref{eq:DirInt}.
We now turn to the discussion of bounded operators in~$\mcB(H)$ and their representation by measurable fields of bounded operators.

\begin{definition}[Measurable fields of bounded operators, decomposable operators]
A field of bounded operators~$\zeta\mapsto B_\zeta\in\mcB(H_\zeta)$ is called \emph{$\nu$-measurable} (\emph{with underlying space~$S$}) if~$\zeta\mapsto B_\zeta u_\zeta\in H_\zeta$ is a $\nu$-measurable vector field with underlying space~$S$ for every $\nu$-measurable vector field~$u$ with underlying space~$S$.
A $\nu$-measurable vector field of bounded operators is called \emph{$\nu$-essentially bounded} if~$\nu$-$\esssup_{\zeta\in Z}\norm{B_\zeta}_{\op,\zeta}<\infty$.
In this case, the direct integral operator~$B\colon H\to H$ of~$\zeta\mapsto B_\zeta$ given by
\begin{align}\label{eq:BoundedOp}
B\colon \class[H]{\rep u}\longmapsto \class[H]{\zeta\mapsto B_\zeta \rep u_\zeta}
\end{align}
is well-defined (in the sense that it does not depend on the choice of the representative~$\rep u\in S$ of~$\class[H]{\rep u}\in H$), and a bounded operator in~$\mcB(H)$.
Its operator norm~$\norm{B}_\op$ satisfies~$\norm{B}_\op=\nu\text{-}\esssup_{\zeta\in Z}\norm{B_\zeta}_{\op,\zeta}$,~\cite[\S{II.2.3},~Prop.~2,~p.~181]{Dix81}.
Conversely, a bounded operator~$B\in\mcB(H)$ is called \emph{decomposable},~\cite[\S{II.2.3}, Dfn.~2, p.~182]{Dix81}, if~$B$ is represented by a $\nu$-essentially bounded $\nu$-measurable field of bounded operators~$\zeta\mapsto B_\zeta$ in the sense of~\eqref{eq:BoundedOp}, in which case we write
\begin{align*}
B=\dint{Z} B_\zeta \diff\nu(\zeta) \fstop
\end{align*}
\end{definition}

The next statement is readily deduced from e.g.~\cite[Thm.~2]{ChoGil71} or~\cite[Thm.~1.10]{Lan75}. For the reader's convenience, a short proof is included.

\begin{lemma}\label{l:DecomposabilityCFC} Let~$H$ be defined as in~\eqref{eq:DirInt},~$B\in \mcB(H)$ be decomposable, and~$D_B$ be the closed disk of radius~$\norm{B}_\op$ in the complex plane. Then, for every~$\phi\in \mcC(D_B)$ the continuous functional calculus~$\phi(B)$ of~$B$ is decomposable and
\begin{align}\label{eq:CFC}
\phi(B)=\dint{Z} \phi(B_\zeta) \diff\nu(\zeta) \fstop
\end{align}
\end{lemma}
\begin{proof} Well-posedness follows by~\cite[\S{II.2.3},~Prop.~2,~p.~181]{Dix81}. The proof is then a
straightforward application of~\cite[\S{II.2.3},~Prop.~3, p.~182]{Dix81} and~\cite[\S{II.2.3},~Prop.~4(ii), p.~183]{Dix81} by approximation of~$\phi$ with suitable polynomials, since~$\sigma(B)$ is compact.
\end{proof}

\begin{remark} Arguing with Tietze Extension Theorem, it is possible to show that the direct-integral representation~\eqref{eq:CFC} of~$\phi(B)$ only depends on the values of~$\phi$ on the spectrum~$\sigma(B)$ of~$B$.
\end{remark}

\subsection{Direct integrals of quadratic forms}\label{ss:DirIntQ}
The main object of our study are direct integrals of quadratic forms. Before turning to the case of Dirichlet forms on concrete Hilbert spaces ($L^2$-spaces), we give the main definitions in the general case of quadratic forms on abstract Hilbert spaces.

\begin{definition}[Direct integral of quadratic forms]\label{d:assQ}
Let~$(Z,\mcZ,\nu)$ be a $\sigma$-finite countably generated measure space. For~$\zeta\in Z$ let~$(Q_\zeta,D_\zeta)$ be a closable (densely defined) quadratic form on a Hilbert space~$H_\zeta$. We say that~$\zeta\mapsto (Q_\zeta,D_\zeta)$ is a \emph{$\nu$-measurable field of quadratic forms on~$Z$} if
\begin{enumerate}[$(a)$]
\item\label{i:assQ:a} $\zeta\mapsto H_\zeta$ is a $\nu$-measurable field of Hilbert spaces on~$Z$ with underlying space~$S_H$;
\item\label{i:assQ:b} $\zeta\mapsto \dom{Q_\zeta}_1$ is a $\nu$-measurable field of Hilbert spaces on~$Z$ with underlying space~$S_Q$;
\item\label{i:assQ:c} $S_Q$ is a linear subspace of~$S_H$ under the identification of~$\dom{Q_\zeta}$ with a subspace of~$H_\zeta$ granted by Lemma~\ref{l:WLSC}.
\end{enumerate}

We denote by
\begin{align*}
Q=\dint[S_Q]{Z} Q_\zeta \diff\nu(\zeta)
\end{align*}
the \emph{direct integral} of~$\zeta\mapsto (Q_\zeta,\dom{Q_\zeta})$, i.e. the quadratic form defined on~$H$ as in~\eqref{eq:DirInt} given by
\begin{equation}\label{eq:DIntQF}
\begin{aligned}
\dom{Q}\eqdef& \set{\class[H]{\rep u} : \rep u \in S_Q , \int_Z Q_{\zeta,1}(\rep u_\zeta) \diff\nu(\zeta) <\infty} \comma
\\
Q(u,v)\eqdef& \int_Z Q_\zeta(u_\zeta,v_\zeta) \diff\nu(\zeta)\comma \qquad u,v\in\dom{Q}\fstop
\end{aligned}
\end{equation}
\end{definition}

\begin{remark}[Separability]\label{r:Separability} It is implicit in our definition of $\nu$-measurable field of Hilbert spaces that~$H_\zeta$ is \emph{separable} for every~$\zeta\in Z$. Therefore, when considering $\nu$-measurable fields of domains as in Definition~\ref{d:assQ}\iref{i:assQ:b}, $\dom{Q_\zeta}_1$ is taken to be $(Q_\zeta)^{1/2}_1$-separable \emph{by assumption}.
\end{remark}

\begin{proposition}\label{p:DirInt}
Let~$(Q,\dom Q)$ be a direct integral of quadratic forms. Then,
\begin{enumerate}[$(i)$]
\item\label{i:p:DirInt1} $\ttonde{Q,\dom Q}$ is a densely defined closed quadratic form on~$H$;
\item\label{i:p:DirInt2} $\zeta\mapsto G_{\zeta,\alpha}$,~$\zeta\mapsto T_{\zeta,t}$ are $\nu$-measurable fields of bounded operators for every~$\alpha,t>0$;
\item\label{i:p:DirInt3} $Q$ has resolvent and semigroup respectively defined by
\begin{equation}\label{eq:p:DirInt:0}
\begin{aligned}
G_\alpha \eqdef & \dint[S_H]{Z} G_{\zeta,\alpha} \diff\nu(\zeta) \comma\qquad \alpha>0\semicolon
\\
T_t \eqdef & \dint[S_H]{Z} T_{\zeta,t} \diff\nu(\zeta) \comma\;\;\qquad t>0\fstop
\end{aligned}
\end{equation}
\end{enumerate}
\end{proposition}
\begin{proof} \iref{i:p:DirInt1}
Since~$\zeta\mapsto H_\zeta$ is a $\nu$-measurable family of Hilbert spaces by Definition~\ref{d:assQ}\iref{i:assQ:a}, the map~$\zeta\mapsto \norm{u_\zeta}_{\zeta}$ is $\nu$-measurable for every~$u\in S_H$ by Definition~\ref{d:DirInt}\iref{i:d:DirInt1}.
Analogously, the map~$\zeta\mapsto Q_{\zeta,1}^{1/2}(u_\zeta)$ is $\nu$-measurable for every~$u\in S_Q$. Together with the polarization identity for~$\dom{Q}_1$, this yields the measurability of the maps
\begin{align*}
\zeta\mapsto Q_{\zeta,\alpha}(u_\zeta,v_\zeta)\comma \qquad u,v\in\dom{Q}\comma \qquad \alpha>0 \fstop
\end{align*}
As a consequence,~$\zeta\mapsto \dom{Q_\zeta}_\alpha$ is a $\nu$-measurable field of Hilbert spaces (on~$Z$, with underlying space~$S_Q$) for every~$\alpha>0$. Thus, it admits a direct integral of Hilbert spaces
\begin{align*}
D_\alpha\eqdef \dint[S_Q]{Z} \dom{Q_{\zeta}}_{\alpha} \diff\nu(\zeta) \comma \qquad \alpha>0\fstop
\end{align*}

For~$\alpha>0$ let~$\seq{u^\alpha_n}_n$ be a fundamental sequence of $\nu$-measurable vector fields for~$D_\alpha$ and~$\seq{u_n}_n$ be a fundamental sequence of $\nu$-measurable vector fields for~$H$.
Since~$(Q_\zeta,D_\zeta)$ is closable on~$H_\zeta$ for every~$\zeta\in Z$, the extension of the canonical inclusion~$\iota_{\zeta,1}\colon \dom{Q_\zeta}_1\rar H_\zeta$ is injective and non-expansive\footnote{We say that a linear map~$\iota\colon H_1\to H_2$ between Hilbert spaces~$(H_1,\norm{\emparg}_1)$ and~$(H_2,\norm{\emparg}_2)$ is \emph{non-expansive} if~$\norm{\iota h}_2\leq \norm{h}_1$ for every~$h\in H_1$.} for every~$\zeta\in Z$ by Lemma~\ref{l:WLSC}. 
By Definition~\ref{d:assQ},~$D_\alpha$ and~$H$ are defined on the same underlying space~$S$. Therefore, the maps
\begin{align*}
\zeta\mapsto\scalar{\iota_{\zeta,\alpha} u^\alpha_{i,\zeta}}{u_{j,\zeta}}_{\zeta}=\scalar{u^\alpha_{i,\zeta}}{u_{j,\zeta}}_{\zeta}\comma \qquad i,j\in\N\comma \qquad \alpha>0\comma
\end{align*}
are $\nu$-measurable. Together with the uniform boundedness of~$\iota_{\zeta,\alpha}$ in~$\zeta\in Z$, this yields the decomposability of the operator~$\iota_\alpha\colon D_\alpha\rar H$ defined by
\begin{align*}
\iota_\alpha\eqdef \dint[S_Q]{Z} \iota_{\zeta,\alpha} \diff\nu(\zeta)\comma \qquad \alpha>0 \fstop
\end{align*}

By~\cite[\S{II.2.3}, Example,~p.~182]{Dix81} and the injectivity of~$\iota_{\zeta,\alpha}$ for every~$\zeta\in Z$ and every~$\alpha>0$, the operator~$\iota_\alpha\colon D_\alpha\to H$ is injective. In particular, the composition of~$\iota_1$ with the inclusion of~$\dom{Q}$ into~$H$ is injective, thus~$Q$ is closed.
Finally, since~$\tseq{u^\alpha_{n,\zeta}}_n$ is $Q_{\zeta,\alpha}^{1/2}$-total in~$\dom{Q_\zeta}_\alpha$ for every~$\zeta\in Z$ by Definition~\ref{d:DirInt}\iref{i:d:DirInt3}, it is additionally $H_\zeta$-total for every~$\zeta\in Z$ by $H_\zeta$-density of~$\dom{Q_\zeta}$ in~$H_\zeta$. As a consequence,~$\tseq{u^\alpha_n}_n$ is fundamental also for~$H$, thus~$\dom{Q}$ is $H$-dense in~$H$.

\iref{i:p:DirInt2} For fixed~$\alpha>0$ consider the field of linear operators~$\zeta\mapsto G_{\zeta,\alpha}$. The map (cf.~\eqref{eq:Hille--YosidaG})
\begin{align*}
\zeta\mapsto Q_{\zeta,\alpha}(G_{\zeta,\alpha} u^\alpha_{i,\zeta},u^\alpha_{j,\zeta})=\scalar{u^\alpha_{i,\zeta}}{u^\alpha_{j,\zeta}}_{\zeta}
\end{align*}
is $\nu$-measurable for every~$i,j\in \N$ since~$u^\alpha_n$ is a $\nu$-measurable vector field.
Since $\norm{G_{\alpha,\zeta}}_\op\leq \alpha^{-1}$ and~$\seq{u^\alpha_n}_n$ is a fundamental sequence of $\nu$-measurable vector fields for~$H$, then~$\zeta\mapsto G_{\zeta,\alpha}$ is a $\nu$-measurable field of bounded operators by~\cite[\S{II.2.1}, Prop.~1, p.~179]{Dix81} and the operator~$G_\alpha$ defined in~\eqref{eq:p:DirInt:0} is decomposable for every~$\alpha>0$.

By Lemma~\ref{l:DecomposabilityCFC} any image of~$G_\alpha$ via its continuous functional calculus is itself decomposable.

For every~$\zeta\in Z$ one has~$T_{\zeta, t}=\lim_{\beta\rar \infty} e^{t\beta(\beta G_{\zeta,\beta}-1)}$ strongly in~$H_\zeta$ by~\eqref{eq:Hille--YosidaT}, hence
\begin{align*}
\zeta\mapsto \scalar{T_{\zeta,t} u^\alpha_{i,\zeta}}{u^\alpha_{j,\zeta}}_{\zeta}=\lim_{\beta\rar \infty}\scalar{e^{t\beta(\beta G_{\zeta,\beta}-1)} u^\alpha_{i,\zeta}}{u^\alpha_{j,\zeta}}_{\zeta}
\end{align*}
is a pointwise limit of $\nu$-measurable functions, thus it is $\nu$-measurable, for every~$i,j\in\N$ and every~$t>0$. As a consequence,~$\zeta\mapsto T_{\zeta,t}$ is a $\nu$-measurable field of bounded operators for every~$t>0$, again by~\cite[\S{II.2.1}, Prop.~1, p.~179]{Dix81}.
%

\iref{i:p:DirInt3} It suffices to show~\eqref{eq:Hille--Yosida} for~$(Q,\dom{Q})$,~$G_\alpha$ and~$T_t$ defined in~\eqref{eq:p:DirInt:0}.
Now, by definition of~$(Q,\dom{Q})$ one has for every~$\alpha>0$
\begin{align*}
Q_\alpha(G_\alpha u,v)=&\int_Z Q_\zeta\ttonde{(G_\alpha u)_\zeta,v_\zeta} \diff\nu(\zeta) +\alpha\int_Z \scalar{(G_\alpha u)_\zeta}{v_\zeta}_{\zeta} \diff\nu(\zeta)
\\
=&\int_Z Q_{\zeta,\alpha}\ttonde{(G_\alpha u)_\zeta,v_\zeta} \diff\nu(\zeta)
\\
=&\int_Z Q_{\zeta,\alpha}\ttonde{G_{\alpha,\zeta} u_\zeta,v_\zeta} \diff\nu(\zeta) \fstop
\intertext{By~\cite[\S{II.2.3}, Cor.,~p.~182]{Dix81} and decomposability of~$G_\alpha$, one has~$G_{\alpha,\zeta}=G_{\zeta,\alpha}$ for~$\nu$-a.e.~$\zeta\in Z$, whence, by~\eqref{eq:Hille--YosidaG} applied to~$(Q_\zeta,\dom{Q_\zeta})$ and~$G_{\zeta,\alpha}$,}
Q_\alpha(G_\alpha u,v)=&\int_Z \scalar{u_\zeta}{v_\zeta}_{\zeta} \diff\nu(\zeta)
=\scalar{u}{v}\comma
\end{align*}
which concludes the proof of~\eqref{eq:Hille--YosidaG} for~$G_\alpha$.

Let us show~\eqref{eq:Hille--YosidaT} for~$T_t$. 
Define the operators
\begin{align*}
T^\sym{\beta}_t\eqdef \dint[]{Z} e^{t\beta(\beta G_{\zeta,\beta}-1)}\diff\nu(\zeta)\comma \qquad \beta,t>0\comma
\end{align*}
and note that~$\sup_{\beta>1}\norm{T^\sym{\beta}_t}_\op<\infty$ for every~$t>0$.
By~\eqref{eq:Hille--YosidaT} applied to~$T_{\zeta,t}$ for every~$\zeta\in Z$ and every~$t>0$ we have that~$T_{\zeta, t}=\lim_{\beta\rar \infty} e^{t\beta(\beta G_{\zeta,\beta}-1)}$ strongly in~$H_\zeta$.
On the one hand, we may now apply~\cite[\S{II.2.3}, Prop.~4(ii), p.~183]{Dix81} to conclude that
\begin{align}\label{eq:p:DirInt:1}
H\text{-}\lim_{\beta\to\infty}T^\sym{\beta}_t = \dint[]{Z} \tonde{H_\zeta\text{-}\lim_{\beta\to\infty} e^{t\beta(\beta G_{\zeta,\beta}-1)} }\diff\nu(\zeta) = \dint[]{Z} T_{\zeta,t}\diff\nu(\zeta)\defeq T_t\fstop
\end{align}
strongly in~$H$.
On the other hand, by Lemma~\ref{l:DecomposabilityCFC} we have that
\begin{align}\label{eq:p:DirInt:2}
T^\sym{\beta}_t\eqdef \dint[]{Z} e^{t\beta(\beta G_{\zeta,\beta}-1)}\diff\nu(\zeta)=e^{t\beta(\beta G_\beta-1)}\comma \qquad \beta,t>0 \fstop
\end{align}
Taking the strong $H$-limit as~$\beta\to\infty$ in~\eqref{eq:p:DirInt:2} yields the assertion by comparison with~\eqref{eq:p:DirInt:1}.
\end{proof}

\begin{remark}[{cf.~\cite[\S{II.1.3}, Rmk.~3\ p.~166 and \S{II.1.4}, Rmk.\ p.~168]{Dix81}}]\label{r:Dixmier}
Each of the above statements holds with identical proof if one substitutes `$\nu$-measurable' with `measurable'.
\end{remark}

\begin{remark}\label{r:S}
Under the assumptions of Proposition~\ref{p:DirInt}, assertion~\ref{i:p:DirInt1} of the same Proposition implies that the space of $\nu$-measurable vector fields~$S_H$ is uniquely determined by~$S_Q$ as a consequence of Proposition~\ref{p:Dix}. Thus, everywhere in the following when referring to a direct integral of quadratic forms we shall ---~with abuse of notation~--- write $S$ in place of both~$S_H$ and~$S_Q$.
\end{remark}

The next proposition completes the picture, by providing a direct-integral representation for the generator of the form~$(Q,\dom{Q})$ in~\eqref{eq:DIntQF}. Since we shall not need this result in the following, we omit an account of direct integrals of unbounded operators, referring the reader to~\cite[\S1]{Lan75}. Once the necessary definitions are established, a proof is straightforward.

\begin{proposition} Let~$(Q,\dom{Q})$ be defined as in~\eqref{eq:DIntQF}. Then, its generator $(L,\dom{L})$ has the direct-integral representation
\begin{align}\label{eq:DirIntGenerators}
L=\dint[]{Z} L_\zeta \diff\nu(\zeta) \fstop
\end{align}
\end{proposition}

\begin{remark}[Comparison with~\cite{AlbRoe90}] As for quadratic forms,~\eqref{eq:DirIntGenerators} is understood as a direct-integral representation of the Hilbert space~$\dom{L}$, endowed with the graph norm, by the measurable field of Hilbert spaces~$\zeta\mapsto \dom{L_\zeta}$, each endowed with the relative graph norm. The set-wise identification of~$\dom{L}$ as a linear subspace of~$H$ as in~\eqref{eq:DirInt} is already shown in~\cite[Prop.~1.6]{AlbRoe90}.
\end{remark}

\subsection{Dirichlet forms}
We recall a standard setting for the theory of Dirichlet forms, following~\cite{MaRoe92}.

\begin{assumption}\label{ass:Main} The quadruple $(X,\T,\mcX,\mu)$ is so that
\begin{enumerate}[$(a)$]
\item\label{i:ass:Main:1} $(X,\T)$ is a metrizable Luzin space with Borel $\sigma$-algebra~$\mcX$;
\item\label{i:ass:Main:2} $\hat\mu$ is a Radon measure on~$(X,\T,\mcX^\mu)$ with full support.
\end{enumerate}
\end{assumption}

By~\cite[415D(iii), 424G]{Fre00} any space~$(X,\mcX,\mu)$ satisfying Assumption~\ref{ass:Main} is, in particular, $\sigma$-finite standard. The \emph{support} of a ($\mu$-)measurable function~$f\colon X\rar \R$ (possibly defined only on a $\mu$-conegligible set) is defined as the measure-theoretical support~$\supp[\abs{f}\cdot \mu]$. Every such~$f$ has support, independent of the $\mu$-representative of~$f$, cf.~\cite[p.~148]{MaRoe92}.

A closed positive semi-definite quadratic form~$(Q,\dom{Q})$ on~$L^2(\mu)$ is a (\emph{symmetric}) \emph{Dirichlet form} if
\begin{align}\label{eq:subMarkov}
f\in\dom{Q} \implies f^+ \wedge 1 \in \dom{Q} \text{~and~} Q(f^+\wedge 1)\leq Q(f) \fstop
\end{align}
We shall denote Dirichlet forms by~$(E,\dom{E})$. A Dirichlet form~$(E,\dom{E})$ is
\begin{itemize}
\item \emph{local} if~$E(f,g)=0$ for every~$f,g\in\dom{E}$ with~$\supp[f]\cap \supp[g]=\emp$;

\item \emph{strongly local} if~$E(f,g)=0$ for every~$f,g\in\dom{E}$ with~$g$ constant on a neighborhood of~$\supp[f]$;

\item \emph{regular} if~$(X,\T)$ is (additionally) locally compact and there exists a \emph{core}~$\mcC$ for~$(E,\dom{E})$, i.e.\ a subset $\mcC\subset \dom{E}\cap \Cz(\T)$ both $E^{1/2}_1$-dense in~$\dom{E}$ and uniformly dense in~$\Cz(\T)$.
\end{itemize}

On spaces that are not necessarily locally compact, the interplay between a Dirichlet form~$(E,\dom{E})$ and the topology~$\T$ on~$X$ is specified by the following definitions.
For an increasing sequence~$\seq{F_k}_k$ of Borel subsets~$F_k\subset X$ set
\begin{align*}
\dom{E,\seq{F_k}_k}\eqdef \set{f\in \dom{E}: f\equiv 0 \as{\mu} \text{ on } F_k^\complement \text{ for some } k\in \N}\fstop
\end{align*}
The sequence~$\seq{F_k}_k$ is called an $E$-\emph{nest} if each~$F_k$ is closed and~$\msD(E,\seq{F_k}_k)$ is dense in~$\dom{E}_1$. A set $N\subset X$ is \emph{$E$-exceptional} if $N\subset \cap_k F_k^\complement$ for some $E$-nest~$\seq{F_k}_k$. A property of points in~$X$ holds \emph{$E$-quasi-everywhere} if it holds for every point in~$N^\complement$ for some $E$-exceptional set~$N$. A function~$f\colon X\rar \R$ is \emph{$E$-quasi-continuous} if there exists an $E$-nest~$\seq{F_k}_k$ so that the restriction of~$f$ to~$F_k$ is continuous for every~$k\in \N$.
Finally, a form~$(E,\dom{E})$ is
\begin{itemize}
\item \emph{quasi-regular} on~$(X,\T)$ if there exist:
\begin{enumerate*}[$(\mathsc{qr}_1)$]
\item an $E$-nest~$\seq{F_k}_k$ of compact sets;
\item an $E_1^{1/2}$-dense subset of~$\dom{E}_1$ the elements of which all have an $E$-quasi-continuous $\mu$-version;
\item an $E$-ex\-ception\-al set~$N\subset X$ and a countable family~$\seq{f_n}_n$ of $E$-quasi-continuous functions~$f_n\in\dom{E}$ separating points in~$N^\complement$.
\end{enumerate*}
\end{itemize}
We refer to~\cite{CheMaRoe94} or~\cite[\S{A.4}]{FukOshTak11} for the notion of \emph{quasi-homeomorphism} of Dirichlet forms.

We say that~$(E,\dom{E})$ has \emph{carr\'e du champ operator}~$(\cdc,\dom{E})$, if $\cdc\colon \dom{E}^\otym{2}\rar L^1(\mu)$ is a non-negative definite symmetric continuous bilinear operator so that
\begin{align}\label{eq:CdC}
E(f,gh)+E(fh,g)-E(fg,h)=2\int_X h\, \cdc(f,g)\diff\mu \comma \qquad f,g,h\in\dom{E}\cap L^\infty(\mu) \fstop
\end{align}

Finally, let~$\domext{E}$ be the linear space of all functions~$u\in L^0(\mu)$ so that there exists an $E^{1/2}$-Cauchy sequence~$\seq{u_n}_n\subset \dom{E}$ with~$\nlim u_n=u$ $\mu$-a.e. We denote by~$(\domext{E}, E)$ the space~$\domext{E}$ endowed with the extension of~$E$ to~$\domext{E}$ called the \emph{extended Dirichlet space} of~$(E,\dom{E})$. For proofs of well-posedness in this generality, see~\cite[p.~693]{Kuw98}.
If~$(E,\dom{E})$ has semigroup~$T_\bullet\colon L^2(\mu)\rar L^2(\mu)$, we denote as well by~$T_\bullet\colon L^\infty(\mu)\rar L^\infty(\mu)$ the extension of the semigroup to~$L^\infty(\mu)$.
We say that~$(E,\dom{E})$ is
\begin{itemize}
\item \emph{conservative} if~$T_t \car=\car$ $\mu$-a.e.\ for all~$t\geq 0$;
\item \emph{transient} if~$\domext{E}$ is a Hilbert space with inner product~$E$;
\item \emph{recurrent} if~$\car\in\domext{E}$ and~$E(\car)=0$.
\end{itemize}
These definitions are equivalent to the standard ones (e.g.~\cite[p.~55]{FukOshTak11}) by~\cite[Thm.s~1.6.2,~1.6.3, p.~58]{FukOshTak11}, a proof of which may be adapted to the case of spaces satisfying Assumption~\ref{ass:Main}.

\subsection{Direct-integral representation of \texorpdfstring{$L^2$}{L2}-spaces}\label{ss:DirIntE}
In order to introduce direct-integral representations of Dirichlet forms, we need to construct direct integrals of concrete Hilbert spaces in such a way to additionally preserve the Riesz structure of Lebesgue spaces implicitly used to phrase the sub-Markovianity property~\eqref{eq:subMarkov}. To this end, we shall need the concept of a disintegration of measures.

\paragraph{Disintegrations} Let~$(X,\mcX,\mu)$ and~$(Z,\mcZ,\nu)$ be (non-trivial) measure spaces. A map~$s\colon (X,\mcX)\rar (Z,\mcZ)$ is \emph{inverse-measure-preserving} if~$s_\pfwd \mu\eqdef \mu \circ s^{-1}=\nu$.
Hereafter, we fix an inverse-measure-preserving map~$s\colon (X,\mcX)\rar (Z,\mcZ)$.

\begin{definition}[Disintegrations,~{\cite[452E]{Fre00}}]\label{d:Disint}
A \emph{pseudo-disintegration of~$\mu$ over~$\nu$} is any family of non-zero measures~$\seq{\mu_\zeta}_{\zeta\in Z}$ on~$(X,\mcX)$ so that~$\zeta\mapsto \mu_\zeta A$ is $\nu$-measurable and
\begin{align*}
\mu A=\int_Z\mu_\zeta A \diff\nu(\zeta)\comma \qquad A\in\mcX\fstop
\end{align*}
A pseudo-disintegration is
\begin{itemize}
\item \emph{separated} if there exists a family of pairwise disjoint sets~$\seq{A_\zeta}_{\zeta\in Z}\subset \mcX^\mu$ so that~$A_\zeta$ is $\mu_\zeta$-conegligible for $\nu$-a.e.~$\zeta\in Z$, henceforth called a \emph{separating family} for~$\seq{\mu_\zeta}_{\zeta\in Z}$.
\end{itemize}

A \emph{disintegration of~$\mu$ over~$\nu$} is a pseudo-disintegration additionally so that~$\mu_\zeta$ is a sub-probability measure for every~$\zeta\in Z$.
A disintegration is
\begin{itemize}
\item $\nu$-\emph{essentially unique} if the measures~$\mu_\zeta$ are uniquely determined for $\nu$-a.e.~$\zeta\in Z$;

\item \emph{consistent with~$s$} if
\begin{align}\label{eq:Disint:0}
\mu \ttonde{A\cap s^{-1}(B)}=\int_B\mu_\zeta A \diff\nu(\zeta)\comma \qquad A\in\mcX\comma B\in \mcZ\semicolon
\end{align}

\item \emph{strongly consistent with~$s$} if $s^{-1}(\zeta)$ is~$\mu_\zeta$-conegligible for $\nu$-a.e.~$\zeta\in Z$.
\end{itemize}
\end{definition}

If~$\seq{\mu_\zeta}_{\zeta\in Z}$ is a pseudo-disintegration of~$\mu$ over~$\nu$, then
\begin{align}\label{eq:DisintF}
\int_X g\diff\mu=\int_Z\int_X g(x) \diff\mu_\zeta(x) \diff\nu(\zeta)
\end{align}
whenever the left-hand side makes sense,~\cite[452F]{Fre00}. We note that a disintegration~$\seq{\mu_\zeta}_{\zeta\in Z}$ of~$\mu$ over~$\nu$ strongly consistent with a map~$s$ is automatically separated, with separating family~$\seq{s^{-1}(\zeta)}_{\zeta\in Z}$.

\paragraph{Direct integrals and disintegrations} Let~$(X,\mcX,\mu)$ be $\sigma$-finite standard,~$(Z,\mcZ,\nu)$ be $\sigma$-finite countably generated, and~$\seq{\mu_\zeta}_{\zeta\in Z}$ be a pseudo-disintegration of~$\mu$ over~$\nu$. Denote by
\begin{itemize}
\item $\mcL^0(\mu)$ the space of $\mu$-measurable real-valued functions (\emph{not}: $\mu$-classes) on~$X$;
\item $\mcL^\infty(\mu)$ the space of uniformly bounded (\emph{not}: $\mu$-essentially uniformly bounded) functions in~$\mcL^0(\mu)$;
\item $\mcL^p(\mu)$ the space of $p$-integrable functions in~$\mcL^0(\mu)$.
\end{itemize}

For a family~$\mcA\subset\mcL^0(\mu)$, let~$\class[\mu]{\mcA}$ denote the family of the corresponding $\mu$-classes.

\smallskip

Let now~$F\eqdef \prod_{\zeta\in Z} L^2(\mu_\zeta)$. We denote by~$\delta\colon \mcL^2(\mu) \rar F$ the diagonal embedding of~$\mcL^2(\mu)$ into~$F$, regarded up to $\mu_\zeta$-classes, viz.~$\delta\colon f\mapsto \ttseq{\zeta\mapsto \delta(f)_\zeta}$, where
\begin{align}\label{eq:Delta}
\delta(f)_\zeta\eqdef \begin{cases} \class[\mu_\zeta]{f} & \text{if~} f\in \mcL^2(\mu_\zeta) \\ \zero_{L^2(\mu_\zeta)} & \text{otherwise}\end{cases} \fstop
\end{align}
In general, it does not hold that~$f\in \mcL^2(\mu_\zeta)$ for \emph{every}~$\zeta\in Z$, thus we need to adjust the obvious definition of~$\delta(f)$ as above in such a way that~$\delta(f)\in F$, that is~$\delta(f)_\zeta\in\mcL^2(\mu_\zeta)$ for every~$\zeta\in Z$. Note that~$\delta$ is thus \emph{not} linear. However, since~$f\in \mcL^2(\mu)$, then~$\delta(f)_\zeta=\class[\mu_\zeta]{f}$ for~$\nu$-a.e.~$\zeta\in Z$ by~\eqref{eq:DisintF}. It will be shown in Proposition~\ref{p:DIntL^2} below that~$\delta$ is well-defined as linear morphism mapping~$\mu$-classes to $H$-classes.

Further let~$\mcA$ be satisfying
\begin{equation}\label{eq:AssA}
\text{$\mcA$ is a linear subspace of~$\mcL^2(\mu)$, and~$\class[\mu]{\mcA}$ is dense in~$L^2(\mu)$.}
\end{equation}

Since~$\class[\mu]{\mcA}$ is dense in $L^2(\mu)$ and the latter is separable, then there exists a countable family~$\mcU\subset \mcA$ so that~$\class[\mu_\zeta]{\mcU}$ is total in~$L^2(\mu_\zeta)$ for $\nu$-a.e.~$\zeta\in Z$.
Thus for every~$\mcA$ as in~\eqref{eq:AssA} there exists a unique space of $\nu$-measurable vector fields~$S=S_\mcA$ containing~$\delta(\mcA)$, generated by~$\delta(\mcA)$ in the sense of Proposition~\ref{p:Dix}. We denote by~$H$ the corresponding direct integral of Hilbert spaces
\begin{align}\label{eq:H}
H\eqdef \dint[S]{Z} L^2(\mu_\zeta) \diff\nu(\zeta)\fstop
\end{align}

Since~$S$ is unique, it is in fact independent of~$\mcA$. Indeed, let~$\mcA_0$,~$\mcA_1$ be satisfying~\eqref{eq:AssA} and note that~$\mcA\eqdef\mcA_0\oplus \mcA_1$ satisfies~\eqref{eq:AssA} as well. Thus,~$\delta(\mcA_0)$, $\delta(\mcA_1)\subset S_{\mcA}$, and so~$S_\mcA= S_{\mcA_0}=S_{\mcA_1}$ by uniqueness.

\begin{remark}[{cf.~\cite[\S7.2, p.~84]{HayMirYve91}}]\label{r:Order}
The direct integral $H$ constructed in~\eqref{eq:H} is a \emph{Banach lattice} (e.g.~\cite[354A(b)]{Fre00}) for the order
\begin{align*}
h\geq  \zero_H \qquad \iff \qquad h_\zeta\geq \zero_{L^2(\mu_\zeta)} \quad \forallae{\nu} \zeta\in Z \fstop
\end{align*}
%
In particular, for every~$g,h\in H$, the fields~$h^+$,~$h^-$,~$g\wedge h$, and~$g\vee h$, respectively defined by
\begin{align*}
h^\pm \colon \zeta\mapsto h_\zeta^\pm\comma \qquad g\wedge h\colon \zeta\mapsto (g_\zeta\wedge h_\zeta)\comma \qquad g\vee h\colon \zeta\mapsto (g_\zeta\vee h_\zeta)\comma
\end{align*}
are $\nu$-measurable fields representing elements of~$H$.
In the following, we shall occasionally write ---~here,~$\uno_H$ is merely a shorthand~---
\begin{align*}
\zero_H\leq h \leq \uno_H 
\end{align*}
to indicate that
\begin{align*}
0 \leq h_\zeta\leq 1 \as{\mu_\zeta} \quad \forallae{\nu} \zeta\in Z \fstop
\end{align*}
\end{remark}

\begin{remark}\label{r:Algebra}
For arbitrary measurable spaces, the standard choice for~$\mcA$ is the algebra of $\mu$-integrable simple functions. If~$(X,\T,\mcX,\mu)$ were a locally compact Polish Radon measure space, one might take for instance~$\mcA=\Cc(\T)$, the algebra of continuous compactly supported functions.
In fact, for the purposes of the present section, we might as well choose~$\mcA=\mcL^2(\mu)$, as the largest possible choice, or~$\mcA$ a countable $\mbbQ$-vector subspace of~$\mcL^1(\mu)\cap \mcL^\infty(\mu)$, as a smallest possible one. When dealing with direct integrals of regular Dirichlet forms however, the natural choice for~$\mcA$ is that of a special standard core~$\mcC$ for the resulting direct-integral form.
\end{remark}

\begin{remark}[Comparison with~{\cite{AlbRoe90}}]\label{r:AlbRoe90}
We note that for every~$\mcA$ as in~\eqref{eq:AssA},~$\class[\mu]{\mcA}$ is a \emph{determining class} in the sense of~\cite[p.~402]{AlbRoe90}. Conversely, every determining class~$\msL_0$ is contained in a minimal linear space of functions~$\class[\mu]{\mcA}$ satisfying~\eqref{eq:AssA}.
\end{remark}

\begin{remark}[Caveat]\label{r:Caveat} Whereas the space~$H$ does not depend on~$\mcA$, in general it \emph{does} depend on~$S_\mcA$, cf.~Rmk.~\ref{r:S0}. Furthermore,~$H$ depends on the chosen pseudo-disintegration too, and thus~$H$ need \emph{not} be isomorphic to~$L^2(\mu)$, as shown in the next example.
\end{remark}

\begin{example}
Let~$\set{*}$ denote the one-point space, set~$\mu\eqdef 2\delta_*$, and note that $L^2(\mu)\cong \R$. On the other hand, if~$Z\eqdef (\set{0,1},\nu)$ is the two-point space with uniform measure~$\nu$, and~$\mu_\zeta\eqdef \delta_*$ for~$\zeta\in Z$, then~$\seq{\mu_\zeta}_{\zeta\in Z}$ is a (pseudo-)disintegration of~$\mu$, yet~$H\cong L^2(\mu_0)\oplus L^2(\mu_1)\cong \R^2$. 
\end{example}

\begin{proposition}\label{p:DIntL^2}
Let~$(X,\mcX,\mu)$ be $\sigma$-finite standard,~$(Z,\mcZ,\nu)$ be $\sigma$-finite countably generated, and~$\seq{\mu_\zeta}_{\zeta\in Z}$ be a pseudo-disintegration of~$\mu$ over~$\nu$. Then, the morphism
\begin{equation}\label{eq:DirIntL^2:0}
\begin{aligned}
\iota\colon L^2(\mu) &\longrar H\eqdef \dint[S]{Z} L^2(\mu_\zeta) \diff\nu(\zeta)
\\
\class[\mu]{f} &\longmapsto \class[H]{\delta(f)}
\end{aligned}
\end{equation}
\begin{enumerate}[$(i)$]
\item\label{i:p:DIntL^2:1} is well-defined, linear, and an isometry of Hilbert spaces, additionally unitary if~$\seq{\mu_\zeta}_{\zeta\in Z}$ is separated;

\item\label{i:p:DIntL^2:2} is a Riesz homomorphism (e.g.~\cite[351H]{Fre00}). In particular,
\begin{itemize}
\item for each~$f\in\mcL^2(\mu)$, it holds that~$(\iota\class[\mu]{f})_\zeta \geq \zero_{L^2(\mu_\zeta)}$ for~$\nu$-a.e.~$\zeta\in Z$ if and only if~$f\geq 0$ $\mu$-a.e.;

\item for each~$f,g\in\mcL^2(\mu)$, it holds that~$(\iota\class[\mu]{f \wedge g})_\zeta= (\iota\class[\mu]{f})_\zeta \wedge (\iota\class[\mu]{g})_\zeta$ for~$\nu$-a.e.~$\zeta\in Z$.
\end{itemize}
\end{enumerate}
\end{proposition}

\begin{proof} As usual, we denote by~$\norm{\emparg}$ the norm on~$H$, and by~$\norm{\emparg}_2$, resp.~$\norm{\emparg}_{2,\zeta}$, the norm on~$L^2(\mu)$, resp.~$L^2(\mu_\zeta)$. 
Let~$\mcA$ be satisfying~\eqref{eq:AssA}, and define a map~$\hat\iota\colon \mcA\rar H$ by~$\hat\iota\colon f\mapsto \class[H]{\delta(f)}$.

By definition of~$\hat\iota$ and~$\delta$, by definition~\eqref{eq:NormH} of~$\norm{\emparg}$, and by the property~\eqref{eq:DisintF} of the disintegration,
\begin{equation*}
\begin{aligned}
\norm{\hat\iota f_1-\hat\iota f_2}^2=&\ \int_Z \norm{f_1-f_2}^2_{2,\zeta} \diff\nu(\zeta)
\\
=&\ \int_Z\int_X \abs{f_1-f_2}^2 \diff\mu_\zeta\diff\nu(\zeta)
\\
=&\norm{f_1-f_2}_2^2
\end{aligned}\qquad\qquad f_1,f_2\in\mcA\fstop
\end{equation*}
As a consequence,~$\hat\iota\colon \mcA\rar H$ descends to a \emph{linear} isometry~$\iota\colon \class[\mu]{\mcA}\rar H$, and the latter extends to the non-relabeled (linear) isometry~\eqref{eq:DirIntL^2:0} by density of~$\class[\mu]{\mcA}$ in~$L^2(\mu)$.

Assume now that~$\seq{\mu_\zeta}_{\zeta\in Z}$ is separated with separating family~$\seq{A_\zeta}_{\zeta\in Z}$, and fix~$h\in (\im\,\iota)^\perp$. Let~$\rep h\in S$ be an $H$-representative of~$h$. For each~$\zeta\in Z$, let~$\reptwo h_\zeta\in \mcL^2(\mu_\zeta)$ be a representative of~$\rep h_\zeta$, and define a function~$\reptwo h\colon X\rar \R$ by
\begin{align*}
\reptwo h(x)\eqdef \begin{cases} \reptwo h_\zeta(x) & \text{if~} x\in A_\zeta\comma \zeta\in Z\comma \\ 0 & \text{otherwise}\end{cases} \fstop
\end{align*}
This definition is well-posed since the sets~$A_\zeta$'s are pairwise disjoint.

\smallskip

\emph{Claim: $\reptwo h\equiv 0$ $\mu$-a.e.} With slight abuse of notation, set~$\delta(\reptwo h)_\zeta\eqdef \ttclass[\mu_\zeta]{\reptwo h}$ for~$\zeta\in Z$, and~$\delta(\reptwo h)\eqdef \ttseq{\zeta\mapsto \delta(\reptwo h)_\zeta}$. 
By construction,~$\delta(\reptwo h)=\rep h$, therefore~$\delta(\reptwo h)\in S$, and so
\begin{align*}
0=\tscalar{h}{\iota \class[\mu]{f}}=\int_Z \int_X \reptwo h f \diff\mu_\zeta \diff\nu(\zeta) \comma \qquad f\in \mcL^2(\mu) \comma
\end{align*}
where the right-hand side is well-defined since~$\rep h\in S$. As a consequence,
\begin{align*}
f\mapsto \int_Z\int_X \reptwo h f\diff\mu_\zeta \diff\nu(\zeta)
\end{align*}
is the $\zero$-functional on~$L^2(\mu)$. By the Riesz Representation Theorem for~$L^2(\mu)$, and by arbitrariness of~$\class[\mu]{f}\in L^2(\mu)$, we thus have~$\reptwo h\equiv 0$ $\mu$-a.e.

\smallskip

As a consequence of the claim,~$h=\ttclass[H]{\delta(\reptwo h)}=\rep\iota(\reptwo h)=\zero_H$. By arbitrariness of~$h\in (\im\,\iota)^\perp$, we may conclude that~$(\im\,\iota)^\perp=\set{\zero_H}$, i.e.\ that~$\iota$ is surjective.

\medskip

We show the first assertion in~\iref{i:p:DIntL^2:2}. A proof of the second assertion is similar, and therefore it is omitted.
Argue by contradiction that there exists~$f\in \mcL^2(\mu)$ with~$f\geq 0$ $\mu$-a.e., yet such that~$\tquadre{(\iota\class[\mu]{f})_\zeta}^-\neq \zero_{L^2(\mu_\zeta)}$ for all~$\zeta$ in some~$B\in\mcZ$ with~$\nu B>0$.
In particular, since~$\zeta\mapsto\tquadre{(\iota\class[\mu]{f})_\zeta}^-$ is a $\nu$-measurable field by Remark~\ref{r:Order}, the following integral is well-defined and strictly positive
\begin{align*}
\int_B \norm{\tquadre{(\iota\class[\mu]{f})_\zeta}^-}_{L^2(\mu_\zeta)}^2 \diff\nu(\zeta)>0\fstop
\end{align*}
Then, by~\eqref{eq:DisintF},
\begin{align*}
\int_X f^2\diff\mu =& \int_X (f^+)^2\diff\mu = \int_Z \int_X (f^+)^2 \diff\mu_\zeta \diff\nu(\zeta)=\int_Z \norm{\tclass[\mu_\zeta]{f^+}}_{L^2(\mu_\zeta)}^2 \diff\nu(\zeta)
\\
<& \int_Z \norm{\tclass[\mu_\zeta]{f^+}}_{L^2(\mu_\zeta)}^2 \diff\nu(\zeta)
+ \int_B \norm{\tquadre{(\iota\class[\mu]{f})_\zeta}^-}_{L^2(\mu_\zeta)}^2 \diff\nu(\zeta) \fstop
\end{align*}
Since~$\class[\mu_\zeta]{f^+}=\class[\mu_\zeta]{f}^+$ for every~$\zeta\in Z$, and by definition of~$\iota$, continuing from the previous inequality, we have that
\begin{align*}
\int_X f^2\diff\mu <& \int_Z \norm{\tclass[\mu_\zeta]{f^+}}_{L^2(\mu_\zeta)}^2 \diff\nu(\zeta)
+ \int_B \norm{\tquadre{(\iota\class[\mu]{f})_\zeta}^-}_{L^2(\mu_\zeta)}^2 \diff\nu(\zeta)
\\
=& \int_Z \norm{\tquadre{\class[\mu_\zeta]{f}}^+}_{L^2(\mu_\zeta)}^2 \diff\nu(\zeta)
+ \int_B \norm{\tquadre{\class[\mu_\zeta]{f}}^-}_{L^2(\mu_\zeta)}^2 \diff\nu(\zeta)
\\
\leq & \int_Z \tonde{\norm{\tquadre{\class[\mu_\zeta]{f}}^+}_{L^2(\mu_\zeta)}^2+\norm{\tquadre{\class[\mu_\zeta]{f}}^-}_{L^2(\mu_\zeta)}^2} \diff\nu(\zeta)
\\
= & \int_Z \tnorm{\class[\mu_\zeta]{f}}_{L^2(\mu_\zeta)}^2\diff\nu(\zeta)
\\
=&\ \tnorm{\iota\class[\mu]{f}}_H^2
\end{align*}
by definition of~$H$.
The inequality contradicts the fact, shown in~\iref{i:p:DIntL^2:1}, that~$\iota\colon L^2(\mu)\to H$ is an isometry, and therefore~$\tnorm{\class[\mu]{f}}_{L^2(\mu)}^2=\tnorm{\iota\class[\mu]{f}}_H^2$.
\end{proof}

\subsection{Direct integrals of Dirichlet forms} Let~$(X,\mcX,\mu)$ be $\sigma$-finite standard,~$(Z,\mcZ,\nu)$ be $\sigma$-finite countably generated, and~$\seq{\mu_\zeta}_{\zeta\in Z}$ be a pseudo-disintegration of~$\mu$ over~$\nu$. Further let~$\zeta\mapsto(Q_\zeta, \dom{Q_\zeta})$ be a $\nu$-measurable field of quadratic forms, each densely defined in~$L^2(\mu_\zeta)$ with separable domain, and denote by~$(Q,\dom{Q})$ their direct integral in the sense of Definition~\ref{d:assQ}.

\begin{definition}\label{d:Compat}
We say that~$(Q,\dom{Q})$ is \emph{compatible} with the pseudo-disintegra\-tion $\seq{\mu_\zeta}_{\zeta\in Z}$ if the space~$S_Q$ underlying~$\zeta\mapsto \dom{Q_\zeta}_1$ is of the form~$S_\mcA$ for some~$\mcA$ as in~\eqref{eq:AssA} and additionally satisfying~$\mcA\subset \dom{Q}$.
\end{definition}

Note that, if~$S_Q$ is of the form~$S_\mcA$ for~$\mcA\subset\dom{Q}$ and satisfying~\eqref{eq:AssA}, then~$S_H$ is of the form~$S_\mcA$ as well by Remark~\ref{r:S}.

\begin{definition}[Diagonal restriction]\label{d:DRestr}
Let~$(Q,\dom{Q})$ be a direct integral of quadratic forms compatible with a pseudo-disintegration~$\seq{\mu_\zeta}_{\zeta\in Z}$. The form
\begin{align*}
Q_\res=Q\eqdef \dint[S]{Z} Q_\zeta \diff\nu(\zeta)\comma \qquad \dom{Q_\res}\eqdef \dom{Q}\cap \iota(L^2(\mu))
\end{align*}
is a closed (densely defined) quadratic form on~$\iota(L^2(\mu))$, called the \emph{diagonal restriction} of~$(Q,\dom{Q})$.
\end{definition}

\begin{remark}[Comparison with{~\cite{AlbRoe90}}]\label{r:AlbRoe90:2} We note that the form~$(Q_\res,\dom{Q_\res})$ coincides with the form $(\msE,D(\msE))$ defined in~\cite[Thm.~1.2]{AlbRoe90}. As a consequence, at least in this case, the closability of~$\msE$ in~\cite[Thm.~1.2]{AlbRoe90} follows from our Proposition~\ref{p:DirInt}.
\end{remark}

Our first result on direct integrals of concrete quadratic forms is as follows.

\begin{proposition}\label{p:DirIntE}
Let~$(X,\mcX,\mu)$ be $\sigma$-finite standard,~$(Z,\mcZ,\nu)$ be $\sigma$-finite countably generated, and~$\seq{\mu_\zeta}_{\zeta\in Z}$ be a separated pseudo-disintegration of~$\mu$ over~$\nu$. Further let~$(E,\dom{E})$ be a direct integral of closed quadratic forms~$\zeta\mapsto(E_\zeta, \dom{E_\zeta})$ compatible with~$\seq{\mu_\zeta}_{\zeta\in Z}$. Then, $(E,\dom{E})$ is a Dirichlet form on~$L^2(\mu)$ if and only if~$(E_\zeta, \dom{E_\zeta})$ is so on~$L^2(\mu_\zeta)$ for $\nu$-a.e.~$\zeta\in Z$.
\end{proposition}

\begin{remark}
Since~$\seq{\mu_\zeta}_{\zeta\in Z}$ is separated, the isometry~$\iota$ in Proposition~\ref{p:DIntL^2} is a unitary operator, thus there exists~$\iota^{-1}\colon H\to L^2(\mu)$, where~$H$ is as in~\eqref{eq:DirIntL^2:0}.
For the sake of clarity, \emph{only} in the proof of Proposition~\ref{p:DirIntE} below, we distinguish between the quadratic form~$(E,\dom{E})$ on~$H$ and the quadratic form~$(\iota^*E,\dom{\iota^* E})$ on~$L^2(\mu)$ defined by
\begin{align*}
\iota^* E(f,g)\eqdef E(\iota f,\iota g)\comma \qquad f,g\in\dom{\iota^*E}\eqdef \iota^{-1}(\dom{E}) \fstop
\end{align*}
By definition of~$\iota^*E$ and since~$\iota\colon L^2(\mu)\to H$ is unitary, we have that~$(\iota^* E)_1=\iota^*E_1$. Therefore, $(\iota^*E,\dom{\iota^* E})$ is closed on~$L^2(\mu)$, since~$(E,\dom{E})$ is closed on~$H$ by Proposition~\ref{p:DirInt}\iref{i:p:DirInt1}.
Furthermore, the Hilbert spaces~$\dom{E}_1$ and~$\dom{\iota^*E}_1$ are intertwined via the unitary isomorphism~$\iota$.
In the statement of Proposition~\ref{p:DirIntE} above and everywhere after its proof ---~with a slight abuse of notation~--- these two quadratic forms are identified.
Again for the sake of clarity, the statement of the proposition equivalently reads as follows: $(\iota^*E,\dom{\iota^* E})$ is a Dirichlet form on~$L^2(\mu)$ if and only if~$(E_\zeta,\dom{E_\zeta})$ is a Dirichlet form on~$L^2(\mu_\zeta)$ for $\nu$-a.e.~$\zeta\in Z$. 
\end{remark}

\begin{proof}[Proof of Proposition~\ref{p:DirIntE}]
By e.g.~\cite[Thm.~1.4.1]{FukOshTak11}, the \emph{closed} quadratic form $(\iota^*E,\dom{\iota^*E})$ on~$L^2(\mu)$, resp.~$(E_\zeta,\dom{E_\zeta})$ on~$L^2(\mu_\zeta)$, is a Dirichlet form if and only if the associated semigroup $T^\iota_\bullet\colon L^2(\mu)\to L^2(\mu)$, resp.~$T_{\zeta,\iota}\colon L^2(\mu_\zeta)\to L^2(\mu_\zeta)$, is sub-Markovian, viz.
\begin{subequations}
\begin{align}
\label{eq:p:DirIntE:1a}
0\leq T^\iota_t u\leq&\ 1 \as{\mu}\comma & u\in L^2(\mu) : 0\leq&\ u\leq 1 \as{\mu}\comma && t>0\comma
\\
\label{eq:p:DirIntE:1b}
\text{resp.} \quad 0\leq T_{\zeta,t} v_\zeta\leq&\ 1 \as{\mu_\zeta}\comma & v_\zeta\in L^2(\mu_\zeta) : 0\leq&\ v_\zeta\leq 1 \as{\mu_\zeta}\comma && t>0\fstop
\end{align}
\end{subequations}
Thus, it suffices to show that~$T^\iota_\bullet$ is sub-Markovian if and only if~$T_{\zeta,\bullet}$ is so for $\nu$-a.e.~$\zeta\in Z$.

Since~$(E,\dom{E})$ and~$(\iota^*E,\dom{\iota^*E})$ are intertwined by the unitary isomorphism~$\iota$, it is not difficult to show that their semigroups~$T_\bullet$ and~$T^\iota_\bullet$ are intertwined as well, viz.
\begin{align*}
T^\iota_t=\iota^*T_t\eqdef T_t\circ \iota\comma \qquad t>0\fstop
\end{align*}
Furthermore, since~$\iota\colon L^2(\mu)\to H$ is a Riesz homomorphism by Proposition~\ref{p:DIntL^2}\iref{i:p:DIntL^2:2}, 
\begin{align}\label{eq:p:DirIntE:2}
T^\iota_\bullet \text{ is sub-Markovian } \iff \zero_H\leq T_t h\leq \uno_H\comma \quad h\in H : \zero_H\leq h \leq \uno_H \comma
\end{align}
where, by definition of the Banach lattice structure on~$H$ in Remark~\ref{r:Order},
\begin{align}\label{eq:p:DirIntE:3}
\zero_H\leq T_t h\leq \uno_H \iff \begin{aligned} &0\leq (T_t h)_\zeta\leq 1 \as{\mu_\zeta}\comma \quad t>0 \\ &\forallae{\nu} \zeta\in Z\comma \quad h\in H : \zero_H\leq h \leq \uno_H\fstop \end{aligned}
\end{align}
By Proposition~\ref{p:DirInt}\iref{i:p:DirInt3} we have that~$(T_t h)_\zeta= T_{\zeta,t} h_\zeta$ for $\nu$-a.e.~$\zeta\in Z$ for every~$h\in H$.
Thus, combining~\eqref{eq:p:DirIntE:2} and~\eqref{eq:p:DirIntE:3} we may conclude that
\begin{align}\label{eq:p:DirIntE:4}
T^\iota_\bullet \text{ is sub-Markovian } \iff \begin{aligned} &0 \leq T_{\zeta,t} h_\zeta \leq 1 \as{\mu_\zeta} \comma \quad t>0 \\ & \forallae{\nu} \zeta\in Z\comma \ h\in H : \zero_H\leq h \leq \uno_H \fstop\end{aligned}
\end{align}

The reverse implication in~\eqref{eq:p:DirIntE:4} together with~\eqref{eq:p:DirIntE:1b} immediately show that, if~$T_{\zeta,\bullet}$ is sub-Markovian for $\nu$-a.e.~$\zeta\in Z$, then~$T^\iota_\bullet$ is sub-Markovian.
The converse implication is not immediate, since the right-hand side of~\eqref{eq:p:DirIntE:4}, to be compared with~\eqref{eq:p:DirIntE:1b}, contains the additional consistency constraint that~$\zeta\mapsto h_\zeta$ be a measurable field representing an element~$h\in H$.

In order to show that, if~$T^\iota_\bullet$ is sub-Markovian, then~$T_{\zeta,\bullet}$ is sub-Markovian for $\nu$-a.e.~$\zeta\in Z$ is sub-Markovian, we argue as follows.
Since~$T_{\zeta,t}\colon L^2(\mu_\zeta)\to L^2(\mu_\zeta)$ is bounded, and since the unit contraction operator~$v\mapsto 0\vee v \wedge 1$ operates continuously on~$L^2(\mu_\zeta)$, it suffices to show that, for any sequence~$\seq{v^\zeta_n}_n\subset L^2(\mu_\zeta)$ total in~$L^2(\mu_\zeta)$,
\begin{align*}
0\leq T_{\zeta, t} (0\vee v^\zeta_n \wedge 1) \leq 1 \comma \qquad n\in \N \fstop
\end{align*}

Let~$\seq{u_n}_n\subset H$ be a fundamental sequence of $\nu$-measurable vector fields for~$H$, and recall that~$\seq{u_{n,\zeta}}_n$ is total in~$L^2(\mu_\zeta)$ for every~$\zeta\in Z$ by  Definition~\ref{d:DirInt}\iref{i:d:DirInt3}.
Applying~\eqref{eq:p:DirIntE:4} to each element~$u_n$ of this sequence proves the assertion.
\end{proof}

Proposition~\ref{p:DirIntE} motivates the next definition. A simple example follows.

\begin{definition}\label{d:DirIntE}
A quadratic form~$(E,\dom{E})$ on~$L^2(\mu)$ is a \emph{direct integral of Dirichlet forms}~$\zeta\mapsto (E_\zeta,\dom{E_\zeta})$ on~$L^2(\mu_\zeta)$ if it is a direct integral of the Dirichlet forms~$\zeta\mapsto (E_\zeta,\dom{E_\zeta})$ additionally compatible with the \emph{separated} pseudo-disintegration~$\seq{\mu_\zeta}_\zeta$ in the sense of Definition~\ref{d:Compat}.
\end{definition}

\begin{example}\label{ese:R2}
Let~$X=\R^2$ with standard topology, Borel $\sigma$-algebra, and the $2$-dimensional Lebesgue measure~$\Leb^2$. Consider a Dirichlet form measuring energy only in the first coordinate, viz.
\begin{align*}
E(f)\eqdef \int_{\R^2} \abs{\partial_1 f(x_1,x_2)}^2 \diff\Leb^2(x_1,x_2)
\end{align*}
with~$f\in L^2(\R^2)$ and $f(\emparg, x_2)\in W^{1,2}(\R)$ for~$\Leb^1$-a.e.~$x_2\in \R$. Then,~$(E,\dom{E})$ is the direct integral~$x_2\mapsto \ttonde{E_{x_2}, W^{1,2}(\R)}$, where $x_2$ ranges in~$Z=\R$ the real line,~$(X,\mcX,\mu_{x_2})$ is again the standard real line for every~$x_2\in \R$, and
\begin{align*}
E_{x_2}(f)\eqdef \int_{\R} \abs{\diff f(\emparg, x_2)}^2 \diff\Leb^1\comma \qquad \forallae{\Leb^1} x_2\in \R \fstop
\end{align*}
\end{example}

\subsection{Superposition of Dirichlet forms}\label{ss:Superposition} We recall here the definition of a \emph{superposition of Dirichlet forms} in the sense of~\cite[\S{V.3.1}]{BouHir91}.
Let~$(Z,\mcZ,\nu)$ be $\sigma$-finite,~$(X,\mcX,\mu_\zeta)$ be $\sigma$-finite and~$(E_\zeta,\dom{E_\zeta})$ be a Dirichlet form on~$L^2(\mu_\zeta)$, for every~$\zeta\in Z$. Assume that
\begin{enumerate}[$(\mathsc{sp}_1)$]
\item\label{i:Superpos:1} $\zeta\mapsto \mu_\zeta A$ is $\nu$-measurable for every~$A\in\mcX$;
\item\label{i:Superpos:2} $\zeta\mapsto E_\zeta(\rep f)\in [0,\infty]$ is $\nu$-measurable for every measurable~$\rep f\colon X\rar [-\infty,\infty]$.
\end{enumerate}

Let us now consider 
\begin{itemize}
\item a measure~$\lambda\ll \nu$ on~$(Z,\mcZ)$ so that~$\mu\eqdef \int_Z \mu_\zeta\diff\lambda(\zeta)$ is $\sigma$-finite on~$(X,\mcX)$ (therefore: $(X,\mcX,\mu)$ is standard);
\item the subspace~$\mbbD$ of all functions~$f\in L^2(\mu)$ so that
\begin{align}\label{eq:FormDirInt}
\mcE(f)\eqdef \int_Z E_\zeta(f)\diff\nu(\zeta)<\infty \comma
\end{align}
\end{itemize}
and let us further assume that
\begin{enumerate}[$(\mathsc{sp}_1)$]\setcounter{enumi}{2}
\item\label{i:Superpos:3} $\mbbD$ is dense in~$L^2(\mu)$.
\end{enumerate}
Then, it is claimed in~\cite[p.~214]{BouHir91} that~\eqref{eq:FormDirInt} is well-defined and depends only on the $\mu$-class of~$f$, and it is shown in~\cite[Prop.~V.3.1.1]{BouHir91} that
\begin{definition}
$(\mcE,\mbbD)$ is a Dirichlet form on~$L^2(\mu)$, called the \emph{superposition of~$\zeta\mapsto E_\zeta$}.
\end{definition}
Note that we may always choose~$\lambda=\nu$ provided that the integral measure~$\mu$ defined above is $\sigma$-finite. If this is not the case, we may recast the definition by letting~$\nu\eqdef \lambda$. In this way, we may always assume with no loss of generality that~$\mu$ is given, and that~$\seq{\mu_\zeta}$ is a pseudo-disintegration of~$\mu$ over~$\nu$.

\begin{remark}
In fact,~\cite{BouHir91} requires all functions in~\ref{i:Superpos:1}-\ref{i:Superpos:2} to be $\mcZ$-measurable, rather than only $\nu$-measurable. Here, we relax this condition to `$\nu$-measurability' in order to simplify the proof of the reverse implication in the next Proposition~\ref{p:Superpos}. Our definition of `superposition' is equivalent to the one in~\cite{BouHir91}.
\end{remark}

\begin{proposition}\label{p:Superpos}
Let~$(X,\mcX,\mu)$ be $\sigma$-finite standard,~$(Z,\mcZ,\nu)$ be $\sigma$-finite countably generated, $\seq{\mu_\zeta}_{\zeta\in Z}$ be a \emph{separated} pseudo-dis\-in\-te\-gra\-tion of $\mu$ over~$\nu$, and $(E_\zeta,\dom{E_\zeta})$ be a Dirichlet form on~$L^2(\mu_\zeta)$ for every~$\zeta\in Z$.
Then, the following are equivalent:
\begin{enumerate}[$(i)$]
\item\label{i:p:Superpos:3} there exists the superposition~$(\mcE, \mbbD)$ of~$\zeta\mapsto E_\zeta$ and the space~$\mbbD$ in~\eqref{eq:FormDirInt} is $\mcE_1^{1/2}$-separable;
\item\label{i:p:Superpos:4} there exists a direct integral of Dirichlet forms~$(E,\dom{E})$ of the forms~$\zeta\mapsto E_\zeta$.
\end{enumerate}
Furthermore, if either one holds, then~$(E,\dom{E})$ and~$(\mcE, \mbbD)$ are isomorphic Dirichlet spaces.
\end{proposition}
\begin{proof} We only show that~\iref{i:p:Superpos:3} implies~\iref{i:p:Superpos:4}. A proof of the reverse implication is similar, and it is therefore omitted.
For simplicity, set throughout the proof~$H_\zeta\eqdef L^2(\mu_\zeta)$, with norm~$\norm{\emparg}_\zeta$, for every~$\zeta\in Z$.

\smallskip

Assume~\iref{i:p:Superpos:3}. It follows from~\iref{i:Superpos:1} that~$\zeta\mapsto \norm{f}_\zeta$ is $\nu$-measurable for every~$f\in \mbbD$, and thus from~\iref{i:Superpos:2} that~$\zeta\mapsto E_{\zeta,1}(f,g)$ is $\nu$-measurable for every~$f,g\in \mbbD$ by polarization.
By $\mcE_1^{1/2}$-separability of~$\mbbD$, there exists a countable $\Q$-linear space~$\mcU\subset \mcL^2(\mu)$ so that~$\class[\mu]{\mcU}$ is $\mcE_1^{1/2}$-dense in~$\mbbD$, and dense in~$L^2(\mu)$ by~\iref{i:Superpos:3}. Since~$\seq{\mu_\zeta}_{\zeta\in Z}$ is separated by assumption, it follows by Proposition~\ref{p:DIntL^2} that $\class[\mu_\zeta]{\mcU}$ is dense in~$L^2(\mu_\zeta)$ for $\nu$-a.e.~$\zeta\in Z$.
As a consequence, the quadratic form~$(E_{\zeta}, \class[\mu_\zeta]{\mcU})$ is densely defined on~$H_\zeta$. Since~$(E_\zeta,\dom{E_\zeta})$ is closed, the closure~$(E_\zeta^\res,\dom{E_\zeta^\res})$ of~$(E_{\zeta}, \class[\mu_\zeta]{\mcU})$ is well-defined and a Dirichlet form on~$H_\zeta= L^2(\mu_\zeta)$.

Again since~$\seq{\mu_\zeta}_{\zeta\in Z}$ is separated, we may then construct a form~$(E,\dom{E})$ as the direct integral of Dirichlet forms (Dfn.~\ref{d:DirIntE}) of the forms~$\zeta\mapsto (E_\zeta^\res,\dom{E_\zeta^\res})$ with underlying space of measurable vector fields~$S=S_\mcU$ generated by~$\delta(\mcU)$ in the sense of Proposition~\ref{p:Dix}. By construction, the pre-Hilbert spaces~$\ttonde{\class[\mu]{\mcU}, \mcE_1^{1/2}}$ and~$\ttonde{\class[\dom{E}_1]{\delta(\mcU)}, E_1^{1/2}}$ are linearly and latticially isometrically isomorphic. The isomorphism extends to a unitary lattice isomorphism between~$\ttonde{\mbbD, \mcE_1^{1/2}}$ and~$\dom{E}_1$. The last assertion follows provided we show the following claim.

\smallskip

\emph{Claim:~$\dom{E_\zeta}=\dom{E_\zeta^\res}$ for $\nu$-a.e.~$\zeta\in Z$.} Argue by contradiction that there exists~$B\in\mcZ^\nu$, with~$\nu B>0$ and so that~$\dom{E_\zeta^\res}\subsetneq \dom{E_\zeta}$ for every~$\zeta\in B$. We may assume with no loss of generality that~$B\in \mcZ$. Furthermore, since~$(Z,\mcZ,\nu)$ is $\sigma$-finite countably generated, we may and shall assume that~$\nu B<\infty$. Denote by~$\dom{E^\res_\zeta}_1^{\perp_\zeta}$ the $(E_\zeta)^{1/2}_1$-orthogonal complement of~$\dom{E^\res_\zeta}_1$ in~$\dom{E_\zeta}_1$.
By the axiom of choice, there exists~$\rep h\eqdef (\zeta\mapsto \rep h_\zeta)$ with~$\rep h_\zeta\in \dom{E^\res_\zeta}_1^{\perp_\zeta}$ and~$\ttnorm{\rep h_\zeta}_{\dom{E_\zeta}_1}=1$ for all~$\zeta\in B$ and~$\rep h_\zeta\eqdef \zero_{\dom{E_\zeta}}$ for~$\zeta\in B^\complement$.
By closability of~$(E_\zeta,\dom{E_\zeta})$, the domain~$\dom{E_\zeta}_1$ embeds identically via~$\iota_{\zeta,1}$ into~$H_\zeta$, thus~$\rep h_\zeta$ may be regarded as an element~$\iota_{\zeta,1} \rep h_\zeta$ of~$H_\zeta$ for every~$\zeta\in Z$. Since~$\norm{\iota_{\zeta,1}}_\op\leq 1$ for every~$\zeta\in Z$, then
\begin{align}\label{eq:p:Superpos:1}
0<\ttnorm{\iota_{\zeta,1} \rep h_\zeta}_\zeta \leq 1 \comma \qquad \zeta\in B\fstop
\end{align}
Since~$(E,\dom{E})$ is in particular a direct integral of quadratic forms,~$\zeta\mapsto \iota_{\zeta,1}$ is a $\nu$-measurable field of bounded operators, and thus~$\zeta\mapsto\iota_{\zeta,1}\rep h_\zeta$ is a $\nu$-measurable field of vectors.

Now, let~$H$ be defined as in~\eqref{eq:H}. Setting~$\bar h\eqdef (\zeta\mapsto \iota_{\zeta,1}\rep h_\zeta)$, it follows by~\eqref{eq:p:Superpos:1} that
\begin{align*}
\ttnorm{\bar h}^2\eqdef \int_Z \ttnorm{\iota_{\zeta,1} \rep h_\zeta}_\zeta^2\diff\nu(\zeta) =  \int_B \ttnorm{\iota_{\zeta,1} \rep h_\zeta}_\zeta^2\diff\nu(\zeta)\in (0, \nu B] \fstop
\end{align*}
In particular, for the equivalence class~$h\eqdef \ttclass[H]{\bar h}$, we have that~$h\neq \zero_H$. By Proposition~\ref{p:DIntL^2}, there exists~$\reptwo h\in \mcL^2(\mu)$ representing~$h\in H$, and thus satisfying~$\zero_{\mbbD}\neq \ttclass[\mu]{\reptwo h}\in \mbbD$. On the other hand though,
\begin{equation*}
\begin{aligned}
\mcE_1(\ttclass[\mu]{\reptwo h},\ttclass[\mu]{\reptwo u})\eqdef& \int_Z E_{\zeta,1}\ttonde{\ttclass[\mu_\zeta]{\reptwo h},\class[\mu_\zeta]{\reptwo u}}\diff\nu(\zeta)
\\
=& \int_Z E_{\zeta,1}\ttonde{\rep h_\zeta,\delta(\reptwo u)_\zeta}\diff\nu(\zeta)=0\comma
\end{aligned} \qquad \reptwo u\in \mcU\comma
\end{equation*}
by definition of~$\rep h$. By $\mcE_1^{1/2}$-density of~$\class[\mu]{\mcU}$ in~$\mbbD$, the latter implies that~$h=\zero_\mbbD$, the desired contradiction.
\end{proof}

\begin{remark} If the disintegration in Proposition~\ref{p:Superpos} is not separated, then~$(\mcE,\mbbD)$ is still isomorphic, as a quadratic form, to the diagonal restriction~$(Q_\res,\dom{Q_\res})$ (Dfn.~\ref{d:DRestr}) of the direct integral of quadratic forms~$(E,\dom{E})$.
\end{remark}

\section{Ergodic decomposition} 
Everywhere in this section, let~$(X,\T,\mcX,\mu)$ be satisfying Assumption~\ref{ass:Main}.
We are interested in the notion of invariant sets for a Dirichlet form.

\begin{definition}[Invariant sets, irreducibility, {cf.~\cite[p.~53]{FukOshTak11}}]\label{d:Invariant}
Let $(E,\dom{E})$ be a Dirichlet form on~$L^2(\mu)$. We say that~$A\subset X$ is $E$-\emph{invariant} if it is $\mu$-measurable and any of the following equivalent\footnote{See~\cite[Lem.~1.6.1, p.~53]{FukOshTak11}, the proof of which adapts \emph{verbatim} to our more general setting.} conditions holds.
\begin{enumerate}[$(a)$]
\item\label{i:d:Invariant:1} $T_t(\car_A f)=\car_A T_t f$ $\mu$-a.e. for any~$f\in L^2(\mu)$ and~$t>0$;
\item $T_t(\car_A f)=0$ $\mu$-a.e. on $A^\complement$ for any~$f\in L^2(\mu)$ and~$t>0$;
\item $G_\alpha(\car_A f)=0$ $\mu$-a.e. on $A^\complement$ for any~$f\in L^2(\mu)$ and~$\alpha>0$;
\item\label{i:d:Invariant:4} $\car_A f\in \dom{E}$ for any~$f\in \dom{E}$ and
\begin{align}\label{eq:d:Invariant1}
E(f,g)=E(\car_A f,\car_A g)+E(\car_{A^\complement}f,\car_{A^\complement} g)\comma \qquad f,g\in \dom{E}
\end{align}

\item $\car_A f\in \domext{E}$ for any~$f\in \domext{E}$ and~\eqref{eq:d:Invariant1} holds for any~$f,g\in \domext{E}$.
\end{enumerate}

The form~$(E,\dom{E})$ is \emph{irreducible} if, whenever~$A$ is $E$-invariant, then either~$\mu A=0$ or~$\mu A^\complement=0$.
\end{definition}

As shown by Example~\ref{ese:R2}, the form~$(E,\dom{E})$ constructed in Proposition~\ref{p:DirIntE} is hardly ever irreducible, even if~$(E_\zeta,\dom{E_\zeta})$ is so for every~$\zeta\in Z$.

\subsection{The algebra of invariant sets}\label{ss:AlgInvSets}
Invariants sets of symmetric Markov processes on locally compact Polish spaces are studied in detail by H.~\^Okura in~\cite{Oku92}. In particular, he notes the following. For~$A\in\mcX$ set
\begin{align*}
\class[E]{A}\eqdef \set{\tilde A\in\mcX^\mu : \car_{\tilde A} \text{~is an $E$-quasi-continuous version of~} \car_A}\fstop
\end{align*}
For arbitrary~$A\in\mcX$ it can happen that~$\class[E]{A}=\emp$ or that~$A\not\in\class[E]{A}$. If however~$(E,\dom{E})$ is regular, then~$\class[E]{A}$ is non-empty for every $E$-invariant set~$A$. Suppose now~$A_0$,~$A_1\in\mcX$ and~$\class[E]{A_0}\neq\emp$. Then, one has the following dichotomy, \cite[Rmk.~1.1(ii)]{Oku92},
\begin{itemize}
\item $\class[E]{A_0}=\class[E]{A_1}$ if (and only if)~$\mu(A_0\triangle A_1)=0$;
\item $\class[E]{A_0}\cap \class[E]{A_1}=\emp$ if (and only if)~$\mu(A_0\triangle A_1)> 0$. 
\end{itemize}
As a consequence, when describing an $E$-invariant set~$A$ of a regular Dirichlet form~$(E,\dom{E})$, we may use interchangeably the $E$-class~$\class[E]{A}$ ---~i.e.\ the finest object representing~$A$, as far as $E$ is concerned~--- and the $\mu$-class $\class[\mu]{A}$ representing~$A$ in the measure algebra of~$(X,\mcX,\mu)$. This motivates to allow~$A$ in our definition of invariant set to be $\mu$-measurable, rather than only measurable.

We turn now to the study of invariant sets via direct integrals. We aim to show that, under suitable assumptions on~$\mu$, a Dirichlet form~$(E,\dom{E})$ on~$L^2(\mu)$ may be decomposed as a direct integral~$\zeta\mapsto (E_\zeta,\dom{E_\zeta})$ with~$(E_\zeta,\dom{E_\zeta})$ irreducible for every~$\zeta\in Z$.
To this end, we need to construct a \emph{measure} space~$(Z,\mcZ,\nu)$ ``indexing'' $E$-invariant sets. Let us start with a heuristic argument, showing how this cannot be done na\"ively, at least in the general case when~$(X,\mcX,\mu)$ is merely $\sigma$-finite.

Let~$\mcX_0$ be the family of $\mu$-measurable $E$-invariant subsets of~$X$, and note that~$\mcX_0$ is a $\sigma$-sub-algebra of~$\mcX^\mu$, e.g.~\cite[Lem.~1.6.1, p.~53]{FukOshTak11}. Let~$\mu_0$ be the restriction of~$\hat\mu$ to~$(X,\mcX_0)$.
The space~$(X,\mcX_0,\mu_0)$ ---~our candidate for~$(Z,\mcZ,\nu)$~--- is generally \emph{not} $\sigma$-finite, nor even semi-finite. For instance, in the extreme case when~$(E,\dom{E})$ is irreducible and~$\mu X=\infty$, then~$\mcX_0$ is the minimal $\sigma$-algebra on~$X$, the latter is an atom, and thus~$\mu_0$ is purely infinite.
Since~$(X,\mcX,\mu)$ is $\sigma$-finite, every disjoint family of $\mu$-measurable non-negligible subsets is at most countable~\cite[215B(iii)]{Fre00}, thus~$(X,\mcX_0,\mu_0)$ has up to countably many disjoint atoms.
However, even in the case when~$(X,\mcX_0,\mu_0)$ has \emph{no} atoms,~$\mu_0$ might again be purely infinite. This is the case of Example~\ref{ese:R2}, where~$\mcX_0=\set{\emp,\R}\otimes\Bo{\R}^{\Leb^1}$ is the product $\sigma$-algebra of the minimal $\sigma$-algebra on the first coordinate with the Lebesgue $\sigma$-algebra on the second coordinate, and where~$\mu_0$ coincides with the $\hat\mu$-measure of horizontal stripes.
This latter example shows that, again even when~$(X,\mcX_0,\mu_0)$ has no atoms, the complete locally determined version~\cite[213D]{Fre00} of~$(X,\mcX_0,\mu_0)$ is trivial. Thus, in this generality, there is no natural way to make~$(X,\mcX_0,\mu_0)$ into a more amenable measure space while retaining information on $E$-invariant sets.

\medskip

The situation improves as soon as~$(X,\mcX,\mu)$ is a probability space, in which case so is~$(X,\mcX_0,\mu_0)$. The reasons for this fact are better phrased in the language of von Neumann algebras.

\begin{remark}[Associated von Neumann algebras]\label{r:vonNeumann}
Denote by~$\msM$ the space~$L^\infty(\mu)$ regarded as the (commutative, unital) von Neumann algebra of multiplication operators in~$B(L^2(\mu))$. Then,~$L^\infty(\mu_0)$ is a (commutative) von Neumann sub-algebra of~$\msM$, denoted by~$\msM_0$. Two key observations are as follows:
\begin{itemize}
\item since~$(X,\mcX_0,\mu_0)$ is now (semi-)finite,~$\msM_0$ is unital as well;

\item by Definition~\iref{d:Invariant}\iref{i:d:Invariant:4}, the algebra~$\msM_0$ acts by multiplication also on~$\dom{E}$, and the action~$\msM_0\acts L^2(\mu)$ is compatible with the action~$\msM_0\acts \dom{E}$ by restriction.
\end{itemize}
\end{remark}

The next definition, borrowed from~\cite{BiaCar09}, encodes a notion of ``smallness'' of the $\sigma$-algebra~$\mcX$ w.r.t.~$\mu$.

\begin{definition}[{\cite[Dfn.~A.1]{BiaCar09}}]
Let~$\mcX^*\subset \mcX$ be a countably generated $\sigma$-subalgebra. We say that:
\begin{itemize}
\item $\mcX$ is $\mu$-\emph{essentially countably generated by~$\mcX^*$} if for each~$A\in\mcX$ there is~$A^*\in\mcX^*$ with~$\mu(A\triangle A^*)=0$;
\item $\mcX$ is $\mu$-\emph{essentially countably generated} if it is so by some~$\mcX^*$ as above.
\end{itemize}
\end{definition}

By our Assumption~\ref{ass:Main},~$\mcX$ is countably generated, thus~$\mcX_0$ is $\mu_0$-essentially countably generated by $\mcX^*\eqdef \mcX\cap \mcX_0$. We denote by~$\mu_0^*$ the restriction of~$\mu_0$ to~$\mcX^*$.
As noted in~\cite[p.~418]{BiaCar09}, atoms of~$\mcX^*$ are, in general, larger (in cardinality, \emph{not} in measure) than atoms of~$\mcX$. It is therefore natural to pass to a suitable quotient space.
Following~\cite[Dfn.~A.5]{BiaCar09}, we define an equivalence relation~$\sim$ on~$X$ by
\begin{align}\label{eq:Quotient}
x_1\sim x_2 \quad\text{if and only if}\quad x_1\in A\iff x_2\in A \text{~for every~} A\in \mcX^*\fstop
\end{align}
Further let~$p\colon X\rar Z\eqdef X/\sim$ be the quotient map, $\mcZ\eqdef\set{B\subset Z : p^{-1}(B)\in\mcX^*}$ be the quotient $\sigma$-algebra induced by~$p$, and~$\nu\eqdef p_\pfwd \mu_0^*$ be the quotient measure. Similarly to~\cite[p.~416]{BiaCar09}, it follows by definition of~$\sim$ that every~$A\in\mcX^*$ is $p$-saturated. In particular:
\begin{align}\label{eq:BiaCar}
\emp\neq A\subset p^{-1}(p(x)) \implies A=p^{-1}(p(x)) \comma \qquad A\in\mcX^* \fstop
\end{align}

As a consequence~$\mcX^*$ and~$\mcZ$ are isomorphic and thus both are countably generated, since~$\mcX^*$ is by assumption. Furthermore,~$(Z,\mcZ)$ is separable by construction, and thus it is countably separated.

\subsection{Ergodic decomposition of forms: probability measure case}\label{ss:ErgProbab}
We are now ready to state our main result, a decomposition theorem for Dirichlet forms over their invariant sets.

\begin{theorem}[Ergodic decomposition: regular case]\label{t:Stone}
Let~$(X,\T,\mcX,\mu)$ be a locally compact Polish probability space, and~$(E,\dom{E})$ be a $\T$-regular Dirichlet form on $L^2(\mu)$. Then, there exist
\begin{enumerate}[$(i)$]
\item\label{i:t:Stone:1} a probability space~$(Z,\mcZ,\nu)$ and a measurable map~$s\colon (X,\mcX)\rar (Z,\mcZ)$;
\item\label{i:t:Stone:2} a $\nu$-essentially unique disintegration~$\seq{\mu_\zeta}_{\zeta\in Z}$ of~$\hat\mu$ w.r.t.~$\nu$, strongly consistent with~$s$, and so that, when~$s^{-1}(\zeta)$ is endowed with the subspace topology and the trace $\sigma$-algebra inherited by~$(X,\T,\mcX^\mu)$, then~$(s^{-1}(\zeta),\mu_\zeta)$ is a Radon probability space for~$\nu$-a.e.~$\zeta\in Z$;
\item\label{i:t:Stone:3} a $\nu$-measurable field~$\zeta\mapsto (E_\zeta,\dom{E_\zeta})$ of $\T$-regular irreducible Dirichlet forms $(E_\zeta,\dom{E_\zeta})$ on~$L^2(\mu_\zeta)$;
\end{enumerate}
so that
\begin{align}\label{eq:t:Ergodic:0}
L^2(\mu)=\dint[]{Z} L^2(\mu_\zeta)\diff\nu(\zeta) \qquad \text{and} \qquad E=\dint[]{Z} E_\zeta \diff\nu(\zeta) \fstop
\end{align}
\end{theorem}
\begin{proof} \iref{i:t:Stone:1} Let~$(Z,\mcZ,\nu)$ be the quotient space of~$(X,\mcX^*,\mu_0^*)$ defined in~\S\ref{ss:AlgInvSets}, and recall that~$(Z,\mcZ)$ is countably separated. Note that $\id_X\colon (X,\mcX^\mu,\hat\mu)\rar (X,\mcX_0^{\mu_0},\hat\mu_0)$ is inverse-measure-preserving~\cite[235H(b)]{Fre00}, thus so is
\begin{align}\label{eq:QuotientMap}
s\eqdef p\circ \id_X\colon (X,\mcX^\mu,\hat\mu)\rar (Z,\mcZ,\nu)\fstop
\end{align}

Since~$(X,\T,\mcX^\mu,\hat\mu)$ is Radon (in the sense of~\cite[411H(b)]{Fre00}) and~$(Z,\mcZ,\nu)$ is a probability space (in particular: strictly localizable~\cite[322C]{Fre00}), there exists a disintegration~$\seq{\mu_\zeta}_{\zeta\in Z}$ of~$\hat\mu$ over~$\nu$ consistent with~$s$, and so that~$(X,\T,\mcX,\mu_\zeta)$ is a Radon probability space~\cite[452O, 452G(a)]{Fre00}. 
Since~$(Z,\mcZ)$ is countably separated, $\seq{\mu_\zeta}_{\zeta\in Z}$ is in fact strongly consistent with~$s$~\cite[452G(c)]{Fre00}. By definition of strong consistency, we may restrict~$\mu_\zeta$ to~$s^{-1}(\zeta)$; the Radon property is preserved by this restriction~\cite[416R(b)]{Fre00}.
Since~$s$ factors through~$p$, one has~$s^{-1}(\zeta)\in\mcX^*\subset \mcX$ for every~$\zeta\in Z$. In particular,~$s^{-1}(\zeta)$ is a Borel subset of the metrizable Luzin space~$(X,\T)$, and thus a metrizable Luzin space itself by~\cite[\S{II.1}, Thm.~2, p.~95]{Sch73}. It follows that~$\supp[\mu_\zeta]$, endowed with the subspace topology inherited from~$(X,\T)$ and the induced Borel $\sigma$-algebra, satisfies Assumption~\ref{ass:Main}.
The disintegration is $\nu$-essentially unique, similarly to~\cite[Thm.~A.7, Step~1, p.~420]{BiaCar09}.

This shows~\iref{i:t:Stone:1}--\iref{i:t:Stone:2}. The proof of~\iref{i:t:Stone:3} is divided into several steps.

\paragraph{1. Measurable fields} Let~$\mcC$ be a special standard core~\cite[p.~6]{FukOshTak11} for~$(E,\dom{E})$, and~$N\subset Z$ be a $\nu$-negligible set so that~$(X,\T,\mcX,\mu_\zeta)$ is Radon for every~$\zeta\in N^\complement$. Then,~$\mcC\restr_{\supp[\mu_\zeta]}$ is dense in $L^2(\mu_\zeta)$ for every~$\zeta\in N^\complement$. In particular, since~$L^2(\mu)$ is separable, there exists a fundamental sequence~$\seq{u_n}_n\subset \mcC$, i.e.\ total in~$L^2(\mu_\zeta)$ for every~$\zeta\in N^\complement$, and additionally total in~$L^2(\mu)$.
Since~$(E,\dom{E})$ is regular, $\dom{E}_1$ is separable by~\cite[Prop.~IV.3.3(i)]{MaRoe92}, and therefore we can and will assume, with no loss of generality that~$\seq{u_n}_n$ is additionally $E^{1/2}_1$-total in~$\dom{E}_1$.
Moreover,~$\mcC$ is an algebra, thus~$\zeta\mapsto \scalar{f}{g}_\zeta=\mu_\zeta(fg)$ is $\nu$-measurable by definition of disintegration for every~$f$,~$g\in\mcC$. As a consequence, by Proposition~\ref{p:Dix} there exists a \emph{unique} $\nu$-measurable field of Hilbert spaces~$\zeta\mapsto L^2(\mu_\zeta)$ making $\nu$-measurable all functions of the form~$\zeta\mapsto \mu_\zeta f$ with $f\in\mcC$. We denote by~$S_\mcC$ the underlying space of $\nu$-measurable vector fields. Everywhere in the following, we identify~$\class[\mu_\zeta]{f}$ with a fixed continuous representative~$f\in\mcC$, thus writing~$f$ in place of~$\delta(f)$.

\paragraph{2. $L^2$-isomorphism} Since~$\seq{\mu_\zeta}_{\zeta\in Z}$ is strongly consistent with~$s$, it is separated. Therefore, the first isomorphism in~\eqref{eq:t:Ergodic:0} follows now by~\eqref{eq:DirIntL^2:0} with underlying space~$S_\mcC$. In the following, set~$H\eqdef L^2(\mu)$ and~$H_\zeta\eqdef L^2(\mu_\zeta)$.

\paragraph{3. Semigroups} Let~$T_t$ be the semigroup associated to~$(E,\dom{E})$ and consider the natural complexification~$T_t^\C$ of~$T_t$ defined on~$H_\C\eqdef H\otimes_\R \C$. For~$g\in L^\infty(\nu)$ denote by
\begin{align*}
M_g\eqdef \dint[S_\mcC]{Z} g(\zeta)\uno_{H_\zeta} \diff\nu(\zeta)
\end{align*}
the associated \emph{diagonalizable operator} in~$B(H)$,~\cite[\S{II.2.4} Dfn.~3, p.~185]{Dix81}. For~$B\in\mcZ$ set~$M_B\eqdef M_{\car_B}$.

\emph{Claim: the commutator~$[T_t^\C, M_g^\C]$ vanishes for~$g\in L^\infty_\C(\nu)$.} 
By~\cite[\S{II.2.3}, Prop.~4(ii), p.~183]{Dix81} and the norm-density of simple functions in~$L^\infty(\nu)$, it suffices to show that~$[T_t, M_B]=0$ for every~$B\in\mcZ$. 
To this end, recall the discussion~\cite[p.~165]{Dix81} on $\nu$-measurable structures induced by $\nu$-measurable subsets of~$Z$. Since~$B\in\mcZ$, then~$A\eqdef p^{-1}(B)\in \mcX_0^{\mu_0}$, and~$A\in\mcX^\mu$ as well~\cite[235H(c)]{Fre00}.
 Note further that, since~$H$ is reconstructed as a direct integral with underlying space~$S_\mcC$, for every~$h\in \mcC$ the representative~$h_\zeta$ of~$h$ in~$H_\zeta$ may be chosen so that~$h_\zeta=h$ for every~$\zeta\in Z$. Thus, for all~$f,g\in\mcC$,
\begin{align*}
\scalar{M_B f}{g}_H=& \int_Z \car_{B}(\zeta) \scalar{\uno_{H_\zeta} f_\zeta}{g_\zeta}_\zeta \diff\nu(\zeta)
\\
=& \int_B \int_X f g \diff \mu_\zeta \diff\nu(\zeta)=\int_A fg \diff\mu
\\
=& \scalar{\car_A f}{g}_H\fstop
\end{align*}
By density of~$\mcC$ in~$H$ and since~$M_B$ is bounded, it follows that~$M_B=\car_A$ as elements of~$B(H)$. Thus,~$[T_t,M_B]=[T_t, \car_A]=0$ for every~$t>0$ by Definition~\ref{d:Invariant}\iref{i:d:Invariant:1}, since~$A$ is $E$-invariant.

\smallskip

By the characterization of decomposable operators via diagonalizable operators~\cite[\S{II.2.5}, Thm.~1, p.~187]{Dix81},~$T_t$ is decomposable, and represented by a $\nu$-measurable field of contraction operators~$\zeta\mapsto T_{\zeta,t}$.
Finally, in light of~\cite[\S{II.2.3}, Prop.~4, p.~183]{Dix81}, it is a straightforward verification that~$T_{\zeta,t}$, $t>0$, is a strongly continuous symmetric contraction semigroup on~$H_\zeta$ for~$\nu$-a.e.~$\zeta\in Z$, since so is~$T_t$.
Analogously to the proof of Proposition~\ref{p:DirIntE},~$T_{\zeta,t}$ is sub-Markovian for~$\nu$-a.e.~$\zeta\in Z$, since so is~$T_t$.

\paragraph{4. Forms: construction} Denote by~$(E_\zeta,\dom{E_\zeta})$ the Dirichlet form on~$L^2(\mu_\zeta)$ associated to the sub-Markovian semigroup~$T_{\zeta,t}$ for~$\nu$-a.e.~$\zeta\in Z$. Let further~$G_{\zeta,\alpha}$,~$\alpha>0$, be the associated strongly continuous contraction resolvent.

We claim that~$\mcC\subset \dom{E_\zeta}$ for~$\nu$-a.e.~$\zeta\in Z$. 
Firstly, note that~$\zeta\mapsto E_\zeta(f,g)$ is~$\nu$-measurable, since it is the $\nu$-a.e.-limit of the measurable functions~$\zeta\mapsto E^{(\beta)}_\zeta(f,g)\eqdef \scalar{f-\beta G_{\zeta,\beta} f}{g}_\zeta$ as~$\beta\rar\infty$ by~\eqref{eq:Hille--YosidaE}. By~\cite[p.~27]{MaRoe92},
\begin{align*}
E_\zeta(f,g)=\lim_{\beta\rar\infty} \scalar{f-\beta \int_0^\infty e^{-\beta t} T_{\zeta,t} f \diff t}{g}_\zeta\comma \qquad f,g\in\mcC\fstop
\end{align*}

Now,
\begin{align*}
\int_Z E_\zeta(f) \diff\nu(\zeta)=&\int_Z\lim_{\beta\rar\infty} \scalar{f-\beta \int_0^\infty e^{-\beta t} T_{\zeta,t} f \diff t}{f}_\zeta \diff\nu(\zeta)
\\
\leq& \liminf_{\beta\rar\infty} \int_Z \scalar{f-\beta \int_0^\infty e^{-\beta t} T_{\zeta,t} f \diff t}{f}_\zeta \diff\nu(\zeta)
\end{align*}
by Fatou's Lemma. It is readily checked that, since~$\norm{T_{\zeta,t}}_{\op}\leq 1$, we may exchange the order of both integration and $H_\zeta$-scalar products by Fubini's Theorem. Thus,
\begin{align*}
\int_Z E_\zeta(f) \diff\nu(\zeta)\leq& \liminf_{\beta\rar\infty} \int_Z \scalar{f-\beta \int_0^\infty e^{-\beta t} T_{\zeta,t} f \diff t}{f}_\zeta \diff\nu(\zeta)
\\
=&\int_Z \norm{f}_\zeta^2 \diff\nu(\zeta)-\limsup_{\beta\rar\infty} \int_0^\infty \int_Z \beta e^{-\beta t}\scalar{T_{\zeta,t} f}{f}_\zeta \diff\nu(\zeta) \diff t\fstop
\end{align*}
By the representation of~$T_t$ via~$\zeta\mapsto T_{\zeta,t}$ established in \emph{Step~3},\begin{align*}
\int_Z E_\zeta(f) \diff\nu(\zeta)\leq&\int_Z \int_X f^2 \diff\mu_\zeta \diff\nu(\zeta)-\limsup_{\beta\rar\infty} \int_0^\infty \beta e^{-\beta t}\scalar{T_t f}{f}_H \diff t\fstop
\end{align*}
Finally, by~\eqref{eq:Disint:0},~\cite[p.~27]{MaRoe92} and~\eqref{eq:Hille--YosidaE},
\begin{align*}
\int_Z E_\zeta(f) \diff\nu(\zeta)\leq&\liminf_{\beta\rar\infty} \beta \int_0^\infty \scalar{f-e^{-\beta t}T_t f}{f}_H \diff t
\\
=&\liminf_{\beta\rar\infty} \beta\scalar{f-\int_0^\infty e^{-\beta t} T_t f  \diff t }{f}_H=E(f)<\infty \fstop
\end{align*}
This shows that~$E_\zeta(f)<\infty$ for every~$f\in\mcC$ for~$\nu$-a.e.~$\zeta\in Z$, thus~$\mcC\subset\dom{E_\zeta}$ $\nu$-a.e.

\emph{Claim:~$\mcC$ is a core for~$(E_\zeta,\dom{E_\zeta})$ for~$\nu$-a.e.~$\zeta\in Z$.}
It suffices to show that the inclusion~$\mcC\subset \dom{E_\zeta}$ is $(E_\zeta)^{1/2}_1$-dense for~$\nu$-a.e.~$\zeta\in Z$.
Argue by contradiction that there exists a $\nu$-measurable non-negligible set~$B$ so that the inclusion~$\mcC\subset \dom{E_\zeta}$ is not $(E_\zeta)^{1/2}_1$-dense for every~$\zeta\in B$, and let~$\mcC^{\perp_\zeta}$ be the $(E_\zeta)^{1/2}_1$-orthogonal complement of~$\mcC$ in~$\dom{E_\zeta}$. By the axiom of choice we may construct~$h\in \prod_{\zeta\in Z} H_\zeta$ so that~$h_\zeta\in \mcC^{\perp_\zeta} \setminus\set{0}$ for~$\zeta\in B$ and~$h_\zeta=0$ for~$\zeta\in B^\complement$. Further let~$\seq{u_n}_n\subset \mcC$ be as in \emph{Step~1}. Then,~$\zeta\mapsto \scalar{u_n}{h_\zeta}_\zeta=0$ is $\nu$-measurable for every~$n$. As a consequence,~$\zeta\mapsto h_\zeta$ is $\nu$-measurable (i.e., it belongs to~$S_\mcC$) by~\cite[\S{II.1.4}, Prop.~2, p.~166]{Dix81}. By the first isomorphism in~\eqref{eq:t:Ergodic:0}, established in \emph{Step~1},~$\zeta\mapsto h_\zeta$ represents an element~$h$ in~$H$. Since~$\seq{u_n}_n$ is total in~$H$, there exists~$n$ so that
\begin{align*}
0\neq \scalar{u_n}{h}_H=\int_Z\scalar{u_n}{h_\zeta}_\zeta\diff\nu(\zeta)=0 \comma
\end{align*}
a contradiction. Since functions in~$\mcC$ are continuous, the form~$(E_\zeta,\dom{E_\zeta})$ is regular for $\nu$-a.e.~$\zeta\in Z$. In particular, $\dom{E_\zeta}_1$ is separable for $\nu$-a.e.~$\zeta\in Z$ by~\cite[Prop.~IV.3.3]{MaRoe92}.

\smallskip

We note that, by the above claim and~\cite[Prop.~IV.3.3(i)]{MaRoe92},~$\dom{E_\zeta}_1$ is separable for every~$\zeta\in Z$, and so the observation in Remark~\ref{r:Separability} is satisfied.

\paragraph{5. Forms: direct integral} By \emph{Step~1}, resp.~\emph{Step~4},~$\zeta\mapsto L^2(\mu_\zeta)$, resp.~$\zeta\mapsto \dom{E_\zeta}_1$, is a $\nu$-measurable field of Hilbert spaces with underlying space~$S_\mcC$. In particular,~$\zeta\mapsto (E_\zeta,\dom{E_\zeta})$ satisfies Definition~\ref{d:assQ}, and we may consider the direct integral of quadratic forms
\begin{align}\label{eq:t:Stone:1}
\tilde E(f)\eqdef \dint[S_\mcC]{Z} E_\zeta(f) \diff\nu(\zeta)
\end{align}
defined by~\eqref{eq:DIntQF}.
We claim that~$(\tilde E, \dom{\tilde E})=(E,\dom{E})$. This is a consequence of Proposition~\ref{p:DirInt}\iref{i:p:DirInt3}, since~\eqref{eq:p:DirInt:0} was shown in \emph{Step 3} for~$T_t$. Definition~\ref{d:Compat} holds with~$\mcA=\mcC$ by construction.

\paragraph{6. Forms: irreducibility} Let~$A_\zeta$ be $E_\zeta$-invariant, with~$\hat\mu_\zeta A_\zeta>0$.
With no loss of generality, we may and will assume that~$A_\zeta\in \mcX$. Up to removing a $\nu$-negligible set of~$\zeta$'s, we have that~$A_\zeta\subset s^{-1}(\zeta)$, by strong consistency of the disintegration. Thus, by~\eqref{eq:Disint:0},
\begin{align}\label{eq:t:Main:1}
\mu A_\zeta=\mu \ttonde{A_\zeta\cap s^{-1}(\zeta)}=\int_{\set{\zeta}} \mu_\zeta A_\zeta \diff\nu(\zeta)=\mu_\zeta A_\zeta \cdot \nu\set{\zeta} \fstop
\end{align}

\emph{Claim:~$A_\zeta\in\mcX_0$.} Assume first that~$\nu\set{\zeta}=0$. Then,~$A_\zeta$ is contained in the $\mu$-negligible invariant set~$s^{-1}(\zeta)$, hence, it is $E$-invariant, i.e.~$A_\zeta\in\mcX_0$.
Assume now~$\nu\set{\zeta}>0$. By~\eqref{eq:t:Main:1},~$\mu A_\zeta>0$, thus~$A_\zeta\neq \emp$ and~$\car_{A_\zeta}\neq \zero_H$. By~\cite[\S{II.2.3}, Prop.~3, p.~182]{Dix81} and the direct-integral representation of~$T_t$ established in \emph{Step 3}, 
\begin{align*}
\car_{A_\zeta} T_t= \dint{Z} \car_{A_\zeta} T_{\zeta', t}\diff\nu(\zeta') \fstop
\end{align*}
By strong consistency, $A_\zeta$ is $\mu_{\zeta'}$-negligible for every~$\zeta'\neq \zeta$, thus in fact
\begin{align*}
\car_{A_\zeta} T_t= \nu\!\set{\zeta} \car_{A_\zeta} T_{\zeta, t} \comma
\end{align*}
whence, by $E_\zeta$-invariance of~$A_\zeta$,
\begin{align*}
\car_{A_\zeta} T_t= \nu\!\set{\zeta} \car_{A_\zeta} T_{\zeta, t}= \nu\!\set{\zeta}T_{\zeta, t} \car_{A_\zeta}= T_t \car_{A_\zeta} \comma
\end{align*}
and so~$A_\zeta$ is $E$-invariant, and thus~$A_\zeta\in\mcX_0$.

Now, since~$A_\zeta\in\mcX$ by assumption, then~$A_\zeta\in \mcX^*\eqdef \mcX\cap \mcX_0$. Together with~$A_\zeta\subset s^{-1}(\zeta)$, this implies that either~$A_\zeta=\emp$, or~$A_\zeta=s^{-1}(\zeta)$ by~\eqref{eq:BiaCar}. Thus, it must be~$A_\zeta=s^{-1}(\zeta)$, since~$\hat\mu_\zeta A_\zeta>0$ by assumption.
Since~$s^{-1}(\zeta)$ is~$\hat\mu_\zeta$-conegligible, this shows that~$(E_\zeta,\dom{E_\zeta})$ is irreducible.
\end{proof}

In the statement of Theorem~\ref{t:Stone}, we write that each~$(E_\zeta,\dom{E_\zeta})$ is a regular Dirichlet form on~$L^2(\mu_\zeta)$ with underlying space~$(X,\T,\mcX,\mu_\zeta)$ to emphasize that the topology of the space is the given one. As it is well-known however, in studying the potential-theoretic and probabilistic properties of a Dirichlet form~$(E,\dom{E})$ on~$L^2(\mu)$, one should always assume that~$\mu$ has full support, which is usually not the case for~$\mu_\zeta$ on~$(X,\T)$. In the present case, the restriction of~$\mu_\zeta$ to~$s^{-1}(\zeta)$ is however harmless, since~$s^{-1}(\zeta)$ is $E$-invariant, and therefore~$s^{-1}(\zeta)^\complement$ is also~$E_\zeta$-exceptional.

\begin{remark}\label{r:QuasiClopen}
As anticipated in \S\ref{s:Intro}, if~$(E,\dom{E})$ is regular and strongly local, then every invariant set admits an $E$-quasi-clopen $\mu$-modification~\cite[Cor.~4.6.3, p.~194]{FukOshTak11}. This suggests that, at least in the local case, one may treat $E$-invariant sets as ``connected components'' of~$X$.
Our intuition can be made rigorous by noting that $E$-invariant subsets of~$X$ are in bijective correspondence to compact open subsets of the spectrum~$\spec{\msM_0}$ of the von Neumann algebra~$\msM_0$ (cf.~Rmk.~\ref{r:vonNeumann}), endowed with its natural weak* topology. In particular,~$\spec{\msM_0}$ coincides with the Stone space of the measure algebra of~$(X,\mcX_0,\mu_0)$, and is thus a totally disconnected Hausdorff space. Its singletons correspond to the ``minimal connected components'' sought after in~\S\ref{s:Intro}.
At this point, we should emphasize that~$(Z,\mcZ)$ and~$\spec{\msM_0}$ are different measure spaces, the points of which index ``minimal invariant sets'' in~$X$. However, points in~$Z$ index ---~via~$s$~--- sets in~$\mcX^*$, whereas points in~$\spec{\msM_0}$ index sets in~$\mcX_0$. In this sense at least,~$Z$ is minimal with the property of indexing such ``minimal invariant sets'', while~$\spec{\msM_0}$ is maximal. For this reason, one might be tempted to use~$\spec{\msM_0}$ in place of~$(Z,\mcZ)$ in Theorem~\ref{t:Stone}. The issue is that~$\spec{\msM_0}$ is nearly always too large for the disintegration to be \emph{strongly} consistent with the indexing map.
\end{remark}

In the next result we show that the regularity of the Dirichlet form~$(E,\dom{E})$ in Theorem~\ref{t:Stone} may be relaxed to quasi-regularity. As usual, a proof of this result relies on the so-called transfer method.

\smallskip

Let~$(X,\T,\mcX,\mu)$ and~$(X^\sharp,\T^\sharp,\mcX^\sharp,\mu^\sharp)$ be measure spaces satisfying Assumption~\ref{ass:Main}.
We note \emph{en passant} that a Hilbert isomorphism of Dirichlet spaces $\iota\colon \dom{E}_1\rar \dom{E^\sharp}_1$, additionally preserving the $L^\infty$-norm on~$\dom{E}\cap L^\infty(\mu)$, is automatically a lattice isomorphism, e.g.~\cite[Lem.~A.4.1, p.~422]{FukOshTak11}.

\begin{theorem}[Ergodic decomposition: quasi-regular case]\label{t:QRegular}
The conclusions of Theorem~\ref{t:Stone} remain valid if~$(X,\T,\mcX,\mu)$ is a topological probability space satisfying Assumption~\ref{ass:Main} and ``\emph{regular}'' is replaced by ``\emph{quasi-regular}''.
\end{theorem}
\begin{proof} By the general result~\cite[Thm.~3.7]{CheMaRoe94}, there exist a locally compact Polish, Radon probability space $(X^\sharp, \T^\sharp, \mcX^\sharp, \mu^\sharp)$ and a quasi-homeomorphism 
\begin{equation*}
j\colon (X,\T,\mcX,\mu)\longrar (X^\sharp, \T^\sharp, \mcX^\sharp, \mu^\sharp)
\end{equation*}
so that~$(E,\dom{E})$ is quasi-homeomorphic, via~$j$, to a regular Dirichlet form $(E^\sharp,\dom{E^\sharp})$ on~$(X^\sharp, \T^\sharp, \mcX^\sharp, \mu^\sharp)$. 
Applying Theorem~\ref{t:Stone} to~$(E^\sharp,\dom{E^\sharp})$ gives a disintegration~$\tseq{\mu_\zeta^\sharp}_\zeta$ of~$\mu$ w.r.t.~$\nu$ and a direct-integral representation
\begin{align*}
E^\sharp=\dint{Z} E_\zeta^\sharp \diff\nu(\zeta) \comma
\end{align*}
where~$(E_\zeta^\sharp,\dom{E_\zeta^\sharp})$ is a regular Dirichlet form on~$L^2(\mu_\zeta^\sharp)$ for~$\nu$-a.e.~$\zeta\in Z$.

\paragraph{1. Forms}
In the following, whenever~$\seq{F_k}_k$ is a nest, let us set~$F\eqdef \bigcup_k F_k$. With no loss of generality by~\cite[Lem.~2.1.3, p.~69]{FukOshTak11}, we may and will always assume that every nest is \emph{increasing}, and regular w.r.t.\ a measure apparent from context.

Let~$\seq{F_k}_k$, resp.~$\tseq{F_k^\sharp}_k$, be an $E$-, resp.~$E^\sharp$-, nest, additionally so that~$j\colon F\rar F^\sharp$ restricts to a homeomorphism $j\colon F_k\rar F_k^\sharp$ for every~$k$. Since~$\seq{F_k}_k$ is increasing,~$j\colon F\rar F^\sharp$ is bijective. Let~$N_1$ be $\nu$-negligible so that~$(E_\zeta^\sharp, \dom{E_\zeta^\sharp})$ is regular by Theorem~\ref{t:Stone}.
Let~$X_\partial\eqdef X\cup \set{\partial}$, where~$\partial$ is taken to be an isolated point in~$X_\partial$. Since~$j$ may be not surjective, in the following we extend~$j^{-1}$ on~$X^\sharp\setminus j(F)$ by setting~$j^{-1}(x^\sharp)=\partial$. Note that this extension is $\mcX^\sharp$-to-$\mcX$-measurable (having care to extend~$\mcX$ on~$X_\partial$ in the obvious way).
Since~$j_\pfwd \mu=\mu^\sharp$ the set $N_2\eqdef \ttset{\zeta\in Z: \mu_\zeta^\sharp j(F)<1}$ is $\nu$-measurable, since~$j$ is measurable on~$F$, and thus it is $\nu$-negligible. In particular,~$j^{-1}_\pfwd \mu^\sharp \set{\partial}=0$, and~$j^{-1}_\pfwd \mu_\zeta^\sharp \set{\partial}=0$ for every~$\zeta\in N_2^\complement$.
Set now~$N\eqdef N_1\cup N_2$. For~$\zeta\in N^\complement$ set~$\mu_\zeta\eqdef j^{-1}_\pfwd \mu_\zeta^\sharp$ and denote by~$(E_\zeta,\dom{E_\zeta})$ the image form of~$(E_\zeta^\sharp,\dom{E_\zeta^\sharp})$ via~$j^{-1}$ on~$L^2(\mu_\zeta)$, cf.~\cite[Eqn.~(3.2)]{CheMaRoe94}. For~$\zeta\in N$  let~$(E_\zeta,\dom{E_\zeta})$ be the $0$-form on $L^2(\mu)$.
For~$f^\sharp\colon X^\sharp\rar \R$ denote further by~$j^*f\eqdef f^\sharp\circ j\colon X\rar\R$ the pullback of~$f^\sharp$ via~$j$, and recall~\cite[Eqn.~(3.3)]{CheMaRoe94}:
\begin{align}\label{eq:CheMaRoe}
G_\alpha (j^*f^\sharp)=j^*G^\sharp_\alpha f^\sharp \comma\qquad f^\sharp\in L^2(\mu^\sharp) \fstop
\end{align}

\paragraph{2. Nests} For~$\zeta\in N^\complement$ let~$\tseq{F_{\zeta,k}^\sharp}_k$ be a $\mu_\zeta^\sharp$-regular $E_\zeta^\sharp$-nest witnessing the (quasi-)regularity of the form, i.e.\ verifying~\cite[Dfn.~2.8]{CheMaRoe94}. With no loss of generality, up to intersecting~$F_{\zeta, k}^\sharp$ with~$F_h^\sharp$ if necessary, we may assume that for every~$k$ there exists~$h\eqdef h_k$ so that~$F_{\zeta,k}^\sharp\subset F_h^\sharp$. In particular,~$j^{-1}\colon F_{\zeta,k}^\sharp\rar F_{\zeta,k}\eqdef j(F_{\zeta,k}^\sharp)$ is a homeomorphism onto its image. 
Let~$X_\zeta^\sharp\eqdef \supp\tquadre{\mu_\zeta^\sharp}$ and note that~$F_{\zeta,k}^\sharp\subset X_\zeta^\sharp$ since~$\tseq{F_{\zeta,k}^\sharp}_k$ is $\mu_\zeta^\sharp$-regular. Denote by~$j^{-1}_\zeta$ the restriction of~$j^{-1}$ to~$X_\zeta^\sharp$.

\smallskip

\emph{Claim: $j^{-1}_\zeta$ is a quasi-homeomorphism for $\nu$-a.e.~$\zeta\in Z$.} It suffices to show that~$\seq{F_{\zeta,k}}_k$ is an $E_\zeta$-nest for $\nu$-a.e.~$\zeta\in Z$, which holds by construction.

\smallskip

Finally, set~$j_\zeta\eqdef j\restr_{j^{-1}(X_\zeta^\sharp)}$ and note that, again by~\cite[Eqn.~(3.3)]{CheMaRoe94},
\begin{align}\label{eq:CheMaRoeZeta}
G_{\zeta,\alpha} (j_\zeta^*f^\sharp)=j_\zeta^*G^\sharp_{\zeta,\alpha} f^\sharp \comma\qquad f^\sharp\in L^2(\mu_\zeta^\sharp) \comma \qquad \zeta\in N^\complement\fstop
\end{align}

\paragraph{3. Direct integral representation}
By~\eqref{eq:p:DirInt:0} for the resolvent applied to~$(E^\sharp,\dom{E^\sharp})$,
\begin{align}\label{eq:CheMaRoe:0}
G_\alpha^\sharp=\dint{Z} G^\sharp_{\zeta,\alpha} \diff\nu(\zeta) \fstop
\end{align}

By \emph{Step 1} and~\cite[Lem.~3.3(ii)]{CheMaRoe94},~$j_\zeta^*\colon L^2(\mu^\sharp)\rar L^2(\mu)$ is an isomorphism for $\nu$-a.e.~$\zeta\in Z$, with inverse~$(j_\zeta^{-1})^*$, and
\begin{align}\label{eq:CheMaRoe:1}
j^*=\dint{Z} j_\zeta^*\diff\nu(\zeta) \fstop
\end{align}

Now, by a subsequent application of~\eqref{eq:CheMaRoe},~\eqref{eq:CheMaRoe:0}, \eqref{eq:CheMaRoe:1} and~\cite[\S{II.2.3}, Prop.~3, p.~182]{Dix81}, and \eqref{eq:CheMaRoeZeta},
\begin{align*}
G_\alpha \circ j^*=&\ j^*\circ G_\alpha^\sharp=j^*\circ\, \dint{Z} G^\sharp_{\zeta,\alpha} \diff\nu(\zeta) = \dint{Z} j_\zeta^*\circ G^\sharp_{\zeta,\alpha} \diff\nu(\zeta)
\\
=&\ \dint{Z} G_{\zeta,\alpha}\circ j_\zeta^* \diff\nu(\zeta) \fstop
\end{align*}
Thus, by a further application of~\cite[\S{II.2.3}, Prop.~3, p.~182]{Dix81},
\begin{align*}
G_\alpha \circ j^*=\dint{Z} G_{\zeta,\alpha} \diff\nu(\zeta)\circ j^*\fstop
\end{align*}
Cancelling~$j^*$ by its inverse~$(j^{-1})^*$, this yields the direct-integral representation of~$G_\alpha$ via~$\zeta\mapsto G_{\alpha,\zeta}$. By~\eqref{eq:p:DirInt:0} for the resolvent, this shows
\begin{align*}
E=\dint{Z} E_\zeta \diff\nu(\zeta) \fstop
\end{align*}
 
\paragraph{4. Quasi-regularity and irreducibility} By \emph{Step 2} and~\cite[Thm.~3.7]{CheMaRoe94}, the form $(E_\zeta,\dom{E_\zeta})$ is quasi-regular for $\nu$-a.e.~$\zeta\in Z$. Again by \emph{Step 2}, it is also irreducible, since it is isomorphic to the irreducible form~$(E^\sharp_\zeta,\dom{E^\sharp_\zeta})$.
\end{proof}

\subsection{Ergodic decomposition of forms: \texorpdfstring{$\sigma$}{sigma}-finite measure case}\label{ss:ErgSigmaFinite}
Under some additional assumptions, we may now extend the results in Theorem~\ref{t:Stone} to the case when~$\mu$ is only $\sigma$-finite. The main idea ---~borrowed from~\cite{Buf14}~--- is to reduce the $\sigma$-finite case to the probability case.

\paragraph{General strategy} Let~$(E,\dom{E})$ be a strongly local regular Dirichlet form on a locally compact Polish, Radon measure space~$(X,\T,\mcX^\mu,\hat\mu)$ with full support. In particular,~$(X,\mcX,\mu)$ is $\sigma$-finite,~\cite[415D(iii)]{Fre00}. Assume further that~$(E,\dom{E})$ has carr\'e du champ~$(\cdc,\dom{E})$.
Let~$\phi\in \dom{E}$, with~$\phi>0$ $\mu$-a.e. and $\norm{\phi}_{L^2(\mu)}=1$, and set~$\mu^\phi\eqdef \phi^2\cdot\mu$. Here, we understand~$\phi$ as a fixed $E$-quasi-continuous representative of its $\mu$-class. Note that~$(X,\mcX,\mu^\phi)$ is a probability space and that~$\mu^\phi$ is equivalent to~$\mu$. Therefore,~$\mu$-classes and~$\mu^\phi$-classes coincide. On~$L^2(\mu^\phi)$ we define a bilinear form
\begin{align}\label{eq:Girsanov}
\dom{E^\phi}\eqdef \set{\int (\cdc(f)+f^2)\diff\mu^\phi<\infty}\comma \qquad
E^\phi(f,g)\eqdef \int \cdc(f,g) \diff\mu^\phi\fstop
\end{align}
Suppose now that
\begin{enumerate}[$(a)$]
\item\label{i:StrategyA} the form~$(E^\phi,\dom{E^\phi})$ is a (closed) regular Dirichlet form on~$L^2(\mu^\phi)$.
\end{enumerate}
Then, we may apply Theorem~\ref{t:Stone} to obtain 
\begin{itemize}
\item a probability space~$(Z_\phi,\mcZ_\phi,\nu_\phi)$ and a measurable map~$s_\phi\colon X\rar Z$;
\item a $\nu_\phi$-essentially unique disintegration~$\tseq{\mu^{(\phi)}_\zeta}_{\zeta\in Z}$ of~$\mu^\phi$ over~$\nu_\phi$, strongly consistent with~$s_\phi$;
\item a family of regular strongly local Dirichlet forms $\ttonde{E^{(\phi)}_\zeta,\dom{E^{(\phi)}_\zeta}}$ on~$L^2\ttonde{\mu^{(\phi)}_\zeta}$;
\end{itemize}
satisfying the direct-integral representation
\begin{align}\label{eq:GirsanovDirInt}
E^\phi=\dint{Z_\phi} E^{(\phi)}_\zeta \diff\nu_\phi(\zeta)\fstop
\end{align}

For~$\zeta\in Z$ let now~$\mu_\zeta^{[\phi]}\eqdef \phi^{-2}\cdot \mu^{(\phi)}_\zeta$ be a measure on~$(X,\mcX)$ and suppose further that
\begin{enumerate}[$(a)$]\setcounter{enumi}{1}
\item\label{i:StrategyB} the form $\ttonde{E^{(\phi)}_\zeta,\dom{E^{(\phi)}_\zeta}}$ has carr\'e du champ operator~$\ttonde{\cdc^{(\phi)}_\zeta, \dom{\cdc^{(\phi)}_\zeta}}$ for $\nu_\phi$-a.e.~$\zeta\in Z_\phi$;

\item\label{i:StrategyC} the form
\begin{align}\label{eq:Girsanov2}
\begin{aligned}
\dom{E^{[\phi]}_\zeta}\eqdef&\ \set{\int \ttonde{\cdc^{(\phi)}_\zeta(f)+f^2}\diff\mu^{[\phi]}_\zeta<\infty}\comma 
\\
E^{[\phi]}_\zeta(f,g)\eqdef& \int \cdc^{(\phi)}_\zeta(f,g) \diff\mu^{[\phi]}_\zeta\comma\qquad f,g\in \dom{E^{[\phi]}_\zeta}\comma
\end{aligned}
\end{align}
is a (closed) regular Dirichlet form on~$L^2(\mu_\zeta^{[\phi]})$ for $\nu_\phi$-a.e.~$\zeta\in Z$.
\end{enumerate}

Then, finally, we may expect to have a direct-integral representation
\begin{align}\label{eq:DirIntPhi}
E=\dint{Z_\phi} E^{[\phi]}_\zeta \diff\nu_\phi(\zeta) \fstop
\end{align}

As it turns out, the properties of the Girsanov-type transformation~\eqref{eq:Girsanov} are quite delicate.
Before discussing the technical details, let us note here that, provided we have shown the direct-integral representation in~\eqref{eq:DirIntPhi}, it should not be expected that the latter is (essentially) unique, but rather merely essentially \emph{projectively} unique ---~as it is the case for other ergodic theorems, e.g.~\cite[Thm.~2]{Buf14}. In the present setting, projective uniqueness is understood in the following sense.

\begin{definition}\label{d:EPU} We say that the direct integral representation~\eqref{eq:DirIntPhi} is \emph{essentially projectively unique} if, for every~$\phi$,~$\psi$ as above:
\begin{enumerate}[$(a)$]
\item\label{i:d:EPU:1} the measurable space~$(Z,\mcZ)\eqdef (Z_\phi,\mcZ_\phi)=(Z_\psi,\mcZ_\psi)$ is uniquely determined;
\item\label{i:d:EPU:2} the measures~$\nu_\phi$,~$\nu_\psi$ are equivalent (i.e., mutually absolutely continuous);
\item\label{i:d:EPU:3} the forms~$E^{[\phi]}_\zeta$, $E^{[\psi]}_\zeta$ are multiple of each other for $\nu_\phi$- \mbox{(hence~$\nu_\psi$-)}a.e.~$\zeta\in Z$.
\end{enumerate}
\end{definition}

As it is clear, the definition only depends on the $\sigma$-ideal of $\nu_\phi$-negligible sets in~$\mcZ$. By condition~\ref{i:d:EPU:2}, this ideal does not, in fact, depend on~$\nu_\phi$, hence the omission of the measure in the designation. The lack of uniqueness is shown as follows. Since~$\mu_\zeta^{[\phi]}$ is merely a $\sigma$-finite (as opposed to: probability) measure, the family~$\tseq{\mu_\zeta^{[\phi]}}_{\zeta\in Z}$ is merely a pseudo-disintegration (as opposed to: a disintegration). Thus, for every measurable $g\colon Z\rar (0,\infty)$,
\begin{align*}
\mu A=\int_Z \mu_\zeta^{[\phi]}A \; \diff\nu_\phi(\zeta)=\int_Z g(\zeta) \cdot \mu_\zeta^{[\phi]}A \; \diff \ttonde{g^{-1} \cdot \nu_\phi}(\zeta)\comma \qquad A\in\mcX \fstop
\end{align*}
Since~$g$ is defined on~$Z$, the pullback function~$f\eqdef (s_\phi)^*g$ is $\mcX_0$-measurable, i.e.\ all its level sets are $E$-invariant; by strong locality of~$(E,\dom{E})$, $f$ is $E$-quasi-continuous, and therefore an element of the extended domain~$F_e$ of~$(E,\dom{E})$.
As soon as~$f\in L^2(\mu)$, then we have the direct-integral representation
\begin{align*}
E=\dint{Z_\phi} g(\zeta) E^{[\phi]}_\zeta \diff\ttonde{g^{-1}\cdot \nu_\phi}(\zeta)\comma\qquad g(\zeta)\, E^{[\phi]}_\zeta(u)=\int_X \cdc_\zeta^{(\phi)}(u) \; \diff\ttonde{f\cdot\mu^{[\phi]}_\zeta} \fstop
\end{align*}

\paragraph{Proofs' summary} The Girsanov-type transformations~\eqref{eq:Girsanov} are thoroughly studied by A.~Eberle in~\cite{Ebe96}, where~\iref{i:StrategyA} is proved. We shall therefore start by showing~\iref{i:StrategyB} above, Lemma~\ref{l:CdC} below.
Informally, in the setting of Theorem~\ref{t:Stone}, if~$(E,\dom{E})$ has carr\'e du champ~$\cdc$, then
\begin{align}\label{eq:CdCdirInt}
\cdc=\dint{Z} \cdc_\zeta \diff\nu(\zeta) \comma
\end{align}
where~$\cdc_\zeta$ is the carr\'e du champ of~$(E_\zeta,\dom{E_\zeta})$.
Since the range of~$\cdc$ is a Banach (\emph{not} Hilbert) space, we shall need the concept of direct integrals of Banach spaces. In particular, we shall need an analogue of Proposition~\ref{p:DIntL^2} for $L^1$-spaces, an account of which is given in the Appendix, together with a proof of the next lemma.

\begin{lemma}\label{l:CdC} Under the assumptions of Theorem~\ref{t:Stone} suppose further that~$(E,\dom{E})$ admits carr\'e du champ operator~$(\cdc,\dom{E})$. Then,~$(E_\zeta,\dom{E_\zeta})$ admits carr\'e du champ operator~$(\cdc_\zeta,\dom{E_\zeta})$ for $\nu$-a.e.~$\zeta\in Z$.
\end{lemma}

\begin{lemma}\label{l:Local} Under the assumptions of Theorem~\ref{t:Stone} suppose further that~$(E,\dom{E})$ is strongly local. Then,~$(E_\zeta,\dom{E_\zeta})$ is strongly local for $\nu$-a.e.~$\zeta\in Z$. 
\end{lemma}
\begin{proof}
Note: In this proof we shall make use of results in~\cite{BouHir91}. We recall that a \emph{regular} form is `strongly local' in the sense of~\cite[p.~6]{FukOshTak11} if and only if it is `local' in the sense of~\cite[Dfn.~I.V.1.2, p.~28]{BouHir91}. This is noted e.g.\ in~\cite[\S2, p.~78]{tElRobSikZhu06}, after~\cite[Prop.~1.4]{Sch95}. In this respect, we always adhere to the terminology of~\cite{MaRoe92, FukOshTak11}.
Let~$(\mcE,\dom{\mcE})$ be a regular Dirichlet form.
By~\cite[Cor.~I.5.1.4 and Rmk.~I.5.1.5, p.~31]{BouHir91},~$(\mcE,\dom{\mcE})$ is strogly local if and only if~$\abs{u}\in\dom{\mcE}$ and~$\mcE(\abs{u})=\mcE(u)$ for every~$u\in\dom{\mcE}$.
Further note that~$\abs{u}$ is a normal contraction of~$u\in\dom{\mcE}$ for every~$u\in\dom{\mcE}$.
As a consequence,~$\abs{u}\in\dom{\mcE}$ and~$\mcE(\abs{u})\leq \mcE(u)$ for every~$u\in\dom{\mcE}$, see e.g.~\cite[Thm.~1.4.1]{FukOshTak11}.
In particular, a Dirichlet form~$(\mcE,\dom{\mcE})$ is \emph{not} strongly local if and only if there exists~$u\in\dom{\mcE}$ with~$\mcE(u)>\mcE(\abs{u})$.

Since~$\mu_\zeta X\leq 1$ for every~$\zeta\in Z$, it is not difficult to show, arguing by contradiction, that
\begin{align}
\car_X\in\dom{E_\zeta}\comma \qquad E_\zeta(\car_X)=0 \quad \forallae{\nu} \zeta\in Z \fstop
\end{align}

Let~$\seq{u_n}_n\subset \mcC$ be the fundamental sequence constructed in~\emph{Step~1} in the proof of Theorem~\ref{t:Stone}\iref{i:t:Stone:3}.
The Dirichlet form~$(E_\zeta,\dom{E_\zeta})$ on~$L^2(\mu_\zeta)$ is strongly local if and only if~$E_\zeta(u_n)=E_\zeta(\abs{u_n})$ for every~$n\geq 1$. The same holds for~$(E,\dom{E})$.

Now, argue by contradiction that there exists a $\nu$-measurable non-negligible set~$B\subset Z$ so that the form~$(E_\zeta,\dom{E_\zeta})$ is \emph{not} strongly local for each~$\zeta\in B$. Let~$B_n\eqdef \set{\zeta\in Z: E_\zeta(u_n)> E_\zeta(\abs{u_n})}$ and note that~$B_n\subset B$ is $\nu$-measurable for every~$n\geq 1$ since~$\seq{u_n}_n\subset\mcC$.
Since~$B=\cup_n B_n$ and~$\nu B>0$, there exists some fixed~$n_*$ so that~$\nu B_{n_*}>0$. Without loss of generality, up to relabeling, we may choose~$n_*=1$.
Analogously to the proof of the \emph{Claim} in \emph{Step~3} of Theorem~\ref{t:Stone}\iref{i:t:Stone:3}, set~$A\eqdef p^{-1}(B_1)$ and note that it is $E$-invariant.
Thus, finally,~$\car_A u_1\in\dom{E}$ and
\begin{align*}
E(u_1\car_A)=\int_{B_1} E_\zeta(u_1)\diff\nu(\zeta)>\int_{B_1} E_\zeta(\abs{u_1}) \diff\nu(\zeta)=E(\abs{u_1}\car_A)=E(\abs{u_1 \car_A})\comma
\end{align*}
which contradicts the strong locality of~$(E,\dom{E})$.
\end{proof}

\begin{remark}\label{r:Converse}
The converse implications to Lemmas~\ref{l:CdC},~\ref{l:Local} are true in a more general setting, viz.
\begin{enumerate}[$(a)$]
\item if~$(E,\dom{E})$ is a direct integral of Dirichlet forms~$\zeta\mapsto (E_\zeta,\dom{E_\zeta})$ each with carr\'e du champ operator~$(\cdc_\zeta,\dom{\cdc_\zeta})$, then~$(E,\dom{E})$ has carr\'e du champ given by~\eqref{eq:CdCdirInt}, see~\cite[Ex.~V.3.2, p.~216]{BouHir91};
\item if~$(E,\dom{E})$ is a direct integral of strongly local Dirichlet forms~$\zeta\mapsto (E_\zeta,\dom{E_\zeta})$, then~$(E,\dom{E})$ is strongly local, see~\cite[Ex.~V.3.1, p.~216]{BouHir91}.
\end{enumerate}
\end{remark}

We are now ready to prove the main result of this section.

\begin{theorem}[Ergodic decomposition: $\sigma$-finite case]\label{t:SFinite}
Let~$(X,\T,\mcX,\mu)$ be a locally compact Polish Radon measure space, and~$(E,\dom{E})$ be a regular strongly local Dirichlet form on~$L^2(\mu)$ with carr\'e du champ operator~$(\cdc,\dom{\cdc})$. Then, there exist
\begin{enumerate}[$(i)$]
\item a probability space~$(Z,\mcZ,\nu)$ and a measurable map~$s\colon (X,\mcX)\rar (Z,\mcZ)$;
\item an essentially projectively unique family of measures~$\seq{\mu_\zeta}_{\zeta\in Z}$ so that, when $s^{-1}(\zeta)$ is endowed with the subspace topology and the trace $\sigma$-algebra inherited by $(X,\T,\mcX^\mu)$, then~$(s^{-1}(\zeta),\hat\mu_\zeta)$ is a Radon measure space for~$\nu$-a.e.~$\zeta\in Z$;
\item a $\nu$-measurable field~$\zeta\mapsto (E_\zeta,\dom{E_\zeta})$ of regular irreducible Dirichlet forms $(E_\zeta,\dom{E_\zeta})$ on~$L^2(\mu_\zeta)$;
\end{enumerate}
so that
\begin{align*}
L^2(\mu)=\dint[]{Z} L^2(\mu_\zeta)\diff\nu(\zeta) \qquad \text{and} \qquad E=\dint[]{Z} E_\zeta \diff\nu(\zeta) \fstop
\end{align*}
\end{theorem}
\begin{proof}
Let~$\phi\in \dom{E}$ with~$0<\phi<1$ $\mu$-a.e. and~$\norm{\phi}_{L^2}=1$. Since~$(E,\dom{E})$ is regular strongly local on~$L^2(\mu)$ and admits carr\'e du champ~$(\cdc,\dom{E})$, then the Girsanov-type transform~$(E^\phi,\dom{E^\phi})$ defined in~\eqref{eq:Girsanov} is a quasi-regular strongly local Dirichlet form on~$L^2(\mu^\phi)$ by~\cite[Thm.s~1.1 and~1.4(ii)]{Ebe96}, and admits carr\'e du champ~$(\cdc, \dom{E^\phi})$ by construction.
We note that we are applying the results in~\cite{Ebe96} in the context of~\cite[Example 1), p.~501]{Ebe96}. In particular, Assumption~(D3) in~\cite[p.~501]{Ebe96} holds by definition.

\paragraph{1.~Constructions}
Let now~$\mcC$ be a core for~$(E,\dom{E})$. Since~$\phi\leq 1$ $\mu$-a.e., then~$E^\phi_1\leq E_1$, and the form~$(E^\phi,\dom{E^\phi})$ is in fact regular, with same core~$\mcC$.
Since~$\mu^\phi$ is a fully supported probability measure by construction, we may apply Theorem~\ref{t:Stone} to obtain the direct integral representation~\eqref{eq:GirsanovDirInt}. For $\nu_\phi$-a.e.~$\zeta\in Z_\phi$, the form~$\ttonde{E_\zeta^{(\phi)},\dom{E_\zeta^{(\phi)}}}$ is regular with core~$\mcC$ and irreducible by Theorem~\ref{t:Stone}, strongly local by Lemma~\ref{l:Local}, and admitting carr\'e du champ~$\ttonde{\cdc_\zeta^{(\phi)},\dom{E^{(\phi)}_\zeta}}$ by Lemma~\ref{l:CdC}.

\smallskip

\emph{Claim: $\phi^{-2}\vee n \in\dom{E^{(\phi)}_\zeta}$ for $\nu_\phi$-a.e.~$\zeta\in Z_\phi$, for every~$n\geq 1$.} Since~$E^\phi_1\leq E_1$, one has~$\phi\in \dom{E^\phi}$, hence~$\phi\in \dom{E_\zeta^{(\phi)}}$ for $\nu_\phi$-a.e.~$\zeta\in Z_\phi$ as a consequence of~\eqref{eq:DIntQF}. The claim then follows by strong locality of~$(E_\zeta,\dom{E_\zeta})$ for~$\nu_\phi$-a.e.~$\zeta\in Z_\phi$.

\smallskip

Let now~$\ttonde{E^{[\phi],n}_\zeta,\dom{E^{[\phi],n}_\zeta}}$ be defined analogously to~\eqref{eq:Girsanov2} with~$\mu_\zeta^{[\phi],n}\eqdef (\phi^{-2}\vee n) \cdot \mu_\zeta^{(\phi)}$ in place of~$\mu_\zeta^{[\phi]}$.
By applying once more~\cite[Thm.s~1.1 and~1.4(ii)]{Ebe96}, the Girsanov-type transform~$\ttonde{E^{[\phi],n}_\zeta,\dom{E^{[\phi],n}_\zeta}}$ defined in~\eqref{eq:Girsanov2} is a regular Dirichlet form on~$L^2(\mu^{[\phi],n}_\zeta)$ with core~$\mcC$, strongly local, irreducible, and admitting carr\'e du champ operator~$\ttonde{\cdc_\zeta^{(\phi)},\dom{E_\zeta^{[\phi],n}}}$ for~$\nu_\phi$-a.e.~$\zeta\in Z_\phi$, for every~$n\geq 1$.

\smallskip

\emph{Claim: the quadratic form~$\ttonde{E^{[\phi]}_\zeta,\dom{E^{[\phi]}_\zeta}}$ defined in~\eqref{eq:Girsanov2} is a regular Dirichlet form on~$L^2(\mu^{[\phi]}_\zeta)$, strongly local, irreducible, and admitting carr\'e du champ operator~$\ttonde{\cdc_\zeta^{(\phi)},\dom{E_\zeta^{[\phi]}}}$ for~$\nu_\phi$-a.e.~$\zeta\in Z_\phi$.}
Firstly, note that~$E^{[\phi]}_\zeta=\sup_n E^{[\phi],n}_\zeta$ is well-defined on~$\dom{E^{[\phi]}_\zeta}=\bigcap_{n\geq 1} \dom{E^{[\phi],n}_\zeta}$, thus~$\ttonde{E^{[\phi]}_\zeta,\dom{E^{[\phi]}_\zeta}}$ is a closable quadratic form by~\cite[Prop.~I.3.7(ii)]{MaRoe92}. The Markov property, the strong locality and the existence and computation of the carr\'e du champ operator are straightforward. Note that~$\mcC\subset \dom{E^{[\phi]}_\zeta}$, so that the latter is dense in~$L^2(\mu_\zeta^{[\phi]})$ for $\nu_\phi$-a.e.~$\zeta\in Z_\phi$. By Dominated Convergence and~\eqref{eq:DisintF}
\begin{align*}
E(u)=&\int\cdc(u)\diff\mu=\nlim \int\cdc(u)\cdot (\phi^{-2}\vee n) \diff\mu^\phi
\\
=&\nlim\int_{Z_\phi}\int_X\cdc_\zeta^{(\phi)}(u) \diff\mu_\zeta^{[\phi],n}\diff\nu_\phi(\zeta)
\\
=&\int_{Z_\phi}\int_X \cdc_\zeta^{(\phi)}(u) \diff\mu^{[\phi]}_\zeta\diff\nu_\phi(\zeta) \comma
\end{align*}
which establishes the direct integral representation~\eqref{eq:DirIntPhi}, with underlying space~$S_\mcC$. The regularity of the forms~$\ttonde{E^{[\phi]}_\zeta,\dom{E^{[\phi]}_\zeta}}$, all with core~$\mcC$, follows from the regularity of~$(E,\dom{E})$, exactly as in the \emph{Claim} in \emph{Step~4} in the proof of Theorem~\ref{t:Stone}.

\paragraph{2.~Projective uniqueness} By~\eqref{eq:GirsanovDirInt}, $E$-invariant sets are also~$E^{\phi}$-invariant. The reverse implication follows since~$\mu$ and~$\mu^\phi$ are equivalent. As a consequence, the $\sigma$-algebras~$\mcX_0$ and~$\mcX^*$ defined w.r.t.~$\mu^\phi$ as in~\S\ref{ss:AlgInvSets} are independent of~$\phi$ and thus so are the space~$(Z_\phi,\mcZ_\phi)$, henceforth denoted simply by~$(Z,\mcZ)$, and the map~$s_\phi$, henceforth denoted simply by~$s$.
Let now~$\phi,\psi\in\dom{E}$, $\phi,\psi>0$ $\mu$-a.e., and replicate every construction for~$\psi$ as well. Note that $\mu^\phi\eqdef \phi^2\cdot \mu$ and~$\mu^\psi\eqdef \psi^2\cdot \mu$ are equivalent.

\smallskip

\emph{Claim 1:~$\nu_\phi\sim \nu_\psi$.} It suffices to recall that~$\nu_\phi=s_\pfwd \mu^\phi$ and analogously for~$\psi$ w.r.t.\ the same map~$s$, hence the conclusion, since~$\mu^\phi\sim \mu^\psi$. 

\smallskip

It follows that the $\sigma$-ideal~$\mcN\eqdef \mcN_\phi$ of $\nu_\phi$-negligible sets in~$\mcZ$ does not in fact depend on~$\phi$. In the following, we write therefore ``$\mcN$-negligible'' in place of~``$\nu_\phi$-negligible'' and~``$\mcN$-a.e.''\ in place of~``$\nu_\phi$-a.e.''.

\smallskip

\emph{Claim 2:~$\mu_\zeta^{(\phi)}\sim \mu_\zeta^{(\psi)}$ for $\mcN$-a.e.~$\zeta\in Z$.} Argue by contradiction that there exist~$B\in\mcZ\setminus\mcN$ and a family~$\seq{A_\zeta}_{\zeta\in B}$ with~$A_\zeta\in\mcX$ and, without loss of generality,~$\mu_\zeta^{(\phi)}A_\zeta>\mu_\zeta^{(\psi)}A_\zeta=0$ for all~$\zeta\in B$. Set further~$\tilde A\eqdef \cup_{\zeta\in B} A_\zeta$ and let~$A\in\mcX$ be its measurable envelope~\cite[132D]{Fre00}. Then, by~\eqref{eq:Disint:0} and strong consistency of~$\ttseq{\mu_\zeta^{(\phi)}}_{\zeta\in Z}$ with~$s$,
\begin{align}\label{eq:t:SFinite:1}
\begin{aligned}
\mu^\phi A=&\ \mu^\phi\ttonde{A\cap s^{-1}(B)}=\int_B \mu_\zeta^{(\phi)} A \diff\nu_\phi(\zeta)=\int_B \mu_\zeta^{(\phi)} \ttonde{A\cap s^{-1}(\zeta)} \diff\nu_\phi(\zeta)
\\
\geq& \int_B \mu_\zeta^{(\phi)} A_\zeta \diff\nu_\phi(\zeta)>0 \fstop
\end{aligned}
\end{align}
Let now~$\mu^{\psi*}$ be the outer measure of~$\mu^\psi$, and analogously for~$\mu_\zeta^{(\psi)*}$. Note that, by Assumption~\ref{ass:Main}, the Carath\'eodory measure induced by~$\mu^{\psi *}$ coincides with the completion measure~$\widehat{\mu^\psi}=\hat\mu^\psi$. Furthermore, by strong consistency of~$\ttseq{\mu_\zeta^{(\psi)}}_{\zeta\in Z}$ with~$s$, one has~$\mu_\zeta^{(\psi)*}\tilde A=\mu_\zeta^{(\psi)*}\ttonde{\tilde A\cap s^{-1}(\zeta)}=\mu_\zeta^{(\psi)*}A_\zeta=0$ by assumption. In particular, the function~$\zeta\mapsto \mu_\zeta^{(\psi)*}\tilde A\equiv 0$ is measurable. Thus, by~\cite[452X(i)]{Fre00},
\begin{align*}
0=\int_B \mu_\zeta^{(\psi)*}\tilde A = \mu^{\psi*} \tilde A=\mu^\psi A\comma
\end{align*}
which contradicts~\eqref{eq:t:SFinite:1} since~$\mu^\phi\sim\mu^\psi$ by the previous claim.

\smallskip

By the last claim,~$\mu_\zeta^{[\phi]}$-classes and $\mu_\zeta^{[\psi]}$-classes coincide. Therefore, the carr\'e du champ operator~$\cdc_\zeta^{[\phi]}=\cdc_\zeta^{(\phi)}$ is independent of~$\phi$, and henceforth denoted by~$\cdc_\zeta$. Thus we have
\begin{align*}
E(u)=\int_Z \int_X \cdc_\zeta(u) \diff\mu_\zeta^{[\phi]} \diff\nu_\phi(\zeta) = \int_Z \int_X \cdc_\zeta(u) \diff\mu_\zeta^{[\psi]} \diff\nu_\psi(\zeta) \comma \qquad u\in\mcC\comma
\end{align*}
and, finally, it suffices to show the following. 

\smallskip

\emph{Claim 3:~$\mu_\zeta^{[\phi]}=\frac{\diff \nu_\psi}{\diff\nu_\phi}\cdot\mu_\zeta^{[\psi]}$.} By construction,~$\tseq{\mu_\zeta^{[\phi]}}_{\zeta\in Z}$, resp.~$\tseq{\mu_\zeta^{[\psi]}}_{\zeta\in Z}$, is a pseudo-disintegration of~$\mu$ over~$\nu_\phi$, resp.~$\nu_\psi$. For fixed~$f\in L^1(\mu)^+$ and every~$t>0$ set~$A_t\eqdef \ttset{f/\phi^2=t}$. By consistency of the disentegration of~$\mu^\phi$ over~$\nu_\phi$,
\begin{align*}
\int_0^\infty \mu^\phi\ttonde{A_t\cap s^{-1}(B)} \diff t=\int_0^\infty \int_B \mu_\zeta^{(\phi)} A_t \diff \nu_\phi(\zeta) \diff t\comma \qquad B\in\mcZ\fstop
\end{align*}
whence, by the level-set representation of the Lebesgue integral and Tonelli's Theorem
\begin{align}\label{eq:Uniqueness:1}
\int_{s^{-1}(B)} f \diff\mu = \int_B \int_X f \diff\mu_\zeta^{[\phi]} \diff\nu_\phi(\zeta) \fstop
\end{align}
Note that the left-hand side does not depend on~$\phi$. Therefore, equating~\eqref{eq:Uniqueness:1} with the same representation for~$\psi$ and using~\emph{Claim 1} yields
\begin{align*}
\int_B \int_X f \diff\mu_\zeta^{[\phi]} \diff\nu_\phi(\zeta) =\int_B \int_X f \diff\mu_\zeta^{[\psi]} \diff\nu_\psi(\zeta) = \int_B \frac{\diff\nu_\psi}{\diff\nu_\phi}(\zeta) \int_X f \diff\mu_\zeta^{[\psi]} \diff\nu_\phi(\zeta)\comma
\end{align*}
and the conclusion follows, since~$f$ and~$B$ were arbitrary.
\end{proof}

\subsection{Some examples}\label{ss:Examples}
Here, we specialize the results in the previous sections to some particular cases.

In order to discuss the next example, we shall need the definition of $1$-\emph{capacity} $\Cap_E$ of a regular Dirichlet form~$(E,\dom{E})$ on a locally compact Polish space~$(X,\T)$, for which we refer the reader to~\cite[\S{2.1}]{FukOshTak11}.
We recall that~$\Cap_E$ is a \emph{Choquet capacity} on~$X$ in the sense of e.g.~\cite[\S{A.1}]{FukOshTak11}, see e.g.~\cite[Thm.~2.1.1]{FukOshTak11}.
Finally, we say that~$A\subset X$ is $E$-capacitable if
\begin{align*}
\Cap_E(A)=\sup_{K\in\msK_\T : K\subset A} \Cap_E(K)\comma
\end{align*}
where~$\msK_\T$ denotes the family of all $\T$-compact subsets of~$X$.

\begin{example}[Ergodic decomposition of forms on product spaces]
Let~$X=Y\times Z$ be a product of locally compact Polish spaces endowed with a probability (hence Radon) measure~$\mu$, and~$\seq{\mu_\zeta}_{\zeta\in Z}$ be a disintegration of~$\mu$ over~$\nu\eqdef \pr^Z_\pfwd \mu$ strongly consistent with the standard projection~$\pr^Z\colon X\rar Z$. This includes the setting of Example~\ref{ese:R2}.
Indeed, let~$(E_\zeta,\dom{E_\zeta})$ be regular irreducible Dirichlet forms on~$L^2(\mu_\zeta)$, all with common core~$\mcC\subset \Cz(Y)$, and assume that~$\zeta\mapsto (E_\zeta,\dom{E_\zeta})$ is a $\nu$-measurable field of quadratic forms in the sense of Definition~\ref{d:assQ} with underlying $\nu$-measurable field~$S=S_\mcC$.
Then, it is readily verified that
\begin{enumerate}[$(i)$]
\item the direct integral~$(E,\dom{E})$ of quadratic forms~$\zeta\mapsto (E_\zeta,\dom{E_\zeta})$ is a direct integral of Dirichlet forms;
\item $(E,\dom{E})$ is a regular Dirichlet form on~$L^2(\mu)$ with core~$\mcC\otimes \Cz(Z)$ and semigroup
\begin{align*}
(T_t u)(y,\zeta)=\ttonde{T_{\zeta, t}\otimes \id_{L^2(\nu)} u }(y,\zeta)= \ttonde{T_{\zeta, t} u(\emparg, \zeta)}(y) \semicolon
\end{align*}
\item if~$A\subset Z$ is $\nu$-measurable and~$U\subset Y$ is~$E_\zeta$-capacitable for every~$\zeta\in Z$, then~$U \times A\subset X$ satisfies
\begin{align*}
\Cap_E(U\times A)\leq \int_A \Cap_{E_\zeta}(U) \diff\nu(\zeta) \fstop
\end{align*}
\end{enumerate}
\end{example}

As a further example, we state here the ergodic decomposition theorem for mixed Poisson measures on the configuration space over a \emph{connected} Riemannian manifold. We refer the reader to~\cite{AlbKonRoe98} for the main definitions.

\begin{example}[Mixed Poisson measures,~{\cite{AlbKonRoe98}}]
Let~$(M,g)$ be a Riemannian manifold with infinite volume, and~$\sigma=\rho\cdot \vol_g$ be a non-negative Borel measure on~$M$ with density~$\rho>0$ $\vol_g$-a.e., and satisfying~$\rho^{1/2}\in W^{1,2}_\loc(M)$. Let further~$\Gamma_M$ be the configuration space over~$M$, endowed with the vague topology and the induced Borel $\sigma$-algebra, and denote by~$\PP_\sigma$ the Poisson measure on~$\Gamma_M$ with intensity measure~$\sigma$. Let now~$\lambda$ be a Borel probability measure on~$\R_+\eqdef (0,\infty)$ with finite second moment. The \emph{mixed Poisson measure with intensity measure~$\sigma$ and L\'evy measure~$\lambda$} is the measure
\begin{align*}
\mu_{\lambda,\sigma}\eqdef \int_{\R_+} \PP_{s \sigma} \diff\lambda(s)\fstop
\end{align*}

In~\cite{AlbKonRoe98}, Albeverio, Kondratiev and~R\"ockner construct a canonical Dirichlet form~$\ttonde{\mcE_{\mu_{\lambda,\sigma}},\dom{\mcE_{\mu_{\lambda,\sigma}}}}$ on~$L^2(\mu_{\lambda,\sigma})$ and show that
\begin{itemize}
\item $\ttonde{\mcE_{\mu_{\lambda,\sigma}},\dom{\mcE_{\mu_{\lambda,\sigma}}}}$ is quasi-regular strongly local,~\cite[Thm.~6.1]{AlbKonRoe98};
\item $\ttonde{\mcE_{\mu_{\lambda,\sigma}},\dom{\mcE_{\mu_{\lambda,\sigma}}}}$ is irreducible if and only if $\lambda=\delta_s$, i.e.~$\mu_{\lambda,\sigma}=\PP_{s\sigma}$, for some~$s\geq 0$, \cite[Thm.~4.3]{AlbKonRoe98};
\item $\PP_{s \sigma} \perp \PP_{r\sigma}$ for all~$r,s\geq 0$, $r\neq s$.
\end{itemize}

Applying Theorem~\ref{t:QRegular} to the form~$\ttonde{\mcE_{\mu_{\lambda,\sigma}},\dom{\mcE_{\mu_{\lambda,\sigma}}}}$ yields the direct-integral representation
\begin{align*}
\mcE_{\mu_{\lambda,\sigma}}=\dint{\R_+} \mcE_{\PP_{s\sigma}} \diff\lambda(s) \comma
\end{align*}
where~$(Z,\mcZ,\nu)=(\R_+,\Bo{\R_+},\lambda)$, and the disintegration of~$\mu_{\lambda,\sigma}$ constructed in the theorem coincides with~$\seq{\PP_{s\sigma}}_{s\in \R_+}$.
\end{example}

\begin{remark}
Other examples are given by~\cite[Thm.~3.7]{AlbKonRoe97} and~\cite{AlbRoe90}, both concerned with strongly local Dirichlet forms on locally convex topological vector spaces, and by~\cite{LzDS17+}, concerned with a particular quasi-regular Dirichlet form on the space of probability measures over a closed Riemannian manifold.
\end{remark}

\subsection{Some applications}\label{ss:Appl}
We collect here some applications of the direct-integral decomposition discussed in the previous sections.

\paragraph{Transience/recurrence}
Let~$(X,\T,\mcX,\mu)$ be satisfying Assumption~\ref{ass:Main}.
For an invariant set~$A\in\mcX_0$, we denote by~$\mu_A$ the restriction of~$\mu$ to~$A$, and by~$(E_A,\dom{E_A})$ the Dirichlet form
\begin{align*}
\dom{E_A}\eqdef \set{\car_A f: f\in\dom{E}}\comma \qquad E_A(f,g)\eqdef E(\car_A f,\car_A g) \comma
\end{align*}
well-defined on~$L^2(\mu_A)$ as a consequence of Definition~\ref{d:Invariant}\iref{i:d:Invariant:4}.
The next result is standard. In the generality of Assumption~\ref{ass:Main}, a proof is readily deduced from the corresponding result for $\mu$-tight Borel right processes, shown with different techniques by K.~Kuwae in~\cite[Thm.~1.3]{Kuw11}, in the far more general setting of quasi-regular semi-Dirichlet forms.

\begin{corollary}[Ergodic decomposition: transience/recurrence]\label{c:Kuwae}
Under the assumptions of Theorem~\ref{t:QRegular}, there exist $E$-invariant subsets~$X_c$,~$X_d$, and a properly $E$-exceptional subset~$N$ of~$X$, so that
\begin{enumerate}[$(i)$]
\item $X= X_c \sqcup X_d \sqcup N$;
\item the restriction~$(E_d,\dom{E_d})$ of~$(E,\dom{E})$ to~$X_d$ is transient;
\item the restriction~$(E_c,\dom{E_c})$ of~$(E,\dom{E})$ to~$X_c$ is recurrent.
\end{enumerate}
\end{corollary}

As an application, we have the following proposition. Similarly to Remark~\ref{r:Converse}, some implications hold for superpositions of arbitrary Dirichlet forms.

\begin{proposition} Let~$(X,\T,\mcX,\mu)$ be a topological measure space satisfying Assumption~\ref{ass:Main}, and $(Z,\mcZ,\nu)$ be $\sigma$-finite countably generated. Further let~$\seq{\mu_\zeta}_{\zeta\in Z}$ be a separated pseudo-disintegration of~$\mu$ over~$\nu$, and~$(E,\dom{E})$ be a direct integral of quasi-regular Dirichlet forms~$\zeta\mapsto (E_\zeta,\dom{E_\zeta})$ on~$L^2(\mu_\zeta)$. Then,
\begin{enumerate}[$(i)$]
\item\label{i:l:RecTran:1} $(E,\dom{E})$ is conservative if and only if~$(E_\zeta,\dom{E_\zeta})$ is conservative for $\nu$-a.e.~$\zeta\in Z$;

\item\label{i:l:RecTran:2} $(E,\dom{E})$ is transient if and only if~$(E_\zeta,\dom{E_\zeta})$ is transient for $\nu$-a.e.~$\zeta\in Z$. Furthermore, one has the direct-integral representation of Hilbert spaces
\begin{align}\label{eq:l:RecTran:0}
\domext{E}=\dint{Z} \domext{E_\zeta}\diff\nu(\zeta) \semicolon
\end{align}

\item\label{i:l:RecTran:3} if $(E,\dom{E})$ is recurrent, then~$(E_\zeta,\dom{E_\zeta})$ is recurrent for $\nu$-a.e.~$\zeta\in Z$. In the situation of Theorem~\ref{t:Stone} or Theorem~\ref{t:SFinite}, the converse implication holds as well.
\end{enumerate}
\end{proposition}
\begin{proof}
Analogously to the proof of Theorem~\ref{t:QRegular} we may restrict to the regular case by the transfer method. Thus we can and will assume with no loss of generality that~$(X,\T,\mcX^\mu,\hat\mu)$ is a locally compact Polish Radon measure space with full support, and~$(E,\dom{E})$ is a direct integral of Dirichlet forms~$\zeta\mapsto (E_\zeta,\dom{E_\zeta})$ with underlying space of $\nu$-measurable vector fields~$S_\mcC$ generated by a core~$\mcC$ for~$(E,\dom{E})$.
By this assumption, the form~$(E,\dom{E})$ and all forms~$(E_\zeta,\dom{E_\zeta})$ are regular, with common core~$\mcC$. Analogously to the proof of Theorem~\ref{t:Stone}, if~$u\in\mcC$, then we may choose~$u$ as a representative for~$u_\zeta$, thus writing~$u$ in place of~$u_\zeta$ for every~$\zeta\in Z$ and every~$u\in\mcC$. Without loss of generality, possibly up to enlargement of~$\mcC$, we may assume that~$\mcC$ is special standard (e.g.,~\cite[p.~6]{FukOshTak11}). In particular,~$\mcC$ is a lattice.

\smallskip

\iref{i:l:RecTran:1} 
Let~$\seq{u_n}_n\subset \Cc(\T)$ with $0\leq u_n\leq 1$ and so that~$\nlim u_n\equiv \car$ locally uniformly on~$(X,\T)$, and note that $\nlim u_n=\car$ both $\mu$- and~$\mu_\zeta$-a.e.\ for every~$\zeta\in Z$. By the direct-integral representation~\eqref{eq:p:DirInt:0} of~$T_\bullet$ and by~\eqref{eq:DirIntL^2:0},
\begin{align*}
\int_X f \, T_t u_n \diff\mu =\int_Z \int_X f_\zeta \, T_{\zeta, t} u_n \diff\mu_\zeta \diff\nu(\zeta) \comma \qquad \begin{gathered} f\in L^1(\mu)\cap L^2(\mu)\comma\\ f=\seq{\zeta\mapsto f_\zeta}\comma\end{gathered} \qquad t> 0\fstop
\end{align*}
Letting~$n$ to infinity, it follows by several applications of the Dominated Convergence Theorem that
\begin{align}\label{eq:l:RecTran:1}
\int_X f \, T_t\! \car \diff\mu =\int_Z \int_X f_\zeta \, T_{\zeta, t}\! \car \diff\mu_\zeta \diff\nu(\zeta) \comma \qquad \begin{gathered} f\in L^1(\mu)\cap L^2(\mu)\comma\\ f=\seq{\zeta\mapsto f_\zeta}\comma\end{gathered} \qquad t> 0\fstop
\end{align}

Now, assume that~$T_\bullet$ is \emph{not} conservative and argue by contradiction that~$T_{\zeta,\bullet}$ is conservative for $\nu$-a.e.~$\zeta\in Z$. Then, choosing~$f\in L^1(\mu)\cap L^2(\mu)$, with~$f>0$ $\mu$-a.e., in~\eqref{eq:l:RecTran:1},
\begin{align*}
\int_X f\diff\mu >& \int_X f\,  T_t\! \car \diff\mu =\int_Z \int_X f_\zeta\, T_{\zeta, t}\!\car \diff\mu_\zeta\diff\nu(\zeta)=\int_Z f_\zeta \car \diff\mu_\zeta\diff\nu(\zeta)
\\
=& \int_X f\diff\mu\comma
\end{align*}
a contradiction. The reverse implication follows from~\eqref{eq:l:RecTran:1} in a similar way.

\smallskip
 
\iref{i:l:RecTran:2}
Assume~$(E_\zeta,\dom{E_\zeta})$ is transient for every~$\zeta\in N^\complement$ for some $\nu$-negligible~$N\subset Z$. 
That is,~$\domext{E_\zeta}$ is a Hilbert space with inner product~$E_\zeta$ for every~$\zeta\in N^\complement$.
By (the proof of)~\cite[Lem.~1.5.5, p.~42]{FukOshTak11}, the space~$\dom{E_\zeta}$ is $E_\zeta^{1/2}$-dense in~$\domext{E_\zeta}$ for every~$\zeta\in N^\complement$. Thus, the space of $\nu$-measurable vector fields~$S_\mcC$ is underlying to each of the direct integrals
\begin{align*}
\dom{E}_1=&\ \dint[S_\mcC]{Z} \dom{E_\zeta}_1 \diff\nu(\zeta)\comma
\\
F_e\eqdef&\ \dint[S_\mcC]{Z} \domext{E_\zeta} \diff\nu(\zeta) \comma
\\
L^2(\mu)=&\ \dint[S_\mcC]{Z} L^2(\mu_\zeta) \diff\nu(\zeta)\fstop
\end{align*}
In particular, there exists a sequence~$\seq{u_n}_n\subset \mcC$ simultaneously $\dom{E}_1$-, $F_e$- and~$L^2(\mu)$-fundamental in the sense of Definition~\ref{d:DirInt}. Denote by~$\iota_{\zeta,e}$ the identity of~$L^2(\mu_\zeta)$, regarded as the continuous embedding~$\iota_{\zeta, e}\colon\dom{E_\zeta}_1\rar \domext{E_\zeta}$ and note that 
\begin{align*}
\zeta\mapsto \scalar{\iota_{\zeta,e} u_n}{ u_m }_{\domext{E_\zeta}} = E_\zeta(u_n, u_m) = \scalar{u_n}{u_m}_{\dom{E_\zeta}_1}- \scalar{u_n}{u_m}_{L^2(\mu_\zeta)}
\end{align*}
is $\nu$-measurable for every~$n,m$. By~\cite[\S{II.1.4},~Prop.~2, p.~166]{Dix81} this implies that~$\zeta\mapsto \iota_{\zeta,e}$ is a $\nu$-measurable field of bounded operators. Writing
\begin{align*}
\iota_e\colon \dom{E}_1\rar F_e\comma u\mapsto \dint{Z} \iota_{\zeta,e} u_\zeta \diff\nu(\zeta) \comma
\end{align*}
and arguing as in the proof of Proposition~\ref{p:DirInt}, the map~$\iota_e$ is injective, and thus it is a continuous embedding of~$\dom{E}_1$ into the Hilbert space~$(F_e, E)$. As a consequence,~$K\eqdef \cl_{F_e}\ttonde{\iota_e \dom{E}_1}$ is a Hilbert space with scalar product~$E$. By definition of~$\domext{E}$, the identity map~$\iota_e$ is a continuous embedding~$(\domext{E},E)\subset (K,E)\subset (F_e,E)$. In particular,~$(\domext{E}, E)$ is a Hilbert space, and the form~$(E,\dom{E})$ is transient by~\cite[Thm.~1.6.2, p.~58]{FukOshTak11}.

\smallskip

\emph{Claim: $\domext{E}=K=F_e$ and~\eqref{eq:l:RecTran:0} holds.} By (the proof of)~\cite[Lem.~1.5.5, p.~42]{FukOshTak11}, the space~$\dom{E}$ is $E^{1/2}$-dense in~$\domext{E}$, thus the same holds for~$\mcC$. It suffices to show that~$\mcC$ is $E^{1/2}$-dense in~$F_e$ as well. We denote by~$\mcC^\perp$ the $E$-orthogonal complement of~$\mcC$ in~$F_e$, resp.\ by~$\mcC^{\perp_\zeta}$ the $E_\zeta$-orthogonal complement of~$\mcC$ in~$\dom{E_\zeta}$. By assumption,~$\mcC^{\perp_\zeta}=\set{0}$ for every~$\zeta\in N^\complement$. Finally, by the direct-integral representation of~$F_e$,
\begin{align*}
\mcC^\perp=\dint{Z} \mcC^{\perp_\zeta}\diff\nu(\zeta)=\set{0} \comma
\end{align*}
similarly to the proof of the \emph{Claim} in \emph{Step~4} of Theorem~\ref{t:Stone}.

\smallskip

We say that~$u,v\in\mcC$ are $E$-equivalent if~$E(u-v)=0$, and we write~$u\sim v$. Let the analogous definition for~$u \sim_\zeta v$ be given. By the direct-integral representation~\eqref{eq:FormDirInt} of~$(E,\dom{E})$, it is readily seen that
\begin{align}\label{eq:l:RecTran:2}
u\sim v \iff u\sim_\zeta v\, \forallae{\nu} \zeta\in Z \fstop
\end{align}

Assume now that~$(E,\dom{E})$ is transient. That is,~$\domext{E}$ is a Hilbert space with inner product~$E$. It suffices to show~\eqref{eq:l:RecTran:0}. Since~$E^{1/2}$ is non-degenerate on~$\domext{E}$, it is non-degenerate on~$\mcC$, thus~$u \sim v$ if and only if~$u=v$ $\mu$-a.e. By~\eqref{eq:l:RecTran:2}, $E_\zeta^{1/2}$ is non-degenerate on~$\mcC$ for~$\nu$-a.e.~$\zeta\in Z$, thus~$(\mcC,E_\zeta^{1/2})$ is a pre-Hilbert space for~$\nu$-a.e.~$\zeta\in Z$. For each~$\zeta\in Z$ denote by~$K_\zeta$ the abstract completion of~$(\mcC,E_\zeta^{1/2})$, endowed with the non-relabeled extension of~$E_\zeta^{1/2}$. It is a straightforward verification that there holds the direct-integral representation
\begin{align}\label{eq:l:RecTran:3}
\domext{E}=\dint[S_\mcC]{Z} K_\zeta \diff\nu(\zeta) \fstop
\end{align}
By definition of~$\domext{E_\zeta}$, the completion embedding~$\iota_\zeta\colon \mcC\rar K_\zeta$ extends to a setwise injection $\bar\iota_\zeta\colon\domext{E_\zeta}\rar K_\zeta$. 
Indeed, let~$u_\zeta\in \domext{E_\zeta}$ and~$\seq{u_n}_n\subset \mcC$ be its approximating sequence. Since~$\seq{u_n}_n$ is, by definition, $E_\zeta^{1/2}$-Cauchy, it converges to some~$h_\zeta\in K_\zeta$ by completeness of~$K_\zeta$. Set~$\bar\iota_\zeta(u_\zeta)\eqdef h_\zeta$, and note that the definition is well-posed since~$E_\zeta^{1/2}$ is a norm in~$K_\zeta$.
Thus,~$\domext{E_\zeta}$, identified with a subset of~$K_\zeta$ via~$\bar\iota_\zeta$, is a pre-Hilbert space with scalar product~$E_\zeta$, and in fact it holds that~$\domext{E_\zeta}=K_\zeta$ by $E_\zeta^{1/2}$-density of~$\mcC$ in~$K_\zeta$. We note that equality~$\domext{E_\zeta}=K_\zeta$ is not a mere isomorphism of Hilbert spaces, but rather an extension of the completion embedding~$\iota_\zeta$, thus preserving the lattice property of~$\mcC$ regarded as a subspace of both~$\domext{E_\zeta}$ and~$K_\zeta$. Together with~\eqref{eq:l:RecTran:3}, this shows~\eqref{eq:l:RecTran:0}.

\smallskip 

\iref{i:l:RecTran:3} Assume~$(E,\dom{E})$ is recurrent. By~\cite[Thm.~1.6.3, p.~58]{FukOshTak11}
there exists a sequence~$\seq{u_n}_n\subset \dom{E}$, so that~$\nlim u_n(x)=1$ for~$\mu$-a.e.~$x\in X$, and~$\nlim E(u_n)=0$. By the Markov property for~$(E,\dom{E})$ we may assume that~$u_n\in [0,1]$. By regularity of~$(E,\dom{E})$, we may assume that~$\seq{u_n}_n \subset \mcC^+\subset \Cz(\T)$.
By e.g.~\cite[Prop.~3.1(iii), p.~690]{Kuw98},~$E(u\vee v)\leq E(u)+ E(v)$, thus, up to passing to a suitable non-relabeled subsequence, we may assume that~$\seq{u_n}_n$ is monotone non-decreasing. Then,~$\nlim u_n\equiv \car$ $\T$-locally uniformly on~$\supp[\mu]=X$ by Dini's Theorem, and therefore~$\nlim u_n(x)=1$ for~$\mu_\zeta$-a.e.~$x\in X$ for every~$\zeta\in Z$.
By the direct integral representation~\eqref{eq:FormDirInt}, it is readily seen arguing by contradiction that~$\nlim E_\zeta(u_n)=0$ for $\nu$-a.e.~$\zeta\in Z$. As a consequence,~$(E_\zeta,\dom{E_\zeta})$ is recurrent for~$\nu$-a.e.~$\zeta\in Z$, again by~\cite[Thm.~1.6.3, p.~58]{FukOshTak11}.

\smallskip

Suppose now that~$(E,\dom{E})$ is given as the direct integral of Dirichlet forms in Theorem~\ref{t:Stone}, and assume that~$(E_\zeta,\dom{E_\zeta})$ is recurrent for $\nu$-a.e.~$\zeta\in Z$. We show that~$(E,\dom{E})$ is recurrent. A proof in the setting of Theorem~\ref{t:SFinite} is nearly identical, and therefore it is omitted.

Recall the notation in~\S\ref{ss:AlgInvSets} and argue by contradiction that~$(E,\dom{E})$ is not recurrent. By Corollary~\ref{c:Kuwae}, there exists an $E$-invariant subset~$X_d$, with~$\mu X_d>0$, so that~$(E_d,\dom{E_d})$ is transient. Since~$\mcX_0$ is $\mu$-essentially countably generated by~$\mcX^*$, we may and shall assume without loss of generality that~$X_d\in \mcX^*$, so that~$B\eqdef s(X_d)\in\mcZ$. Since~$\mu X_d>0$, we have~$\nu B>0$.
It is not difficult to show that the direct-integral decomposition of~$L^2(\mu)$ splits as a direct sum of Hilbert spaces
\begin{align*}
L^2(\mu)\cong L^2(\mu_{X_d})\oplus L^2(\mu_{X_d^\complement})\cong \dint[]{B}L^2(\mu_\zeta)\diff\nu(\zeta) \oplus \dint[]{B^\complement} L^2(\mu_\zeta)\diff\nu(\zeta)\fstop
\end{align*}
Since~$X_d$ is $E$-invariant, a corresponding direct-integral decomposition of~$\dom{E}$ is induced by Corollary~\ref{c:Kuwae}
\begin{align*}
\dom{E}_1\cong \dom{E_d}_1\oplus\dom{E_c}_1\cong \dint[]{B} \dom{E_\zeta}_1 \diff\nu(\zeta) \oplus \dint[]{B^\complement}\dom{E_\zeta}_1 \diff\nu(\zeta)\fstop
\end{align*}

Applying the reverse implication in~\iref{i:l:RecTran:2}, it follows from the transience of $(E_d,\dom{E_d})$ that~$\dom{E_\zeta,\dom{E_\zeta}}$ is transient for $\nu$-a.e.~$\zeta\in B$. Since~$\nu B>0$, this contradicts the assumption.
\end{proof}

\paragraph{Ergodic decomposition of measures}
Let~$(X,\T,\mcX,\mu)$ be a locally compact Polish probability space. Since~$(X,\mcX)$ is a standard Borel space, the space~$\mfM$ of all $\sigma$-finite measures on~$(X,\mcX)$ is a standard Borel space as well when endowed with the $\sigma$-algebra generated by the family of sets
\begin{align*}
\set{\eta\in\mfM: a_1<\eta A < a_2}\comma \qquad a_1,a_2\in \R_+\comma A\in\mcX \fstop
\end{align*}

Let now~$(E,\dom{E})$ be a regular Dirichlet form on~$L^2(\mu)$, and
\begin{equation*}
\mbfM\eqdef\tseq{\Omega, \mcF,\seq{M_t}_{t\geq 0}, \seq{P_x}_{x\in X_\partial},\xi}
\end{equation*}
be the properly associated right process. We set
\begin{align*}
p_t(x,A)\eqdef P_x\set{\omega\in \Omega : M_t(\omega)\in A}\comma \qquad x\in X_\partial\comma t\geq 0\comma A\in\mcX_\partial \fstop
\end{align*}
The semigroup~$T_\bullet$ of~$(E,\dom{E})$ is thus well-defined on bounded Borel measurable functions, by letting 
\begin{equation*}
\begin{aligned}
T_t\colon \mcX_b(\R)&\longrar \mcX_b(\R)
\\
f &\longmapsto \int_X f(y) \,p_t(\emparg,\diff y) 
\end{aligned}\comma \qquad t\geq 0\fstop
\end{equation*}

\begin{definition}
We say that a $\sigma$-finite measure~$\eta$ on~$(X,\mcX)$ is \emph{$T_\bullet$-invariant} if
\begin{align*}
\int_X T_t f\diff\eta=\int_X f\diff\eta \comma \qquad f\in \mcX_b(\R)\comma t\geq 0 \fstop
\end{align*}
An invariant measure~$\eta$ is~\emph{$T_\bullet$-ergodic} if every $E$-invariant set is either $\eta$-negligible or $\eta$-conegligible.
We denote by~$\mfM_\inv$, resp.~$\mfM_\erg$, the set of all $\sigma$-finite $T_\bullet$-invariant, resp.\ $T_\bullet$-ergodic, measures.
\end{definition}

The formulation of the following result is adapted from~\cite[Thm.~1]{Buf14}. In light of Corollary~\ref{c:Kuwae}, we may restrict to the case of \emph{recurrent} Dirichlet forms.

\begin{corollary}
Let~$(X,\T,\mcX,\mu)$ be a probability space satisfying Assumption~\ref{ass:Main}, and let~$(E,\dom{E})$ be a recurrent quasi-regular Dirichlet form on~$L^2(\mu)$. Then, there exists a properly $E$-coexceptional subset~$X_\inv$ of~$X$, and a surjective map~$\pi\colon X_\inv\rar \mfM_\erg$ so that
\begin{enumerate}[$(i)$]
\item\label{i:t:ErgodicM:1} for every~$\lambda\in\mfM_\erg$ the set~$\pi^{-1}(\lambda)$ is $\lambda$-conegligible;
\item\label{i:t:ErgodicM:2} for every~$\eta\in\mfM_\inv$,
\begin{align*}
\eta=\int_{\mfM_\erg} \lambda\, \diff \n\eta(\lambda) \comma \qquad \n\eta\eqdef \pi_\pfwd \eta\semicolon
\end{align*}

\item\label{i:t:ErgodicM:3} the map~$\pi_\pfwd\colon \mfM_\inv\rar \mfM(\mfM_\erg)$ is a Borel isomorphism;

\item\label{i:t:ErgodicM:4} for any~$\eta_1,\eta_2\in \mfM_\inv$ one has~$\eta_1\ll \eta_2$ if and only if~$\pi_\pfwd\eta_1 \ll \pi_\pfwd\eta_2$, and~$\eta_1\perp \eta_2$ if and only if~$\pi_\pfwd\eta_1 \perp \pi_\pfwd\eta_2$.
\end{enumerate}
\end{corollary}
\begin{proof} As a consequence of Theorem~\ref{t:QRegular}, we may restrict to the case when~$(X,\T)$ is a locally compact Polish space. This reduces measurability statements to the case of standard Borel spaces.

By Theorem~\ref{t:Stone}\iref{i:t:Stone:3}, there exists a $\nu$-negligible set~$N\in\mcZ^\nu$ so that, for every~$\zeta\in N^\complement$,
\begin{enumerate*}[$(a)$]
\item\label{i:t:proof:ErgodicM:1} $\mu_\zeta s^{-1}(\zeta)=1$, in particular,~$\mu_\zeta$ is a probability measure (as opposed to: sub-probability);
\item\label{i:t:proof:ErgodicM:2} $(E_\zeta,\dom{E_\zeta})$ is a regular irreducible recurrent Dirichlet form on~$L^2(\mu_\zeta)$ over the space~$\supp[\mu_\zeta]$.
\end{enumerate*}
Set~$X_\inv\eqdef s^{-1}(N^\complement)$ and note that~$X_\inv^\complement$ is properly $E$-exceptional. Further define~$\pi \colon x\mapsto \mu_{s(x)}$. For notational simplicity, we relabel~$Z$ as~$Z\setminus N$, so that~\iref{i:t:proof:ErgodicM:1},~\iref{i:t:proof:ErgodicM:2} hold for every~$\zeta\in Z$, and~$X_\inv=s^{-1}(Z)$.
Assertions~\iref{i:t:ErgodicM:2}--\iref{i:t:ErgodicM:3} are standard, e.g.~\cite[Thm.~6.6]{Var01}.
As a consequence of~\iref{i:t:ErgodicM:3}, assertion~\iref{i:t:ErgodicM:1} is precisely the strong consistency of~$\seq{\mu_\zeta}_{\zeta\in Z}$ with~$s$. The \emph{`only if'} part of assertion~\iref{i:t:ErgodicM:4} is straightforward. The \emph{`if'} part is a consequence of the representation in~\iref{i:t:ErgodicM:2}, together with~\iref{i:t:ErgodicM:1}.
\end{proof}

\section{Appendix}\label{s:Appendix}
The theory of direct integrals of Banach spaces is inherently more sophisticated than the corresponding theory for Hilbert spaces. We discuss here an irreducible minimum after~\cite[Ch.s~5-7]{HayMirYve91} and especially~\cite[\S3]{deJRoz17}. For simplicity, we restrict ourselves to the case of $\sigma$-finite (not necessarily complete) indexing spaces~$(Z,\mcZ,\nu)$.

A \emph{decomposition}~$\seq{Z_\alpha}_{\alpha\in \Alpha}$ of~$(Z,\mcZ,\nu)$ is a family of subsets~$Z_\alpha\subset Z$ so that
\begin{align*}
\mcZ=&\ \set{B\subset Z : B\cap Z_\alpha\in \mcZ \text{~for all~$\alpha\in\Alpha$}}
\qquad\text{~and~}\qquad
\\
\nu B=&\ \sum_{\alpha\in\Alpha} \nu (B\cap Z_\alpha)\comma \qquad B\in\mcZ \fstop
\end{align*}

\begin{definition}[Measurable fields, cf.~{\cite[\S6.1, p.~61f.]{HayMirYve91} and~\cite[\S3.1]{deJRoz17}}]\label{d:DirIntB}
Let~$(Z,\mcZ,\nu)$ be a $\sigma$-finite measure space, and~$V$ be a real linear space. A \emph{$\nu$-measurable family of semi-norms on~$V$} is a family~$\ttseq{\norm{\emparg}_\zeta}_{\zeta\in Z}$ so that
\begin{itemize}
\item $\norm{\emparg}_\zeta$ is a semi-norm on~$V$ for every~$\zeta\in Z$;
\item the map~$\zeta\mapsto \norm{v}_\zeta$ is $\nu$-measurable for every~$v\in V$.
\end{itemize}

Letting~$Y_\zeta$ denote the Banach completion of~$V/\ker\norm{\emparg}_\zeta$, we say that a vector field~$u\in \prod_{\zeta\in Z} Y_\zeta$ is \emph{$\nu$-measurable} if, for each~$B\in \mcZ$ with~$\nu B<\infty$, there exists a sequence~$\seq{u_n}_n$ of simple $V$-valued vector fields on~$B$ so that~$\nlim \norm{u_\zeta-u_{n,\zeta}}_\zeta=0$ $\nu$-a.e.\ on~$B$.
A family~$\seq{Y_\zeta}_{\zeta\in Z}$ of Banach spaces~$Y_\zeta$ is a \emph{$\nu$-measurable field of Banach spaces} if there exist
\begin{itemize}
\item a decomposition~$\seq{Z_\alpha}_{\alpha\in\Alpha}$ of~$(Z,\mcZ,\nu)$ consisting of sets of finite $\nu$-measure;
\item a family of real linear spaces~$\seq{Y^\alpha}_{\alpha\in\Alpha}$;
\item for each~$\alpha\in\Alpha$, a $\nu$-measurable family of norms~$\norm{\emparg}_\zeta$ on~$Y^\alpha$,
\end{itemize}
so that, for each~$\alpha\in\Alpha$ and each~$\zeta\in Z_\alpha$, the space~$Y_\zeta$ is the completion of~$(Y^\alpha,\norm{\emparg}_\zeta)$.
Extending the above definition of $\nu$-measurability, we say that~$u\in \prod_{\zeta\in Z} Y_\zeta$ is $\nu$-measurable if (and only if) the restriction of~$u$ to each~$Z_\alpha$ is $\nu$-measurable.

Let~$p\in [1,\infty]$. A $\nu$-measurable vector field~$u$ is called $L^p(\nu)$-integrable if $\norm{u}_p\eqdef \tnorm{ (\zeta\mapsto \norm{u_\zeta}_\zeta) }_{L^p(\nu)}$ is finite. Two $L^p(\nu)$-integrable vector fields~$u$,~$v$ are \emph{$\nu$-equivalent} if~$\norm{u-v}_p=0$.
The space~$Y_p$ of equivalence classes of $L^p(\nu)$-integrable vector fields modulo $\nu$-equivalence, endowed with the non-relabeled quotient norm~$\norm{\emparg}_p$, is a Banach space~\cite[Prop.~3.2]{deJRoz17}, called the \emph{$L^p$-direct integral of~$\zeta\mapsto Y_\zeta$} and denoted by
\begin{align}\label{eq:DirIntP}
Y_p=\tonde{\dint[]{Z} Y_\zeta \diff\nu(\zeta)}_{p} \fstop
\end{align}
\end{definition}

The following is a generalization of Proposition~\ref{p:DIntL^2} to direct integrals of~$L^p$-spaces. Recall~\eqref{eq:Delta}.

\begin{proposition}\label{p:DIntL^p} Let~$(X,\mcX,\mu)$ be $\sigma$-finite standard,~$(Z,\mcZ,\nu)$ be $\sigma$-finite countably generated, and $\seq{\mu_\zeta}_{\zeta\in Z}$ be a separated pseudo-disintegration of~$\mu$ over~$\nu$.
Further let~$\mcA$ be the lattice algebra of all real-valued $\mu$-integrable simple functions on~$(X,\mcX)$. Then, for every~$p\in[1,\infty)$, the map
\begin{equation}\label{eq:Yp}
\begin{aligned}
\iota\colon \class[\mu]{\mcA} & \longrar Y_p\eqdef \tonde{\dint[]{Z} L^p(\mu_\zeta)\diff\nu(\zeta)}_p
\\
\class[\mu]{s} & \longmapsto \class[Y_p]{\delta(s)}
\end{aligned}
\end{equation}
extends to an isomorphism of Banach lattices~$\iota_p\colon L^p(\mu) \rar Y_p$.
\end{proposition}

A proof of the above Proposition~\ref{p:DIntL^p} is quite similar to that of Proposition~\ref{p:DIntL^2}, and therefore it is  omitted. Alternatively, a proof may be adapted from~\cite[\S4.2]{deJRoz17}, having care that:
\begin{itemize}
\item the algebra~$\mcA$ corresponds to the vector lattice~$V$ in~\cite[p.~694]{deJRoz17};
\item the order on $Y_p$ is defined analogously to Remark~\ref{r:Order}, cf.~\cite[p.~694]{deJRoz17};
\item the map~$\iota$ corresponds to the map defined in~\cite[Eqn.~(4.6)]{deJRoz17};
\item the surjectivity of~$\iota_p$ follows as in~\cite[p.~696]{deJRoz17} since it only depends on the disintegration being separated. In the terminology and notation of~\cite{deJRoz17}, this is accounted by the fact that the decomposition~$\beta$ satisfies~\cite[Thm.~4.2(2)]{deJRoz17}.
\end{itemize}

As an obvious corollary to Proposition~\ref{p:DIntL^p}, we obtain that the direct integral of \emph{Hilbert} spaces~$H$ in~\eqref{eq:H} with underlying space of measurable vector fields generated by~$\mcA$ is \emph{identical} to~$Y_2$ as in~\eqref{eq:Yp}. The specification of the underlying space of $\nu$-measurable vector fields is necessary in light of Remark~\ref{r:Caveat}.

\begin{proof}[Proof of Lemma{~\ref{l:CdC}}]
Retain the notation established in \S\ref{ss:AlgInvSets} and in the proof of Theorem~\ref{t:Stone}. 
Firstly, note that~$L^1(\mu)$ is, trivially, an $L^\infty(\mu_0)$-module, and~$\dom{E}$ is an $L^\infty(\mu_0)$-module too, by Definition~\ref{d:Invariant}\iref{i:d:Invariant:4}. As in \S\ref{ss:AlgInvSets}, let~$p$ be the quotient map of~\eqref{eq:Quotient}. For~$u\in L^\infty(\nu)$ denote by~$p^*u\in L^\infty(\mu_0)$ the pullback of~$u$ via~$p$. Setting~$u.\colon f\mapsto p^*u \cdot f$ defines an action of~$L^\infty(\nu)$ on~$L^2(\mu)$ and~$\dom{E}$.
Thus, since the spaces~$(X,\mcX_0,\mu_0)$ and~$(Z,\mcZ,\nu)$ have the same measure algebra by construction of~$Z$, here and in the following we may replace $L^\infty(\mu_0)$-modularity with $L^\infty(\nu)$-modularity.

Let now~$A\in\mcX_0$. Since~$A$ is $E$-invariant, then~$\car_A f\in\dom{E}$ and
\begin{align}\label{eq:l:CdC:1}
E(\car_A f,g)=E(f,\car_A g)=E(\car_A f,\car_A g)\comma \qquad f,g\in\dom{E}
\end{align}
by Definition~\ref{d:Invariant}.
Replacing~$f$ with~$\car_Af$ in~\eqref{eq:CdC}, and applying~\eqref{eq:l:CdC:1} and again~\eqref{eq:CdC} yields
\begin{align}\label{eq:l:CdC:2}
\car_A \cdc(f,g)=\cdc(\car_A f,g)\comma \qquad f,g\in \dom{E}\cap L^\infty(\mu)\comma
\end{align}
which is readily extended to~$f,g\in\dom{E}$ by approximation. Then,~\eqref{eq:l:CdC:2} shows that~$f\mapsto \cdc(f,g)\colon \dom{E}_1\rar L^1(\mu)$ is, for every fixed~$g\in \dom{E}$, a \emph{bounded} $L^\infty(\nu)$-modular operator in the sense of~\cite[\S5.2]{HayMirYve91}. %
By~\emph{Step~1} in the proof of Theorem~\ref{t:Stone}\iref{i:t:Stone:3},~$\dom{E}_1$ is a countably generated direct integral of Banach spaces, thus we may apply~\cite[Thm.~9.1]{HayMirYve91} to obtain, for every fixed~$g\in \dom{E}$, a direct integral decomposition
\begin{align*}
\cdc(\emparg, g)=\dint[]{Z} \cdc_{\zeta,g} \diff\nu(\zeta)\colon \dint{Z} \dom{E_\zeta}_1 \diff\nu(\zeta) \longrar \tonde{\dint{Z} L^1(\mu_\zeta) \diff\nu(\zeta)}_1\cong L^1(\mu)
\end{align*}
for some family of bounded operators~$\cdc_{\zeta,g}\colon \dom{E_\zeta}_1\rar L^1(\mu_\zeta)$.
Let~$\mcC$ be a core for~$(E,\dom{E})$ underlying the construction of the direct integral representation of~$(E,\dom{E})$ as in \emph{Step 1} in the proof of Theorem~\ref{t:Stone}. It follows by symmetry of~$\cdc$ that~$\cdc_{\zeta,g}(f)=\cdc_{\zeta, f}(g)$ for every~$f,g\in\mcC$ and $\nu$-a.e.~$\zeta\in Z$. In particular, the assignment~$g\mapsto \cdc_{\zeta,g}$ is linear on~$\mcC\subset \dom{E_\zeta}$ for $\nu$-a.e.~$\zeta\in Z$. A symmetric bilinear map is then induced on~$\mcC^\otym{2}$ by setting~$\cdc_\zeta\colon (f,g)\mapsto \cdc_{\zeta,g}(f)$.

Thus, finally, it suffices to show~\eqref{eq:CdC} for~$\cdc_\zeta$ and~$(E_\zeta,\dom{E_\zeta})$ for~$\nu$-a.e.~$\zeta\in Z$ with~$f,g,h\in\mcC$, which is readily shown arguing by contradiction, analogously to the proof of the claim in \emph{Step~4} of Theorem~\ref{t:Stone}.
\end{proof}


\begin{thebibliography}{10}
\providecommand{\url}[1]{{#1}}
\providecommand{\urlprefix}{URL }
\expandafter\ifx\csname urlstyle\endcsname\relax
  \providecommand{\doi}[1]{DOI~\discretionary{}{}{}#1}\else
  \providecommand{\doi}{DOI~\discretionary{}{}{}\begingroup
  \urlstyle{rm}\Url}\fi

\bibitem{AlbKonRoe97}
Albeverio, S.A., Kondrat'ev, {\relax Yu}.G., R{\"{o}}ckner, M.: {Ergodicity for
  the Stochastic Dynamics of Quasi-invariant Measures with Applications to
  Gibbs States}.
\newblock {J.\ Funct.\ Anal.} \textbf{149}, 415--469 (1997).
\newblock \doi{10.1006/jfan.1997.3099}

\bibitem{AlbKonRoe98}
Albeverio, S.A., Kondratiev, {\relax Yu}.G., R{\"{o}}ckner, M.: {Analysis and
  Geometry on Configuration Spaces}.
\newblock {J.\ Funct.\ Anal.} \textbf{154}(2), 444--500 (1998).
\newblock \doi{10.1006/jfan.1997.3183}

\bibitem{AlbRoe90}
Albeverio, S.A., R{\"{o}}ckner, M.: {Classical Dirichlet Forms on Topological
  Vector Spaces --- Closability and a Cameron--Martin Formula}.
\newblock {J.\ Funct.\ Anal.} \textbf{88}, 395--436 (1990).
\newblock \doi{10.1016/0022-1236(90)90113-Y}

\bibitem{BiaCar09}
{Bianchini, S.}, {Caravenna, L.}: {On the extremality, uniqueness, and
  optimality of transference plans}.
\newblock {Bull.\ Inst.\ Math.\ Acad.\ Sin.\ (N.S.)} \textbf{4}(4), 353--454
  (2009)

\bibitem{BouHir91}
{Bouleau, N.}, {Hirsch, F.}: {Dirichlet forms and analysis on Wiener space}.
\newblock {De Gruyter} (1991)

\bibitem{Buf14}
{Bufetov, A.~I.}: {Ergodic decomposition for measures quasi-invariant under a
  Borel action of an inductively compact group}.
\newblock {Sb. Math.+} \textbf{205}(2), 192--219 (2014).
\newblock \doi{10.1070/sm2014v205n02abeh004371}

\bibitem{CheMaRoe94}
{Chen, Z.-Q.}, {Ma, Z.-M.}, {R\"ockner, M.}: {Quasi-homeomorphisms of Dirichlet
  forms}.
\newblock {Nagoya Math. J.} \textbf{136}, 1--15 (1994)

\bibitem{ChoGil71}
{Chow, T.~R.}, {Gilfeather, F}: {Functions of Direct Integrals of Operators}.
\newblock {Proc.\ Amer.\ Math.\ Soc.} \textbf{29}(2), 325--330 (1971)

\bibitem{deJRoz17}
{de Jeu, M.}, {Rozendaal, J.}: {Disintegration of positive isometric group
  representations on $\mathrm{L}^p$-spaces}.
\newblock {Positivity} \textbf{21}, 673--710 (2017)

\bibitem{LzDS17+}
Dello~Schiavo, L.: {The Dirichlet--Ferguson Diffusion on the Space of
  Probability Measures over a Closed Riemannian Manifold}.
\newblock {arXiv:1811.11598}  (2018)

\bibitem{LzDSWir21}
Dello~Schiavo, L., Wirth, M.: {Ergodic Decompositions of Dirichlet Forms under
  Order Isomorphisms}.
\newblock {arXiv:2109.00615}  (2021)

\bibitem{Dix81}
{Dixmier, J.}: {Von Neumann Algebras}.
\newblock {North-Holland} (1981)

\bibitem{Ebe96}
{Eberle, A.}: {Girsanov-type transformations of local Dirichlet forms: An
  analytic approach}.
\newblock {Osaka J.~Math.} \textbf{33}(2), 497--531 (1996)

\bibitem{Fre00}
{Fremlin, D.~H.}: {Measure Theory -- Volume I - IV, V Part I \& II}.
\newblock {Torres Fremlin (ed.)} (2000-2008)

\bibitem{FukOshTak11}
{Fukushima, M.}, {Oshima, Y.}, {Takeda, M.}: {Dirichlet forms and symmetric
  Markov processes}, \emph{{De Gruyter Studies in Mathematics}}, vol.~19,
  extended edn.
\newblock {de Gruyter} (2011)

\bibitem{GreSch00}
{Greschonig, G.}, {Schmidt, K.}: {Ergodic decomposition of quasi-invariant
  probability measures}.
\newblock {Colloq.~Math.} \textbf{84}, 495--514 (2000)

\bibitem{HayMirYve91}
Haydon, R., Levy, M., Raynaud, Y.: {Randomly normed spaces}.
\newblock {Travaux en cours}. Hermann (1991)

\bibitem{Kat95}
{Kato, T.}: {Perturbation Theory for Linear Operators}, {Reprint of the 1980
  Edition. Corrected Printing of the Second Edition} edn.
\newblock {Classics in Mathematics}. {Springer-Verlag} (1995)

\bibitem{Kuw98}
{Kuwae, K.}: {Functional Calculus for Dirichlet Forms}.
\newblock {Osaka J.~Math.} \textbf{35}, 683--715 (1998)

\bibitem{Kuw11}
{Kuwae, K.}: {Invariant sets and ergodic decomposition of local semi-Dirichlet
  forms}.
\newblock {Forum Math.} \textbf{23}, 1259--1279 (2011).
\newblock \doi{10.1515/FORM.2011.046}

\bibitem{Lan75}
{Lance, C.}: {Direct Integrals of Left Hilbert Algebras}.
\newblock {Math.\ Ann.} \textbf{216}, 11--28 (1975)

\bibitem{LenSchWir18}
{Lenz, D.}, {Schmidt, M.}, {Wirth, M.}: {Geometric Properties of Dirichlet
  Forms under Order Isomorphism}.
\newblock {arXiv:1801.08326}  (2018)

\bibitem{MaRoe92}
Ma, Z.M., R{\"{o}}ckner, M.: Introduction to the Theory of (Non-Symmetric)
  Dirichlet Forms.
\newblock {Graduate Studies in Mathematics}. Springer (1992)

\bibitem{Neu49}
von Neumann, J.: {On Rings of Operators. Reduction Theory}.
\newblock {Ann.\ Math.} \textbf{50}(2), 401 (1949).
\newblock \doi{10.2307/1969463}.
\newblock \urlprefix\url{https://doi.org/10.2307%2F1969463}

\bibitem{Oku92}
{{\^O}kura, H.}: {On the Invariant Sets and the Irreducibility of Symmetric
  Markov Processes}.
\newblock In: {Shiryaev, A.~N.}, {Korolyuk, V.~S.}, {Watanabe, S.}, {Fukushima,
  M.} (eds.) {Probability Theory and Mathematical Statistics -- Proceedings of
  the 6th USSR--Japan Symposium}, pp. 238--247. {World Scientific} (1992).
\newblock \doi{10.1142/1780}

\bibitem{Sch95}
{Schmuland, B.}: {On the local property for positivity preserving coercive
  forms}.
\newblock In: {Ma, Z.-M.}, {R{\"{o}}ckner, M.}, {Yan, J.~A.} (eds.) {Dirichlet
  Forms and Stochastic Processes: Proceedings of the International Conference
  held in Beijing, China, October 25-31, 1993}, {De Gruyter Proceedings in
  Mathematics}, pp. 345--354 (1995).
\newblock \doi{10.1515/9783110880052.345}

\bibitem{Sch73}
{Schwartz, L.}: {Radon measures on arbitrary topological spaces and cylindrical
  measures}.
\newblock Oxford University Press (1973)

\bibitem{tElRobSikZhu06}
{ter Elst, A. F. M.}, {Robinson, D. W.}, {Sikora, A.}, {Zhu, Y.}: {Dirichlet
  Forms and Degenerate Elliptic Operators}.
\newblock In: {Koelink, E.}, {van Neerven, J.}, {de Pagter, B.}, {Sweers, G.},
  {Luger, A.}, {Woracek, H.} (eds.) {Partial Differential Equations and
  Functional Analysis: The Philippe Cl{\'e}ment Festschrift}, pp. 73--95.
  {Birkh{\"a}user} (2006)

\bibitem{Tom80}
{Tomisaki, M.}: {Superposition of diffusion processes}.
\newblock {J.~Math.\ Soc.\ Japan} \textbf{32}(4), 671--696 (1980)

\bibitem{Var01}
{Varadhan, S.~R.~S.}: {Probability Theory}, \emph{{Courant Lecture Notes}},
  vol.~7.
\newblock {Amer.\ Math.\ Soc.} (2001)

\end{thebibliography}
\end{document}